\documentclass[11pt]{article}
\usepackage{amssymb,latexsym,amsmath, amsthm}
\usepackage{amsfonts,amssymb, latexsym, amsmath, amsthm,mathrsfs, verbatim,geometry}

\usepackage{tikz}
\usepackage{amsfonts, graphics}

\usetikzlibrary{calc}
\usetikzlibrary{arrows.meta}

\usepackage{tikz}
\usetikzlibrary{shapes,backgrounds,shadows,arrows} 
\usetikzlibrary{arrows.meta,calc,decorations.markings,math}
\definecolor{ubtgreen}{RGB}{3,138,94}

\usepackage{mathtools}
\usepackage{hyperref}
\usepackage{csquotes} %enquote
\usepackage{enumitem,moreenum} %labels in enums
\usepackage{float}
\usepackage{dsfont,bbm} 
\numberwithin{equation}{section}
\def\R{\mathbb R}

\def\N{\mathbb N}

\hypersetup{
    colorlinks,
    citecolor=red,
    filecolor=green,
    linkcolor=blue,
    urlcolor=black
}

\def\R{\mathbb{R}}

\def\N{\mathbb N}
\def\Z{\mathbb Z}

\def\A{\mathscr{A}}

\def\D{\mathcal{D}}
\def\be{\begin{equation}}
\def\ee{\end{equation}}
\def\bea{\begin{eqnarray}}
\def\eea{\end{eqnarray}}
\def\beas{\begin{eqnarray*}}
\def\eeas{\end{eqnarray*}}

\def\pa{\partial }

\def\P{{\mathcal P}}

\def\A{{\mathcal A}}
\def\D{{\mathcal D}}
\def\B{{\mathcal B}}

\def\pa{\partial}

\def\al{\alpha}
\def\de{\delta}

\def\lv{\left\vert}
\def\rv{\right\vert}
\def\Emin{E_{\mathrm{min}}}
\def\Lmax{L_{\mathrm{max}}}
\def\k{\varepsilon}

\def\Rmax{R_{\mathrm{max}}}
\def\Rmin{R_{\mathrm{min}}}

\def\Q{\mathcal{Q}}

\def\al{\alpha}
\def\de{\delta}

\def\dif{\textup{d}}
\def\EminL{E^L_{\min}}
\def\A{\mathcal A}
\def\B{\mathcal B}
\def\P{\mathcal P}
\def\Q{\mathcal Q}

\newcommand{\fhat}{\widehat{f}}
\newcommand{\ghat}{\widehat{g}}

\let\th\relax
\newcommand{\th}{\theta}
\newcommand{\eps}{\varepsilon}
\newcommand\om{\omega}

\newcommand{\beq}{\begin{equation}}
\newcommand{\eeq}{\end{equation}}
\newcommand{\beqs}{\begin{equation*}}
\newcommand{\eeqs}{\end{equation*}}
\newcommand{\beqa}{\begin{equation}\begin{aligned}}
\newcommand{\eeqa}{\end{aligned}\end{equation}}
\newcommand{\beqas}{\begin{equation*}\begin{aligned}}
\newcommand{\eeqas}{\end{aligned}\end{equation*}}

\let\originalleft\left
\let\originalright\right
\renewcommand{\left}{\mathopen{}\mathclose\bgroup\originalleft}
\renewcommand{\right}{\aftergroup\egroup\originalright}

\newcommand*\diff{\mathop{}\!\mathrm{d}} %Write the d in integrals in mathrm\\

\newtheorem{theorem}{Theorem}[section]

\newtheorem{prop}[theorem]{Proposition}

\newtheorem{corollary}[theorem]{Corollary}

\newtheorem{lemma}[theorem]{Lemma}

\newtheorem{remark}[theorem]{Remark}

\def\be{\begin{equation}}
\def\ee{\end{equation}}
\def\Tmax{T_{\mathrm{max}}}
\def\Tmin{T_{\mathrm{min}}}

\def\Om{\Omega}

\def\dif{\textup{d}}

\def\bcr{\begin{color}{red}}
\def\bcb{\begin{color}{blue}}
\def\bcc{\begin{color}{violet}}
\def\ec{\end{color}}

\begin{document}

\title{\vspace*{-2.2em}Quantitative phase mixing for Hamiltonians with trapping}
\author{M.~Had\v zi\'c\thanks{University College London, UK. Email: m.hadzic@ucl.ac.uk}, G.~Rein\thanks{University of Bayreuth, Germany. Email: gerhard.rein@uni-bayreuth.de}, M.~Schrecker\thanks{University of Bath, UK. Email: mris21@bath.ac.uk}, C.~Straub\thanks{University of Bayreuth, Germany. Email: christopher.straub@uni-bayreuth.de}}

\maketitle

\vspace*{-1.2em}
\begin{abstract}
We prove quantitative decay estimates of macroscopic quantities generated by the solutions to linear transport equations driven by a general family of Hamiltonians. The associated particle trajectories are all trapped in a compact region of phase-space and feature a non-degenerate elliptic stagnation point. The analysis covers a large class of Hamiltonians generated by the radially symmetric compactly supported equilibria of the gravitational Vlasov-Poisson system. Working in radial symmetry, our analysis features both the 1+2-dimensional case and the harder 1+1-dimensional case, where all the particles have the same value of the modulus of angular momentum. The latter case is also of importance in both the plasma physics case and two dimensional incompressible fluid flows.
\end{abstract}

\tableofcontents

\section{Introduction}\label{S:INTRO}

We derive quantitative decay-in-time estimates for the potential generated by the phase-space density~$f$ satisfying the transport equation
\begin{align}\label{E:PT0L}
\pa_t f(t,r,w,L) + \D f(t,r,w,L) = 0,
\end{align} 
where the anti-symmetric transport operator $\D$ takes the form
\begin{align}\label{E:DDEFL}
\D \coloneqq w\pa_r - \Psi_L'(r) \pa_w.
\end{align}
Here we shall use potentials  $(0,\infty)\ni r\mapsto \Psi_L(r)\in\R$ with parameter $L>0$  generated by
compactly supported spherical steady states of the gravitational Vlasov-Poisson system,
the characteristic orbits of which are therefore
bounded in both space and velocity.

We shall also be interested in the harder, restricted case
\begin{align}\label{E:PT0}
\pa_t f(t,r,w) + \D f(t,r,w) = 0,
\end{align} 
where we drop the dependence on the angular momentum $L$
and the anti-symmetric transport operator $\D$ takes the form
\begin{align}\label{E:DDEF}
\D \coloneqq w\pa_r - \Psi'(r) \pa_w.
\end{align}
We refer to~\eqref{E:DDEF} as the 1+1-dimensional case.
The problem~\eqref{E:PT0}--\eqref{E:DDEF} arises in the study of several Hamiltonian transport problems, including the gravitational Vlasov-Poisson system when all the stars have a fixed  modulus of  angular momentum. Related questions appear in the Vlasov-Poisson system in the electrostatic case describing space-inhomogeneous plasmas and the two-dimensional incompressible Euler equation. 
Our results apply to a general class of potentials $\Psi_L$ and $\Psi$ that satisfy certain qualitative properties, which will be listed below.
%%%%%%%%%%%%%%%%%%%%%%%%%%%
%%%%%%%%%%%%%%%%%%%%%%%%%%%

\subsection{Motivation: the gravitational Vlasov-Poisson system}

%%%%%%%%%%%%%%%%%%%%%%%%%%%
%%%%%%%%%%%%%%%%%%%%%%%%%%%

The gravitational Vlasov-Poisson system is a fundamental Newtonian model of galaxy dynamics~\cite{BiTr}. The basic unknown is
the phase space density
$f \colon \mathbb R\times\mathbb R^3 \times\mathbb R^3\to[0,\infty)$, i.e.,
$f=f(t,x,v)\geq 0$ is a function of time $t\in \mathbb R$,
position $x\in \mathbb R^3$, and velocity (or momentum) $v\in \mathbb R^3$.
We restrict ourselves to radial symmetry where $f$ can
be expressed as a function of $t$ and
\begin{equation*}
r=|x|, \ w= \frac{x\cdot v}{r}, \ L = |x\times v|^2,
\end{equation*}
and the Vlasov-Poisson system takes the form
\begin{equation}
  \partial_tf + w\pa_ rf-\Psi_L'(t,r) \pa_wf = 0, 
  \label{E:FEQN}
\end{equation}
\begin{equation}
  \left(\pa_{rr}+\frac2r\pa_r\right) U(t,r) = 4\pi\rho(t,r),  \quad
  \rho(t,r) =\frac{\pi }{r^2}\int_{\R}\int_{(0,\infty)} f(t,r,w,L)
  \diff L \diff w. \label{E:POISSONEQN}
\end{equation}
In the above $\Psi_L(t,r)$ is the so-called effective potential which takes the form 
\begin{align}\label{E:PSILDEF}
 \Psi_L(t,r) := U_{\text{ext}}(r) + U(t,r)+\frac{L}{2r^2},
\end{align}
where $U\colon\mathbb R\times[0,\infty)\to\mathbb R$ is the gravitational potential induced by the distribution of matter in the galaxy
and $U_{\text{ext}}$ a time-independent external gravitational potential. The most interesting cases are when $U_{\text{ext}}$ is identically zero
or when $U_{\text{ext}}$ is the potential of a point mass $M>0$ located at the origin $r=0$:
\begin{align}\label{E:DELTA}
U_{\text{ext}}(x) = U_{\text{ext}}(r) = - \frac Mr, \qquad |x|=r>0.
\end{align}
The choice~\eqref{E:DELTA} gives a Newtonian model of a galaxy featuring a large gravitational source (e.g. black hole) in the centre.
Finally, to ensure that we are describing an 
isolated galaxy, we impose the boundary condition
\begin{equation}\label{E:ISOLATED}
\lim_{r\to\infty} U(t,r) = 0.
\end{equation}

\subsubsection{Steady states}

The system~\eqref{E:FEQN}--\eqref{E:PSILDEF}, \eqref{E:ISOLATED} possesses
a plethora of steady states whose stability properties are of great interest 
in the astrophysics literature~\cite{BiTr}. Due to the assumption of
spherical symmetry, 
the phase space density of these steady states is always of the form
\begin{align}\label{E:POLYL}
\mathring f(r,w,L) = \varphi(E,L), 
\ \ E(r,w,L)=\frac12w^2+\frac{L}{2r^2}+U_0(r)+U_{\text{ext}}(r),
\end{align}
where $E(r,w,L)$ is the local particle energy and $U_0$ the self-induced gravitational potential~\cite{BaFaHo86}.
Here the function $\varphi$ can be thought of as a microscopic equation of state. In the so-called {\em regular} case $U_{\text{ext}}=0$, it is well-known~\cite{BaFaHo86,RaRe13} that finite-mass compactly supported solutions of the form~\eqref{E:POLYL} exist 
if, for example,
\begin{align}
& \varphi(E,L) = (E_0-E)^k_+(L-L_0)_+^\ell \qquad \text{ (polytropic galaxy) },\label{E:POLY} \\
& 0<k<3\ell + \frac72, \ \ \ell>-\frac12,
\end{align}
where $E_0<0$ is a suitable cut-off energy whose value depends on the solution itself and $L_0\ge0$ a given cut-off angular momentum. 
We notice that the regularity of the steady states~\eqref{E:POLY} at the vacuum (corresponding to $E=E_0$ or $L=L_0$) is limited by the size of the polytropic exponents $k,\ell$ and this limited regularity 
is a generic feature of the problem.

On the other hand, in the point mass case $U_{\text{ext}}=-\frac Mr$ with $M>0$, the set of steady states is even larger with less stringent regularity constraints at the vacuum boundary.
There exist finite mass steady states supported on a compact shell around the galactic center of the form
\begin{align}
&\varphi(E) = \varepsilon (E_0-E)^k_+(L-L_0)_+^\ell \ \ \text{ (polytropic galaxy with a central delta mass)}, \label{E:POLY2} \\ 
 &k>0, \ \ell>-\frac12, \ L_0>0,
\end{align}
cf.~\cite[Sc.~6.2]{St23}, where the parameter $\eps>0$ controls the size of the steady state compared to the point mass. 
The interior of the phase space support of such steady states is given by
\begin{equation}\label{E:OMDEFL}
\Om \coloneqq \left\{(r,w,L)\mid  E(r,w,L)< E_0, \ L> L_0\right\}.
\end{equation}
The closure of the spatial support corresponds to an interval $[\Rmin,\Rmax]\subset(0,\infty)$.
This implies that these steady states feature a vacuum region around the centre and are hence referred to as {\em shells}.
 In the regular case~\eqref{E:POLY} we also restrict ourselves here to shell-shaped steady states by choosing $L_0>0$, cf.~\cite[Sc.~2.2]{St23}. 
This removes some technical issues arising at the spatial origin $r=0$, but it is straight-forward to relax this assumption.

{\bf Restricted (1+1-d) case.}
We consider now a symmetry reduction wherein all stars in the galaxy are constrained to have 
the same squared modulus of the angular momentum $L=\bar L>0$. This is a dynamically consistent constraint
as the function $L(x,v)$ is conserved along the characteristics.
The quantity $L>0$ becomes a fixed parameter 
and the phase space coordinates are reduced to $(r,w)$.
In this setting, it is straight-forward to show that analogues of the steady states of the form~\eqref{E:POLY} and~\eqref{E:POLY2} also exist provided that $k>\frac12$ due to regularity reasons~\cite{HRSS2023}. 
Specifically, in the case with the point mass potential $U_{\text{ext}}(r)=-\frac M r$, we  consider steady states of the form
\begin{equation}\label{E:POLY3}
\mathring f(r,w) = \varphi(E) = \varepsilon (E_0-E)_+^k, \quad k>\frac12,
\end{equation}
with the interior of its compact support given by
\begin{equation}\label{E:OMDEF}
\Om \coloneqq \left\{(r,w)\mid  E(r,w)< E_0\right\},
\end{equation}
where $E=E(r,w)$ is defined similarly to~\eqref{E:POLYL}.
Just as in the above, the spatial extent of~$\bar{\Om}$ is given by some interval $[\Rmin,\Rmax]$, where $0<\Rmin<\Rmax$, cf.~\cite{HRSS2023}.
We warn the reader that we largely keep the same notation for both the radially symmetric and the 1+1-dimensional case for analogous quantities; no confusion should arise in the process. 

\subsubsection{Stability theory}
A classical condition that implies the spectral stability of the above class of steady states is referred to as the Antonov stability condition~\cite{An1961,BiTr,DoFeBa,KS} and, modulo suitable assumptions on the regularity of $\varphi$, it reads
\begin{equation}\label{E:ASC}
  \varphi'\coloneqq\pa_E\varphi<0 \ \ \text{inside the support of the steady
    state}.
\end{equation}
Examples~\eqref{E:POLY}, \eqref{E:POLY2}, and~\eqref{E:POLY3} clearly satisfy~\eqref{E:ASC}. Due to the energy subcritical nature of the problem, solutions satisfying~\eqref{E:ASC}
are in fact known to be nonlinearly orbitally stable~\cite{GuRe2007,LeMeRa11,LeMeRa12}. 

By contrast, almost nothing is known about the asymptotic stability of the above family of steady states. To understand the problem, we linearise the Vlasov-Poisson system around a steady state~\eqref{E:POLYL} to obtain
\begin{equation}\label{E:FULLLIN}
\pa_t f + \D f + \B f = 0, 
\end{equation}
where $\D$ is given by~\eqref{E:DDEFL} and
%we recall~\eqref{E:PSILDEF}.
the potential~$U$ appearing in~\eqref{E:PSILDEF} is the steady state potential~$U_0$. 
The operator $\B$ is nonlocal, formally of lower order than $\D$, and given by
\begin{equation*}
\B f = |\varphi'(E,L)|\D U_f,
\end{equation*}
where
$U_f$ is the solution of the Poisson equation $\Delta U_f = 4\pi \int f \diff v$ satisfying the boundary condition~\eqref{E:ISOLATED}.
The process of equilibration of macroscopic quantities generated by the solutions of~\eqref{E:FULLLIN} is referred to as {\em gravitational relaxation}~\cite{BiTr,LB1962,LB1967}.
Numerical simulations show that macroscopic damping at the full non-linear level is accurately described by the this effect~\cite{St24}.
 
The natural Hilbert space in the stability analysis is given by 
\begin{equation*}
H\coloneqq\left\{f\colon\Om\to\mathbb R\, \Big| \, \int_{\Om}\frac{f(r,w,L)^2}{|\varphi'(E,L)|}\diff (r,w,L)<\infty\right\},
\end{equation*}
or
\begin{equation*}
\tilde H\coloneqq\left\{f\colon\Om\to\mathbb R\, \Big| \, \int_{\tilde\Om}\frac{f(r,w)^2}{|\varphi'(E)|}\diff (r,w)<\infty\right\}
\end{equation*}
in the 1+1-dimensional case.
The transport operator~$\D$ can be realised as an anti-selfadjoint operator
on~$H$ (or $\tilde H$)~\cite{ReSt20}, and $\B$ is
relatively compact with respect to $\D$ (see, e.g.,~\cite{HRS2022}).
Therefore the essential spectra of $\D$ and $\D+\B$ are the same so that in particular the decay induced by $\D$ accounts for dispersive effects induced by the continuous part of the spectrum. 
We emphasise that the pure transport problem~\eqref{E:PT0L}--\eqref{E:DDEFL} however does {\em not} correspond to the linearised Vlasov-Poisson dynamics except in the special case of a vanishing steady state $\mathring f\equiv0$ with an external potential, e.g., $U_{\text{ext}}=-\frac Mr$.

Moreover, as proven in~\cite{HRS2022,HRSS2023}, there are examples of steady states of the form~\eqref{E:POLY} where the linearised operator $\D+\B$ possesses purely oscillatory eigenmodes (and thus there is no decay for generic initial data), even though the solutions of the pure transport problem~\eqref{E:PT0} will generically decay. 
We refer the reader to~\cite{HRSS2023,St23,St24} for a detailed discussion of the role of the steady state regularity (as measured by the polytropic exponents $k,\ell$ in~\eqref{E:POLY} and~\eqref{E:POLY2} and the value of~$k$ in~\eqref{E:POLY3}) in discerning oscillatory vs.\ damped behaviour.

\begin{remark}\label{R:IC}
Since $(E,L)\mapsto |\varphi'(E,L)|$ is supported on the support $\overline\Om$ of the steady state~$\mathring{f}$, it is natural to consider perturbations of
the form $f = |\varphi'(E,L)|g$. This recasts~\eqref{E:FULLLIN} into the form
%\begin{align}
$\pa_t g + \D g + \D U_{|\varphi'(E,L)|g} = 0$.
%\end{align}
We now see that the linearised dynamics~\eqref{E:FULLLIN} formally preserve the 
property $f(t,\cdot) = |\varphi'(E,L)|g(t,\cdot)$ if imposed initially as $f(0,\cdot) = |\varphi'(E,L)|g(0,\cdot)$.
\end{remark}

%%%%%%%%%%%%%%%%%%%%%%%%%%%
%%%%%%%%%%%%%%%%%%%%%%%%%%%

\subsection{Main results}

%%%%%%%%%%%%%%%%%%%%%%%%%%%
%%%%%%%%%%%%%%%%%%%%%%%%%%%

We now return to the main objectives of this note -- the decay properties of macroscopic quantities generated 
by the flows~\eqref{E:PT0L} and~\eqref{E:PT0}. 
We emphasise again that this is not the true linearised Vlasov-Poisson dynamics around the equilibrium, but instead its leading order part. 
The analysis of~\eqref{E:PT0L} and~\eqref{E:PT0} is both conceptually and technically easier than the analysis of the full linearisation~\eqref{E:FULLLIN}.
Moreover, it is not clear to what extent the linearised dynamics is expected to qualitatively look like the one generated by~\eqref{E:PT0L} and~\eqref{E:PT0}, cf.~\cite[Obs.~4.6]{St24}.

The two canonical macroscopic quantities in the study of gravitational relaxation are the macroscopic density
\begin{equation}\label{E:DENSITYDEF}
  \rho_f(t,r) \coloneqq
  \frac{\pi }{r^2}\int_\R \int_{(0,\infty)}f(t,r,w,L)  \diff L \diff w
\end{equation}
and the gravitational potential $U_f$ generated by the density distribution $\rho_f$. It is defined as the unique solution of the Poisson equation
\begin{equation}\label{E:POISSON}
\Delta U_f = 4\pi \rho_f, \quad \lim_{r\to\infty} U_f(t,r) = 0.
\end{equation}
In the 1+1-dimensional setting~\eqref{E:PT0} the macroscopic density is given by the formula
\begin{equation}\label{E:DENSITYDEF1}
\rho_f(t,r) \coloneqq \frac{\pi }{r^2}\int_\R f(t,r,w) \diff w.
\end{equation}
The goal of this note is to establish quantitative decay rates for these quantities under the flow~\eqref{E:PT0L} or~\eqref{E:PT0}.

\subsubsection{The 1+1-dimensional case}

%%%%%%%%%%%%%%%%%%%%%%%%%%%
%%%%%%%%%%%%%%%%%%%%%%%%%%%

In this case, the transport operator $\D$ has a nontrivial kernel~\cite{HRSS2023,ReSt20} given by
\begin{equation*}
\text{ker}(\D) = \left\{f(r,w) = \phi(E) \ \big| \ \|f\|_{H}<\infty\right\},
\end{equation*} 
where $E=E(r,w)$ is given by~\eqref{E:POLYL}; recall that $L=\bar L>0$ is a fixed constant.
It follows that equation~\eqref{E:PT0} possesses infinitely many steady states that are simply given as any function of $E$ (in the right Hilbert space). Therefore, to meaningfully address the decay of solutions
to~\eqref{E:PT0} we must either mod out the kernel by considering data in its orthogonal complement 
$\text{ker}(\D)^\perp$ in $H$,
or consider the decay of the time-differentiated spatial density $\pa_t\rho_f$ and gravitational potential $\pa_t U_f$. 

\begin{figure}[h!]
	\begin{tikzpicture}
		
		\draw[-{>},black] (-.5,0.) -- (12.6,0.) node[right] {$r$};
		\draw[-{>},black] (0.,-2.8) -- (0.,.5) node[left] {$\Psi(r)$};
		%	\node at (-.1,.13) {$0$};
		
		%\draw[domain=0:2.5, smooth, variable=\x, gray,samples=100,dashed,thick] plot ({5*\x}, {-4./(1+\x*\x)});
		%\node[gray] at (2.2,-3.8) {$U_0$};

		\draw[ubtgreen, densely dotted,thick] (2.6833,-2.4112) -- (2.6833,.1) node[above,ubtgreen] {$r_\ast$};
		\draw[ubtgreen, densely dotted,thick] (2.6833,-2.4112) -- (-.1,-2.4112) node[left,ubtgreen] {$\Emin$};
		
		\draw[purple,densely dotted,thick] (6.00165,-1.5) -- (-.1,-1.5) node[left,purple] {$E$};
		\draw[purple, densely dotted, thick] (1.521035,-1.5) -- (1.521035,.1) node[above,purple] {$r_-(E)$};
		\draw[purple, densely dotted, thick] (6.00165,-1.5) -- (6.00165,.1) node[above,purple] {$r_+(E)$};
		
		\draw[domain=0.225:2.5, smooth, variable=\x, black, samples=100,thick] plot ({5*\x}, {-4./(1+\x*\x) + .2/(\x*\x)});	
	\end{tikzpicture} 
	\caption{The qualitative shape of the effective potential~$\Psi$ in the 1+1-dimensional case.}
	\label{F:PW}
\end{figure}
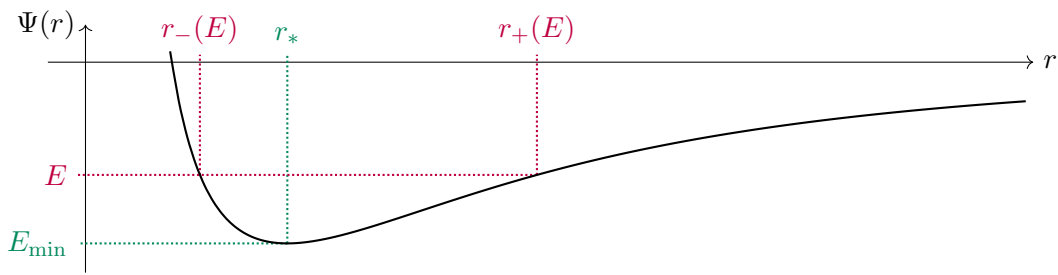

{\bf Assumptions on $\Psi$.}
The steady states of the form~\eqref{E:POLY3} lead to an effective potential 
\be\label{E:PSIDEF}
\Psi(r) = - \frac Mr + U_0(r) + \frac{\bar L}{2r^2}, 
\ee
cf.~\cite{HRSS2023}.
For concreteness, we assume $\Psi$ in~\eqref{E:PT0} to be of the form~\eqref{E:PSIDEF}. 
The graph of~$\Psi$ is qualitatively depicted in Figure~\ref{F:PW}. 
In particular, $\Psi$ has a unique global minimum at a radius $r_*\in(\Rmin,\Rmax)$~\cite{HRSS2023} such that
\begin{equation}
\Emin\coloneqq\Psi(r_*)<0,\qquad\Psi'(r_*)=0,\qquad \Psi''(r_*)>0.  \label{E:ELLIPTIC}
\end{equation}
We additionally assume that
\begin{equation}\label{E:PSIREG}
\Psi \in C^3(0,\infty)\cap C^\infty(\Rmin,\Rmax).
\end{equation}
This is, for example, true when the gravitational potential $U_0$ is generated by steady states of the form~\eqref{E:POLY3} with $k>\frac12$; see~\cite[Sc.~3]{HRSS2023} for similar arguments.
Below, we will impose a few further assumptions on~$\Psi$ via the characteristic flow it generates.

\begin{remark}\label{R:FEATURES}
Assuming the specific form~\eqref{E:PSIDEF} of the effective potential is not essential. 
Modulo regularity assumptions on $\Psi$, our fundamental requirements are 1)~conditions~\eqref{E:ELLIPTIC} and 2)~the requirement that all particle orbits generated by the Hamiltonian $E(r,w)$ with negative energy values are bounded in phase-space and have a strictly monotone period function viewed as a function of~$E$, cf.\ below. 
Assumptions~\eqref{E:ELLIPTIC} guarantee that the characteristic dynamical system features a single nondegenerate elliptic stagnation point. 
Potentials with these features also naturally appear in fluid mechanics~\cite{BCZV19,BrCZMa2022,Na,RhYo1983} and plasma physics (see~\cite{GuLi2017} and the references therein). We point to more references in the work~\cite{BrCZMa2022} where enhanced dissipation in the presence of diffusion and elliptic points is studied.
\end{remark}

%%%%%%%%%%%%%%%%%%%%%%%%%%%%
%%%%%%%%%%%%%%%%%%%%%%%%%%%%
%%%%%%%%%%%%%%%%%%%%%%%%%%%%%%%%%

%%%%%%%%%%%%%%%%%%%%%%%%%%%%
%%%%%%%%%%%%%%%%%%%%%%%%%%%%

\begin{theorem}\label{T:MAIN}
Let $k>\frac12$ and, for any  $g_0\in W^{3,\infty}(\bar{\Omega})$, consider the initial data of the form
\begin{align}
f_0(r,w) = |\varphi'(E(r,w))| g_0(r,w), \label{E:INITIAL}
\end{align}
where $\varphi$ is given by~\eqref{E:POLY3}.
Let the operator $\D$ be given by~\eqref{E:DDEF} with $r\mapsto \Psi(r)$ as described above. 
Let $f(t,\cdot,\cdot)$ be the solution to~\eqref{E:PT0} with initial data $f(0,r,w)=f_0(r,w)$. 
\begin{enumerate}
\item[(a)]There
exists a constant $C>0$ such that, for all $t\in\R_+$, 
\begin{align}
|\pa_t U_f(t,r)| &\le C\|g_0\|_{W^{3,\infty}}\,\frac {1+\lv \log|r-r_\ast|\rv}{1+t^{\min\{2,k\}}},\qquad r\neq r_\ast,\label{E:PARTIALTUFDECAYINTRO} \\
\|\pa_t\pa_r U_f(t,\cdot)\|_{L^\infty_r} & \le \frac {C\|g_0\|_{W^{2,\infty}}}{1+t^{\min\{\frac32,k\}}}, \label{E:GRADPOTENTIALDECAYINTRO}\\
\| \pa_t\rho_f(t,\cdot) \|_{L^\infty_r} &\le \frac {C\|g_0\|_{W^{1,\infty}}}{1+t^{\min\{\frac12,k-\frac12\}}}. \label{E:DENSITYDECAYINTRO}
\end{align}
If additionally $g_0\in W^{4,\infty}(\Om)$ and $k\geq\frac52$, then there exists a constant $C>0$ such that 
\begin{equation}
|\pa_t U_f(t,r)| \le C\|g_0\|_{W^{4,\infty}}\,\frac{(r-r_\ast)^{-1}}{1+|t|^{\frac52}}, \qquad r_\ast<r\le \Rmax.\label{E:PARTIALTUFDECAYEXTERIORINTRO}
\end{equation} 
\item[(b)]
Moreover, if $f_0\in \text{{\em ker}}\,(\D)^\perp$, then for all $t\in\R_+$,
\begin{align}
\| \rho_f(t,\cdot) \|_{L^\infty_r} \le \frac{C\|g_0\|_{W^{1,\infty}}}{1+t^{\min\{\frac12,k-\frac12\}}}, \label{E:DECAYJUSTRHOMAIN} \\
\| \pa_rU_f(t,\cdot) \|_{L^\infty_r} \le \frac{C\|g_0\|_{W^{2,\infty}}}{1+t^{\min\{\frac32,k\}}}. \label{E:PARTIALRUFDECAYEXTRAMAIN}
\end{align}
\end{enumerate}
\end{theorem}

%%%%%%%%%%%%%%%%%%%%%%%%%%%%
%%%%%%%%%%%%%%%%%%%%%%%%%%%%

The theorem is proven in Section~\ref{S:DECAY}.
The proof of~\eqref{E:DENSITYDECAYINTRO} and~\eqref{E:DECAYJUSTRHOMAIN} is the content of Proposition~\ref{P:RHOTDECAY}. The proofs of~\eqref{E:PARTIALTUFDECAYINTRO}
and~\eqref{E:PARTIALTUFDECAYEXTERIORINTRO} are contained in Proposition~\ref{P:DECAY3}. The proofs of~\eqref{E:PARTIALRUFDECAYEXTRAMAIN} and~\eqref{E:GRADPOTENTIALDECAYINTRO}
are contained in Proposition~\ref{P:PARTIALRUFDECAY} and Corollary~\ref{C:PARTIALRUFDECAY}, respectively.

\begin{remark}[$(1+t)^{-\frac32}$-decay]
The theorem shows that for sufficiently smooth data and the potential $\Psi$, the force field decays like $t^{-\frac32}$. The improvement from $t^{-1}$ to the integrable $t^{-\frac32}$ rate is delicate and relies on a precise description of certain derivative singularities of the characteristic flow near the turning points $w=0$. The rate $(1+t)^{-\frac32}$ is likely optimal.
\end{remark}

\begin{remark}
The estimates~\eqref{E:PARTIALTUFDECAYINTRO} and~\eqref{E:PARTIALTUFDECAYEXTERIORINTRO} feature constants that degenerate as the spatial radius
approaches the stationary radius $r_\ast$. This is a quantitative description of the trapping effect coming from the elliptic point $E=\Emin$ of the characteristic flow.
Interestingly, our proof shows that in the outer
region  $\{r>r_\ast\}$ of the galaxy/plasma 
the gravitational/electrostatic potential enjoys a  decay
which is enhanced by a power of $t^{\frac12}$ compared to the inner
region $\{r<r_\ast\}$, at the expense of constants that blow up 
at $r=r_\ast$.
\end{remark}

\begin{remark}
Assumption~\eqref{E:INITIAL} encodes the vanishing (or the absence thereof) of $f_0,\pa_E f_0,\pa_{EE}f_0$ at the vacuum boundary $\pa\Omega$. This assumption is physically natural as explained in Remark~\ref{R:IC}. We can replace the assumption~\eqref{E:INITIAL} by the requirement
\[
\mathrm{supp}\,(f_0)\Subset\Omega, \ \ f_0\in W^{3,\infty}(\Omega),
\] 
which in turn implies infinite-order vanishing of the initial data near the vacuum boundary.
\end{remark}

%%%%%%%%%%%%%%%%%%%%%%%%%%%%%%%
%%%%%%%%%%%%%%%%%%%%%%%%%%%%%%%

\subsubsection{The general radially symmetric case}

%%%%%%%%%%%%%%%%%%%%%%%%%%%%%%%
%%%%%%%%%%%%%%%%%%%%%%%%%%%%%%%

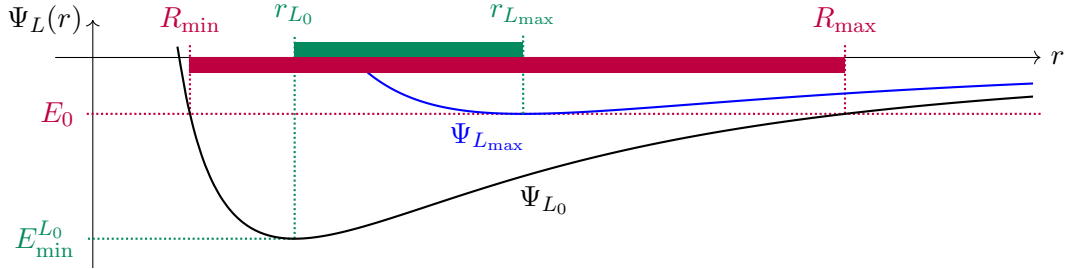
\begin{figure}[h!]
	\begin{tikzpicture}
		\draw[-{>},black] (-.5,0.) -- (12.6,0.) node[right] {$r$};
		\draw[-{>},black] (0.,-2.8) -- (0.,.5) node[left] {$\Psi_L(r)$};
		
		\draw[ubtgreen, densely dotted,thick] (2.6833,-2.4112) -- (2.6833,.3) node[above,ubtgreen] {$r_{L_0}$};
		\draw[ubtgreen, densely dotted,thick] (2.6833,-2.4112) -- (-.1,-2.4112) node[left,ubtgreen] {$\Emin^{L_0}$};
		
		\draw[purple,densely dotted,thick] (12.6,-.75) -- (-.1,-.75) node[left,purple] {$E_0$};
		\draw[purple, densely dotted, thick] (1.29099,-.75) -- (1.29099,.2) node[above,purple] {$\Rmin$};
		\draw[purple, densely dotted, thick] (10,-.75) -- (10,.2) node[above,purple] {$\Rmax$};
		
		\draw[domain=0.225:2.5, smooth, variable=\x, black, samples=100,thick] plot ({5*\x}, {-4./(1+\x*\x) + .2/(\x*\x)});	
		
		\draw[domain=.66:2.5, smooth, variable=\x, blue, samples=100,thick] plot ({5*\x}, {-4./(1+\x*\x) + 1.2859/(\x*\x)});
		
		\draw[ubtgreen, densely dotted,thick] (5.7215,-.75) -- (5.7215,.3) node[above,ubtgreen] {$r_{\Lmax}$};
		
		\node[black] at (6,-1.9) {$\Psi_{L_0}$};
		\node[blue] at (5.25,-1.05) {$\Psi_{\Lmax}$};
		
		\fill[opacity = .2, ubtgreen] (2.6833,.2) -- (2.6833,.0) -- (5.7215,.0) -- (5.7215,.2) -- cycle;
		\fill[opacity = .2, purple] (1.29099,-.2) -- (1.29099,.0) -- (10,.0) -- (10,-.2) -- cycle;
		
	\end{tikzpicture} 
	\caption{Shapes of~$L$-dependent effective potentials~$\Psi_L$. The figure illustrates the definitions of $\Rmin=r_-(E_0,L_0)$ and $\Rmax=r_+(E_0,L_0)$ as well as that $\Lmax>0$ is defined via $\Emin^{\Lmax}=E_0$. On the $r$-axis, it is further visualised that the trapping region $[r_{L_0},r_{\Lmax}]$ (green) is compactly contained in the radial steady state support $[\Rmin,\Rmax]$ (purple).}
	\label{F:PWL}
\end{figure}

Here we work with potentials~$\Psi_L$ of the form
\begin{equation}\label{E:PSIDEFL}
	\Psi_L(r)=-\frac Mr+U_0(r)+\frac L{2r^2},\qquad r>0,
\end{equation}
with the parameter $L>0$ which correspond to effective potentials induced by the steady states of the form~\eqref{E:POLY2}. 
With $M=0$ in~\eqref{E:PSIDEFL}, one obtains the effective potentials corresponding to regular steady states~\eqref{E:POLY}.
For each $L>0$, the graph of~$\Psi_L$ qualitatively resembles the graph of~$\Psi$ defined by~\eqref{E:PSIDEF} (see Figure~\ref{F:PW}).
It is straight-forward to prove~\cite{St23} that for any $L>0$ there exist a unique radius $r_L>0$ such that
\begin{equation*}
	\Psi_L'(r_L)=0, \qquad \EminL\coloneqq\Psi_L(r_L)<0,\qquad\Psi_L''(r_L)>0.
\end{equation*}
Moreover, for any $L>0$ and $\EminL<E<0$, there exist two unique solutions $r_-(E,L)<r_L<r_+(E,L)$ of $\Psi_L(r_\pm(E,L))=E$.
These radii determine the radial support of the steady state via $\Rmin=r_-(E_0,L_0)$ and $\Rmax=r_+(E_0,L_0)$.
The support of the squared modulus of angular momentum~$L$ is of the form $[L_0,\Lmax]$, where $\Lmax$ is uniquely determined via $\Emin^{\Lmax}=E_0$; see Figure~\ref{F:PWL} for a visualisation.

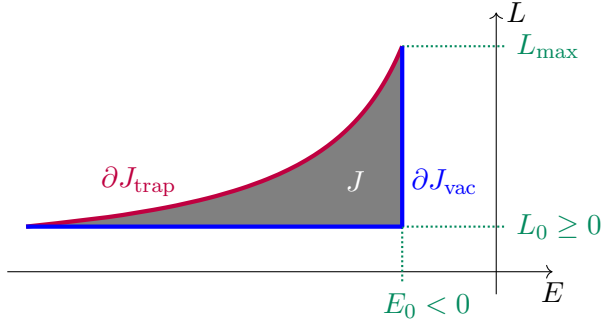
\begin{figure}[h!]
	
	\begin{center}
		\begin{tikzpicture}
			\tikzmath{\xscale = 2.5; \yscale = 3.;}
			
			\draw[-{>},black] (0.,-.1*\yscale) -- (0.,1.15*\yscale) node[right] {$L$};
			\draw[-{>},black] (-2.6*\xscale,0.) -- (.3*\xscale,0.) node[below] {$E$};

			\fill [gray, domain=.2:1.,smooth, variable=\x, samples=20, very thin]
			(-2.5*\xscale,.2*\yscale)
			-- plot ({-\xscale*.5/\x}, {\yscale*\x})
			-- (-.5*\xscale,.2*\yscale)
			-- cycle;
			
			\node[white,align=center] at (-.75*\xscale,.4*\yscale) {$J$};
			
			\draw[ubtgreen, densely dotted, thick] (-.5*\xscale,1.*\yscale) -- (-.5*\xscale,-.05*\yscale) node[below right,ubtgreen,align=center,xshift=-4.5*\xscale] {$E_0<0$};
			\draw[ubtgreen, densely dotted, thick] (-2.5*\xscale,.2*\yscale) -- (.05*\xscale,.2*\yscale) node[right,ubtgreen] {$L_0\geq0$}; 
			\draw[ubtgreen, densely dotted, thick] (-.5*\xscale,\yscale) -- (.05*\xscale,\yscale) node[right,ubtgreen] {$\Lmax$}; 
			
			\draw[purple,domain=.2:1.,smooth, variable=\x, samples=20, ultra thick]  plot ({-\xscale*.5/\x}, {\yscale*\x});
			%		\draw[purple,densely dotted,domain=1.:1.25,smooth, variable=\x, samples=10,thick]  plot ({-\xscale*.5/\x}, {\yscale*\x});
			\node[purple,align=center] at (-1.9*\xscale,.4*\yscale) {$\pa J_{\text{trap}}$};
			
			\draw[blue, ultra thick] (-2.5*\xscale,.2*\yscale) -- (-.5*\xscale,.2*\yscale);
			\draw[blue, ultra thick] (-.5*\xscale,.2*\yscale) -- (-.5*\xscale,\yscale);
			\node[blue,align=center] at (-.275*\xscale,.4*\yscale) {$\pa J_{\text{vac}}$};
			
			%		\draw[purple,densely dotted,domain=1.:1.25,smooth, variable=\x, samples=10,thick]  plot ({-\xscale*.5/\x}, {\yscale*\x});
		\end{tikzpicture}
	\end{center}
	\caption{The $(E,L)$-support~$J $ of a spherically symmetric steady state (gray). The purple and blue parts of the boundary of~$J$ depict~$\pa J_{\text{trap}}$ and $\pa J_{\text{vac}}$, respectively.}
	\label{F:ELtriangle}
\end{figure}

In the $(E,L)$-space, the interior of the steady state support is
\begin{equation}\label{E:DEFJ}
	J\coloneqq\{(E,L)\mid L>L_0,\;\EminL<E<E_0\}.
\end{equation}
As depicted in Figure~\ref{F:ELtriangle}, this set qualitatively resembles a ``triangle".
Its topological boundary naturally splits into
\begin{align}
  \pa J_{\text{vac}} &:= \{E_0\}\times[L_0,\Lmax] \ \cup\ 
  [E_{\min}^{L_0},E_0] \times\{L_0\} \label{E:JVACDEF}\\
  \pa J_{\text{trap}} & := \left\{(\EminL,L) \mid L\in[L_0,\Lmax]\right\}.
  \label{E:JTRAPDEF}
\end{align}
The set $\pa J_{\text{vac}}$ coincides with the vacuum boundary $\pa\Omega$
where the steady state $\mathring f$ vanishes to varying orders depending on the size of $k$ and $\ell$ in~\eqref{E:POLY2}. Geometrically, it is the union of the two straight sides of ~$J$, see Figure~\ref{F:ELtriangle}. 
The curve $\pa J_{\text{trap}}$ is the purple (curved) side of the action support triangle $J$ 
which corresponds to the $(E,L)$-pairs for which the associated dynamical system has a stagnation point.
We note that $k+\ell+\frac32>j\in\N$ implies~\cite{St23}
\begin{align}\label{E:PSIREGL}
\Psi_L\in C^{j+2}([\Rmin,\Rmax]).
\end{align}

Compared to Theorem~\ref{T:MAIN}, in this radially symmetric case, the force field generated by solutions to~\eqref{E:PT0L} decays faster, provided  the data is sufficiently regular. 

%%%%%%%%%%%%%%%%%%%%%%%%%%%%%%%%%%%%%%%%
%%%%%%%%%%%%%%%%%%%%%%%%%%%%%%%%%%%%%%%%

\begin{theorem}\label{T:MAIN2}
Let $k>1$ and $\ell>0$. 
For any  $g_0\in W^{2,\infty}(\Om)$, consider the initial data of the form
\begin{equation}
f_0(r,w,L) = |\varphi'(E(r,w,L),L)|\,g_0(r,w,L), \label{E:INITIALL}
\end{equation}
where $\varphi$ is given by~\eqref{E:POLY2}.
Let the operator $\D$ be given by~\eqref{E:DDEFL} with $r\mapsto \Psi_L(r)$ as described above.  
Let $f(t,\cdot,\cdot,\cdot)$ be the solution to~\eqref{E:PT0L} with initial data $f(0,r,w,L)=f_0(r,w,L)\in \text{{\em ker}}\,(\D)^\perp$. 
\begin{enumerate}
\item[(a)]
There exists a constant $C>0$ such that for all $t\in\R_+$,
\begin{align}\label{E:DECAYL}
|\pa_r U_f(t,r)|\leq C\|g_0\|_{W^{2,\infty}}\,\begin{cases}
	\frac{1+|\log t|}{t^2}, & k\geq 2,\, \ell\geq 1,\\
	\frac{1}{t^{\min\{\ell+1,k\}}}, &\text{else},
\end{cases}
\end{align}
and
\begin{align}\label{E:RHODECAYL}
|\rho_f(t,r)|\leq C\|g_0\|_{W^{1,\infty}}\frac{|\log t|}{t}.
\end{align}
\item[(b)] 
Assume that the data $g_0$ is supported strictly away from both the trapping set $\pa J_{\text{trap}}$ and the vacuum boundary $\pa J_{\text{vac}}$. Then, for any $g_0\in C^{j}(\overline \Om)$, where $j\in\mathbb N$ with $j\ge2$, there exists a constant $C=C_j$ such that 
\begin{equation}
\|\pa_rU_f(t,\cdot) \|_{L^\infty_r} \le \frac{C\|g_0\|_{W^{j,\infty}}}{1+t^{j}}.  \label{E:DECAYBETTER1}
\end{equation}
\end{enumerate}
\end{theorem}

%%%%%%%%%%%%%%%%%%%%%%%%%%%%%%%%%%%%%%%%
%%%%%%%%%%%%%%%%%%%%%%%%%%%%%%%%%%%%%%%%

The proof of Theorem~\ref{T:MAIN2} is presented in Section~\ref{S:MAIN2}.

%%%%%%%%%%%%%%%%%%%%%%%%%%%%%%%%%%%%%%%%
%%%%%%%%%%%%%%%%%%%%%%%%%%%%%%%%%%%%%%%%

%%%%%%%%%%%%%%%%%%%%%%%%%%%%%%%
%%%%%%%%%%%%%%%%%%%%%%%%%%%%%%%

\begin{remark}
Part (b) of the theorem states that once the trapping is removed from the dynamics, the only obstruction to higher algebraic decay comes from the regularity of the initial data and the regularity of the potential $\Psi_L$. In particular, since our potentials are $C^\infty$ inside $J$ (but of finite regularity on $\bar J$), one can prove arbitrarily fast decay for arbitrarily smooth initial data assuming that the data is supported strictly away from the trapping set, or that it vanishes there to infinite order.

Even though the no-trapping assumption in part (b) is preserved by the dynamics of~\eqref{E:PT0}, it is {\bf not} preserved by the dynamics of the true linearised Vlasov-Poisson dynamics~\eqref{E:FULLLIN} around the equilibria~\eqref{E:POLY2}. It is therefore very important to understand the decay without such an assumption; as we see from part (a), the presence of trapping puts a severe limitation on the decay rate. 
\end{remark}

%%%%%%%%%%%%%%%%%%%%%%%%%%%%%%%
%%%%%%%%%%%%%%%%%%%%%%%%%%%%%%%

\subsection{Action-angle variables}

The basic tool to study  decay is the action-angle variables  associated to the steady state. These allow us to 
reparametrise the support~$\Omega$  of the steady state $\mathring f$  and
 provide a convenient way of describing the trapping of the particles in the potential well of the steady state, see Figures~\ref{F:PW} and~\ref{F:PWL}.

%%%%%%%%%%%%%%%%%%%%%%%%%%%%
%%%%%%%%%%%%%%%%%%%%%%%%%%%%
{\bf Full radial case.}
In the action-angle description,
the actions $(E,L)$ range in (the closure of) the set~$J$; recall~\eqref{E:DEFJ} and Figure~\ref{F:ELtriangle}.
For fixed $L>0$ and $\EminL\leq E<0$, let $(R,W)(\cdot,E,L)\colon\R\to(0,\infty)\times\R$ be the global solution to the characteristic system
\begin{equation}\label{E:CHFLOWL}
	\dot r=w,\qquad\dot w=-\Psi_L'(r), \qquad \dot L=0
\end{equation}
satisfying the initial condition 
%\begin{equation*}
	$(R,W)(0,E,L)=(r_-(E,L),0)$. 
%\end{equation*}
In the case $E=E^L_{\min}$ we set $r_-(E^L_{\min},L)\coloneqq r_L$ and there the solution takes the form $(R,W)(\cdot,E^L_{\min},L)\equiv(r_L,0)$.
For $0>E>E^L_{\min}$, the solution $(R,W)(\cdot,E,L)$ is time-periodic with period
\begin{equation*}
	T(E,L)\coloneqq2\int_{r_-(E,L)}^{r_+(E,L)}\frac{\diff r}{\sqrt{2E-2\Psi_L(r)}}.
\end{equation*}
Due to~\eqref{E:PSIREGL}, $T\colon\bar J\to(0,\infty)$ can shown to be in $C^{j+1}(\bar J)$ with
\begin{equation}\label{eq:TminexplicitL}
	T(E^L_{\min},L)=\frac{2\pi}{\sqrt{\Psi''(r_L)}},
\end{equation}
cf.~\cite{Ku21,St23}.
Moreover, in the regime $\k\ll1$ in~\eqref{E:POLY3}, it is shown in~\cite{St23} that $\partial_ET(E,L)>\Tmin'(L)\coloneqq \pa_ET(E^L_{\min},L)>0$ for $(E,L)\in J$. 
This motivates the following crucial monotonicity assumption
\begin{equation}\label{E:MONOTONEPERIODL}
\pa_ET(E,L)>0, \qquad (E,L)\in \bar J. 
\end{equation} 
Additionally, we assume that there exists a sufficiently small constant $\eta>0$ such that 
\begin{align}\label{E:PAL}
|\pa_LT(E,L)|\le\eta, \ \ (E,L)\in\bar J.
\end{align}
Assumption~\eqref{E:PAL} is satisfied for example in the case $0<\varepsilon\ll1$ in~\eqref{E:POLY2} by~\cite{St23}, where it can be shown that $\eta=\eta(\k)\to0$ as $\k\to0$. It is also trivially satisfied for the class of potentials  where the period~$T$ is $L$-independent. This happens for example in the limiting case $\k=0$, which corresponds to the gravitational point mass potential.
\begin{remark}
Numerical simulations~\cite{St24} indicate that \eqref{E:MONOTONEPERIODL}--\eqref{E:PAL}
hold for many steady states of the form~\eqref{E:POLY}.
\end{remark}

We let now for any $L\in[L_0,\Lmax]$
\begin{equation*}
	(r,w)\colon\mathbb S^1\times(\EminL,0)\to(0,\infty)\times\R,\;(r,w)(\theta,E,L)\coloneqq(R,W)(\theta\,T(E,L),E,L).
\end{equation*}
This mapping defines a one-to-one change of variables between $(r,w,L)\in\Omega^\ast\coloneqq\Omega\setminus\{(r_L,0,L)\mid L\in[L_0,\Lmax]\}$ and $(\theta,E,L)\in\mathbb S^1\times J$, cf.~\cite{HRS2022,Ku21}.
The variables $(\theta,E,L)$ are called {\em action-angle variables}.
By slight abuse of notation, we write $g=g(r,w,L)=g(\theta,E,L)$ for functions on~$\Omega^\ast$.
In these variables, equation~\eqref{E:PT0L} takes the form
\beq\label{E:PURETRANSPORTL}
\pa_t f+\om(E,L)\pa_\th f=0,
\eeq
cf.~\cite{HRS2022,Ku21}, where 
\begin{equation*}
\om(E,L) = \frac1{T(E,L)}, \ \  (E,L)\in \bar J,
\end{equation*}
defines the frequency function.

\begin{remark}\label{R:KERD}
As shown in~\cite{HRS2022,Ku21}, the orthogonal complement of the nullspace of~$\D$ in~$H$ is very easily characterised in action-angle variables: for $f\in H$,
\begin{equation}\label{E:ORTH}
f\in  (\mathrm{ker}\,\D)^\perp \ \ \text{{\em iff}} \ \ \int_{\mathbb S^1}f(\theta,E,L)\diff\theta = 0 \text{ for a.e.}\ (E,L)\in J.
\end{equation} 
\end{remark}

{\bf 1+1-dimensional case.}
In this case we define the action-variables analogously to above, where we fix $L=\bar L$.
The energies (actions) $E$ range in the interval $[\Emin,E_0]$ and we denote its interior by
\begin{equation*}
I\coloneqq(\Emin,E_0). 
\end{equation*}
The associated characteristic system coincides with~\eqref{E:CHFLOWL} with a fixed choice $L=\bar L$ and 
for $E\in(\Emin,0)$, we define the radii $r_-(E)<r_+(E)$ as the unique solutions of $\Psi(r_\pm(E))=E$, corresponding to the turning points of the periodic 
orbits.
In the case $E=\Emin$ we set $r_-(\Emin)\coloneqq r_\ast$ and there the solution takes the form $(R,W)(\cdot,\Emin)\equiv(r_\ast,0)$.
Like above, the period is given by the formula
$T(E)\coloneqq2\int_{r_-(E)}^{r_+(E)}\frac{\diff r}{\sqrt{2E-2\Psi(r)}}$
for $\Emin<E<0$.  
Due to the assumptions~\eqref{E:ELLIPTIC}--\eqref{E:PSIREG} and the behaviour of~$\Psi$ close to $\{\Rmin,\Rmax\}$, the mapping $T\colon\bar I\to(0,\infty)$
can be shown~\cite{HRSS2023,Ku21,St23} to be in $C^3(\bar I)$ with
\begin{equation}\label{eq:Tminexplicit}
	\Tmin=\frac{2\pi}{\sqrt{\Psi''(r_\ast)}}.
\end{equation}
Moreover, in the regime $\varepsilon\ll1$ in~\eqref{E:POLY3}, we have shown in~\cite{HRSS2023} that $\Tmin'\coloneqq T'(\Emin)>0$. 
By analogy to~\eqref{E:MONOTONEPERIODL} we impose the monotonicity assumption
\begin{equation}\label{E:MONOTONEPERIOD}
T'(E)>0, \qquad E\in \bar I,
\end{equation} 
which is satisfied~\cite{HRSS2023} in the case $\varepsilon\ll1$ in~\eqref{E:POLY3} because the smallness of $\k$ implies that the period function to the leading order corresponds to the period function of the pure point mass potential, which is proportional to $|E|^{-\frac32}$.
Analogously to the above the map
\begin{equation*}
	(r,w)\colon\mathbb S^1\times(\Emin,0)\to(0,\infty)\times\R,\;(r,w)(\theta,E)\coloneqq(R,W)(\theta\,T(E),E)
\end{equation*}
defines a one-to-one change of variables between $(r,w)\in{\Omega}^\ast\coloneqq \Omega\setminus\{(r_\ast,0)\}$ and $(\theta,E)\in\mathbb S^1\times I$ and we then 
refer to $(\theta,E)$ as the action-angle variables.
Now the equation~\eqref{E:PT0} takes the form
\beq\label{E:PURETRANSPORT}
\pa_t f+\om(E)\pa_\th f=0,
\eeq
cf.~\cite{HRS2022,HRSS2023,Ku21}, where 
\begin{equation*}
\om(E) = \frac1{T(E)},\qquad\Emin\leq E<0,
\end{equation*}
is the frequency function.

%%%%%%%%%%%%%%%%%%%%%%%%%
%%%%%%%%%%%%%%%%%%%%%%%%%

%%%%%%%%%%%%%%%%%%%%%%%%%
%%%%%%%%%%%%%%%%%%%%%%%%%

\subsection{Related literature and ideas of the proof}

%%%%%%%%%%%%%%%%%%%%%%%%%
%%%%%%%%%%%%%%%%%%%%%%%%%

The basic mechanism 
that leads to the decay of solutions to~\eqref{E:PT0L}, called {\em gravitational phase mixing},
was discovered by Lynden-Bell in his pioneering discussion~\cite{LB1962} in 1962. In particular Lynden-Bell used the action-angle variables
to suggest an equilibration of action-averaged observables associated to the pure transport flow (see also~\cite{Kalnajs1971,BaYa}).
At a rigorous level, phase mixing for equations of the form~\eqref{E:PT0L} or~\eqref{E:PT0} was first studied by Rioseco and Sarbach~\cite{RiSa2020} who obtained
a weak convergence result for action averages of the solutions of pure transport equations. Such a result is morally related to
$H^{-1}$-mixing in the spirit of ergodic theory. In 1-d,  Chaturvedi and Luk~\cite{ChLu2022} proved quantitative decay rates for the gravitational/electrostatic
potential and its derivative, where an explicit external potential is imposed. The associated period function satisfies the qualitative features satisfied by our $E\mapsto T(E)$.  The authors  make the assumption that the initial data are supported uniformly away from the elliptic point $E=\Emin$ and use the vector-field method to obtain $\frac1{1+t}$ decay rate for the associated force field. Also using the vector-field method, quantitative decay rates for action averages were obtained
by Moreno, Rioseco, and Van Den Bosch~\cite{MoRiVa2022}; there the authors instead use test functions supported strictly away from the elliptic point. In particular, the trapping effects at $(r,w)=(r_\ast,0)$
are effectively removed from~\cite{ChLu2022,MoRiVa2022}.

%%%%%%%%%%%%%%%%%%%%%%%%%
%%%%%%%%%%%%%%%%%%%%%%%%%

Equations~\eqref{E:PURETRANSPORTL} and~\eqref{E:PURETRANSPORT} can be solved explicitly via the formulas
\begin{align}
f(t,\theta,E,L) &= f_0(\theta-\omega(E,L)t,E,L), \ \ f_0(\theta,E,L): = f(0,\theta,E), \label{E:EXPLICITL}\\
f(t,\theta,E) & = f_0(\theta-\omega(E)t,E), \ \ f_0(\theta,E): = f(0,\theta,E), \label{E:EXPLICIT}
\end{align}
respectively.
We denote the mixing semigroups appearing in~\eqref{E:EXPLICITL} and \eqref{E:EXPLICIT} by $T_t$ and~$\tilde T_t$, respectively, i.e.,  for any $t\ge0$ we let
\begin{align}
T_t(\theta, E,L) &: = (\theta-\omega(E,L)t,E,L), \label{E:SEMIGROUPL} \\
\tilde T_t(\theta, E) &: = (\theta-\omega(E)t,E). \label{E:SEMIGROUP}
\end{align}

One of the main advantages of the action-angle variables is the availability of 
Fourier analysis in the angle $\theta\in\mathbb S^1$.
As usual, for a suitable function~$g=g(\theta,E,L)$ on~$\Omega^\ast$, the Fourier coefficients are given by
\begin{equation}\label{E:DEFFOURIERCOEF}
	\hat g(m,E,L)\coloneqq\int_{\mathbb S^1}g(\theta,E,L)\,e^{-2\pi im\theta}\diff\theta,\qquad m\in\Z,\;(E,L)\in J,
\end{equation}
and the similar definition applies to functions $g=g(\theta,E)$.
Applying the Fourier transform in~$\th$ to equation~\eqref{E:PURETRANSPORT} gives
\begin{equation}\label{E:MODEBYMODEL}
\pa_t\fhat(t,m,E,L)+2\pi i m \om(E,L)\fhat(t,m,E,L)=0.
\end{equation}
We can solve this explicitly to obtain the mode-by-mode formula for the semigroup propagator:
\begin{equation}\label{E:FHATREPL}
\fhat(t,m,E,L) = e^{-2\pi i m \om(E,L) t} \fhat_0(m,E,L), \ \ m\in\Z, \ (E,L)\in J.
\end{equation}
Of course, the analogous formula holds in the 1+1-dimensional case
\begin{equation}\label{E:FHATREP}
\fhat(t,m,E) = e^{-2\pi i m \om(E) t} \fhat_0(m,E), \ \ m\in\Z, \ E\in I.
\end{equation}

%%%%%%%%%%%%%%%%%%%%%%%%%
%%%%%%%%%%%%%%%%%%%%%%%%%

One of the most interesting features of phase-mixing generated by the potential well above is the presence of 
the so-called elliptic point $E=\Emin$ in the 1+1-d case~\eqref{E:PT0} or a curve of elliptic points $L\mapsto \EminL$ in the general case~\eqref{E:PT0L} at the smallest particle energy supported by the steady state of a given angular momentum $L$. 

For example, in the restricted case~\eqref{E:PT0}, the stars that start on the cylinder $r=r_\ast$ with radial velocity $w=0$ always stay there and present a possible obstruction to decay. 
A manifestation of this is 
that the 
change of variables
$
(r,w) \mapsto (\theta,E)
$
becomes singular precisely at the elliptic point. Thus with trapping comes a loss of smoothness at $E=\Emin$, which
places severe restrictions on the decay rates. Various regularity properties of the action-angle variables and quantitative
estimates near $E=\Emin$ and the turning points of the flow $w=0$ are studied in Section~\ref{S:AAREG}. 
An important consequence is that the decay rates stated in Theorem~\ref{T:MAIN} do not improve if we
assume more regularity on the initial density distribution $f_0$.
This is in sharp contrast to the well-known phenomenon of phase-mixing near 
space-homogeneous equilibria of the Vlasov-Poisson system~\cite{MoVi2011}, where improved regularity of initial data leads to faster decay rates.

Our  approach to the proof of Theorem~\ref{T:MAIN} is classical. We use the explicit formulas~\eqref{E:EXPLICIT} and~\eqref{E:FHATREP}
and then resort to a nonstationary phase type argument to obtain quantitative decay rates. 
Schematically speaking, a good model problem is to think of the time-decay of  integrals of the general form
\be\label{E:NSP}
\int_a^b e^{i \om(E) t} F(E)\diff E.
\ee
By a simple integration
by parts argument we obtain $\frac1t$ decay if $\pa_E F(E)$ is in $L^1$, the boundary terms $F(a)$ and $F(b)$ vanish, and $\om'\neq0$ which is where the monotonicity assumption~\eqref{E:MONOTONEPERIOD} comes in crucially.
We will see that the trapping referred to above can cause either of the first two properties to fail, thus placing a limitation on decay rates for generic initial data. 
We comment that the trapping at $E=\Emin$ can be removed by choosing initial data supported away from the minimal energy $E=\Emin$. However, this property is not preserved by the full linearised flow~\eqref{E:FULLLIN} which highlights the need to understand the decay properties in the presence of trapping.
The key tool in our analysis is a precise description of the fundamental solution to the Poisson equation in radial symmetry, using the action-angle variables. This is essentially the content of Section~\ref{S:AAREG1D}, where we carefully track the behaviour of the action-angle $\theta(r,E)$ and its $E$-derivatives near the turning points $E=\Psi(r)$, where they degenerate. 

In the general case~\eqref{E:PT0L} the situation is slightly better, as there is more ``room" for efficient phase-mixing. The schematic analogue of~\eqref{E:NSP}
is now 
\be\label{E:NSPL}
\int_J e^{i \om(E,L) t} F(E,L)\diff (E,L).
\ee
In addition to the ideas highlighted above, we take advantage of
the presence of the angular momentum which allows us to effectively trade $\pa_E$ and $\pa_L$ derivatives to improve upon decay.

{\bf Plan of the paper.} In Section~\ref{S:AAREG} we study various regularity properties of the action-angle variables. In Section~\ref{S:DECAY} we prove Theorem~\ref{T:MAIN} and in Section~\ref{S:MAIN2} we prove Theorem~\ref{T:MAIN2}.

\bigskip

{\bf Acknowledgments.}
M. Had\v zi\'c's research is supported by the EPSRC Early Career Fellowship EP/S02218X/1.
M. Schrecker's research is supported by the EPSRC Post-doctoral Research Fellowship EP/W001888/1.

%%%%%%%%%%%%%%%%%%%%%%%%%%%%
%%%%%%%%%%%%%%%%%%%%%%%%%%%%
%%%%%%%%%%%%%%%%%%%%%%%%%%%%

\section{Regularity in action-angle variables}\label{S:AAREG}

%%%%%%%%%%%%%%%%%%%%%%%%%%%%
%%%%%%%%%%%%%%%%%%%%%%%%%%%%

The behaviour of various key quantities near the elliptic point of the flow is crucial to the proof of quantitative decay rates. 
The purpose of this section is to  capture the singular behaviour near $E=\Emin$ through a series of lemmas which will then be used to prove the main theorems.

\subsection{The 1+1-dimensional case}\label{S:AAREG1D}

We start with the 1+1-d case and emphasise that all results established here can be trivially modified to also capture the $L$-dependent setting. 
Further extensions to the full radial setting will be proven in Section~\ref{S:AAREG3D}. 
We start with a basic observation which clarifies the behaviour of Fourier coefficients~\eqref{E:DEFFOURIERCOEF} of smooth functions near $E=\Emin$.

%%%%%%%%%%%%%%%%%%%%%%%%%%%%%%%%%
%%%%%%%%%%%%%%%%%%%%%%%%%%%%%%%%%

\begin{lemma}\label{lemma:Eminvanishing}
Let $g=g(r,w)\in W^{k,\infty}(\Om)$ for some $k\in\{1,2\}$.
Then there exists a constant $C>0$, depending on the steady state, such that
\begin{equation}\label{E:ghatestimate}
	|\widehat g(m,E)|\leq C\frac{\|g\|_{W^{k,\infty}(\Om)}\,\sqrt{E-E_{\min}}}{|m|^k},\qquad m\in\Z^\ast,\;E\in I.
\end{equation}
\end{lemma}

%%%%%%%%%%%%%%%%%%%%%%%%%%%%%%%%%
%%%%%%%%%%%%%%%%%%%%%%%%%%%%%%%%%

\begin{proof}
We first recall from \cite[Lemma 3.5]{HRSS2023} that there exists $C>0$ as above such that
\begin{equation*}
|r(\th,E)-r_*|+|\pa_\th r(\th,E)|\leq C\sqrt{E-E_{\min}},\qquad(\theta,E)\in\mathbb S^1\times I.
\end{equation*}
Moreover, from the definition of the characteristic flow, we have 
\begin{equation*}
|\pa_\th w(\th,E)|=T(E)\,|\Psi'(r(\th,E))|\leq C|r(\th,E)-r_*|\leq C\sqrt{E-E_{\min}}.
\end{equation*}
For $m\neq 0$, we make the trivial estimate
$|m|^{2k}|\widehat{g}(m,E)|^2\leq\int_{\mathbb{S}^1}|\pa_\th^kg(\th,E)|^2\,\dif\th$. Thus,
\begin{align*}
	|\widehat{g}(m,E)|^2\leq&\frac{1}{m^2}\int_{\mathbb{S}^1}|\pa_\th g(\th,E)|^2\,\dif\th\nonumber\\
	\leq&\,\frac{2}{m^2}\int_{\mathbb{S}^1}\big(|(\pa_r g)(\th,E)|^2|\pa_\th r|^2+|(\pa_w g)(\th,E)|^2|\pa_\th w|^2\big)\,\dif\th\\
	\leq&\,\frac{C}{m^2}\|g\|_{W^{1,\infty}(\Om)}^2\,(E-E_{\min}),
\end{align*}
which shows the estimate in the case $k=1$. 
To prove the case $k=2$ in~\eqref{E:ghatestimate}, we apply the chain rule once more and use that
\begin{equation*}
	|\pa_\th^2r(\th,E)|+|\pa_\th^2w(\th,E)|\leq C\sqrt{E-\Emin}.\qedhere
\end{equation*}
\end{proof}

%%%%%%%%%%%%%%%%%%%%%%%%%%%%%%%%%
%%%%%%%%%%%%%%%%%%%%%%%%%%%%%%%%%

By assuming further regularity of~$g$ (and suitable smoothness of the potential~$\Psi$), one could obtain arbitrarily many additional $\frac1{|m|}$-factors on the right-hand side of~\eqref{E:ghatestimate}.
However, this would not lead to faster decay due to the trapping mechanisms explained before.

Let us next derive suitable representations for the functions whose decay we aim to prove.

%%%%%%%%%%%%%%%%%%%%%%%%%%%%%%%%%
%%%%%%%%%%%%%%%%%%%%%%%%%%%%%%%%%

\begin{lemma}[Representation formulas for the gravitational potential]\label{L:REPRESENTATIONS}
Let $t\mapsto f(t,\cdot,\cdot)$ be a classical solution of the pure transport equation~\eqref{E:PT0} and let the gravitational potential $U_f$ be given
by~\eqref{E:POISSON}. Then the following formulae hold:
\begin{align}
\pa_t U_f(t,R)& = 4\pi^2 \iint  f(t,r,w)\,wg_R(r) \diff(r,w), \label{E:PARTIALTUFFORMULA} \\
\pa_R U_f(t,R) & = \frac{4\pi^2}{R^2} \iint f(t,r,w)\,h_R(r)\diff(r,w), \label{E:PARTIALRUFFORMULA}
\end{align}
where
\begin{align}
g_R(r)&\coloneqq \frac1{r^2}\chi_{r\ge R}(r), \label{E:GRDEF}\\
h_R(r) &\coloneqq \chi_{r\le R}(r) \label{E:HRDEF}.
\end{align}
\end{lemma}

\begin{proof}
The formula~\eqref{E:PARTIALRUFFORMULA} follows classically by radially integrating the Poisson equation.
Differentiating~\eqref{E:PARTIALRUFFORMULA} w.r.t.~$t$ and radially integrating the resulting identity yields~\eqref{E:PARTIALTUFFORMULA}; recall the boundary condition included into~\eqref{E:POISSON}.
\end{proof}

%%%%%%%%%%%%%%%%%%%%%%%%%%%%%%%%
%%%%%%%%%%%%%%%%%%%%%%%%%%%%%%%%

\begin{remark}
As explained in the discussion of the properties of the set $\Om$~\eqref{E:OMDEF} in Section~\ref{S:INTRO}, the spatial extent of $\Om$ is contained in the interval $[\Rmin,\Rmax]$ with $\Rmin>0$.
We will therefore frequently use the bound $\frac1r<C$ on $\Om$.
\end{remark}

%%%%%%%%%%%%%%%%%%%%%%%%%%%%%%%%
%%%%%%%%%%%%%%%%%%%%%%%%%%%%%%%%

Motivated by Lemma~\ref{L:REPRESENTATIONS}, it will be particularly 
important to  describe the set $\{r(\theta,E)\ge R\}$ in the action-angle variables for any given $R\in[\Rmin,\Rmax]$. 
For $E\in[\Psi(R),E_0]$, let
$\theta(R,E)\in[0,\frac12]$ be the unique angle such that 
\begin{align}\label{E:THETADEF}
r(\theta(R,E),E)= R.
\end{align}
In particular, for $E\in(\Emin,E_0]$ and $R\geq r_\ast$,
\begin{equation}\label{E:rthetaEgeqR_Rlowerrast}
r(\theta,E)\geq R \ \ \text{iff } \ E\in[\Psi(R),E_0]\land\theta\in[\theta(R,E) ,1-\theta(R,E)],
\end{equation}
while for $R<r_\ast$,
\begin{equation}\label{E:rthetaEgeqR_Rlargerrast}
	r(\theta,E)\geq R\ \ \text{iff } \ E\leq\Psi(R)\lor\left(E\in[\Psi(R),E_0]\land\theta\in[\theta(R,E),1-\theta(R,E)]\right).
\end{equation}
Note that for all $R\in[\Rmin,\Rmax]$ and $E\in(\Psi(R),E_0]$ we have $\theta(R,E)>0$.
Moreover,
\begin{equation} \label{E:WVANISH}
\theta(R,\Psi(R))=
\begin{cases}
\frac12, & \ R>r_\ast, \\
0, & \ R<r_\ast,
\end{cases}\qquad R\in[\Rmin,\Rmax]\setminus\{r_\ast\}.
\end{equation}
This equation expresses the fact that
the angles $\theta=0$ and $\theta=\frac12$ correspond to the turning points of any characteristic orbit.

%%%%%%%%%%%%%%%%%%%%%%%%%%%%%%%%%%%%%%%
%%%%%%%%%%%%%%%%%%%%%%%%%%%%%%%%%%%%%%%

\begin{lemma}[Regularity of $E\mapsto \theta(R,E) $]\label{L:THETAREG}
For any $R\in(\Rmin,\Rmax)$, the map $(\Psi(R),E_0]\ni E\mapsto \theta(R,E)$ is continuously differentiable and there exists a constant $C>0$, depending only on the steady state, such that
\begin{equation*}
\lv \pa_E\theta(R,E)  \rv \le \frac{C}{\sqrt{E-\Psi(R)}\sqrt{E-\Emin}}, \qquad E\in(\Psi(R),E_0].
\end{equation*}
\end{lemma}

%%%%%%%%%%%%%%%%%%%%%%%%%%%%%%%%%%%%%%%

\begin{proof}
Basic ODE theory shows $r\in C^1(\mathbb S^1\times I)$ and for $(\theta,E)\in\mathbb S^1\times I$ we have $\partial_\theta r(\theta,E)=0$ iff $\theta\in\{0,\frac12\}$.
By~\eqref{E:THETADEF}, \eqref{E:WVANISH}, and the implicit function theorem we hence deduce $\theta(R,\cdot)\in C^1((\Psi(R),E_0])$ with
\begin{equation}\label{E:THETADEF2}
T(E)\,w(\theta(R,E),E)\,\partial_E\theta(R,E)  = - \pa_Er(\theta(R,E)  ,E).
\end{equation}
We now recall the bound~\cite[Lemma 3.5]{HRSS2023}:
\begin{align}\label{E:PARTIALER}
\lv\pa_Er(\theta,E)\rv + \lv\pa_Ew(\theta,E)\rv \le \frac C{\sqrt{E-\Emin}}, \qquad E\in I,\;\theta\in\mathbb S^1;
\end{align}
the bound for $\pa_Ew$ is not explicitly stated in~\cite{HRSS2023}, but follows with the same arguments as the bound for $\pa_Er$ and will be useful later.
Together with $\Tmin>0$, the claim follows.
\end{proof}

%%%%%%%%%%%%%%%%%%%%%%%%%%%%%%%%%%%%%%%
%%%%%%%%%%%%%%%%%%%%%%%%%%%%%%%%%%%%%%%

The next lemma gives a precise description of the leading order rates of degeneration of the characteristic flow as one approaches $E=\Emin$.

\begin{lemma}\label{L:HOLDERONEHALF}
Let $\alpha\coloneqq\sqrt{\Psi''(r_\ast)}$. For $\theta\in\mathbb S^1$, the following bounds hold uniformly in $\theta\in\mathbb S^1$ as $E\searrow\Emin$:
\begin{align}\label{E:HOLDERONEHALFR}
\frac{r(\theta, E) - r_\ast}{(E-\Emin)^{\frac12}} &= -\frac{\sqrt 2 \cos(2\pi\theta)}{\alpha} + \frac{\Psi'''(r_\ast)}{6\alpha^4}\left(\cos(4\pi\theta)-3\right)\sqrt{E-\Emin}+\mathcal O(E-\Emin),\\
\label{E:HOLDERONEHALFW}
\frac{w(\theta,E)}{(E-\Emin)^{\frac12}}&=\sqrt2\sin(2\pi\theta)-\frac{\Psi'''(r_\ast)}{3\alpha^3}\sin(4\pi\theta)\sqrt{E-\Emin}+\mathcal O(E-\Emin).
\end{align}
\end{lemma}
\begin{proof}
We will only show here how to derive the expansions up to the order $\sqrt{E-\Emin}$; the higher order behaviour can be obtained by iterating the arguments below.

For $\Emin\leq E<0$ let $z=z(\cdot,E)\colon\R\to\R$ be the unique global solution of
\begin{equation}\label{eq:ODEz}
	\ddot z=-\Psi''(R(\cdot,E))\,z,\qquad z(0)=1,\;\dot z(0)=0.
\end{equation} 
In particular, because $R(\cdot,\Emin)\equiv r_\ast$ it follows that
\begin{equation*}
z(s,\Emin) = \cos(\alpha s),\qquad s\in\R.
\end{equation*}
We now consider the function $Z(s,E)\coloneqq z(s,E) - \cos(\alpha s)$ for $s\in\R$ and $E\in [\Emin,0)$ which solves the nonhomogeneous linear ODE 
\begin{equation*}
	\ddot Z(s,E) = - \alpha^2 Z (s,E)  + \left(\alpha^2 - \Psi''(R(s,E))\right)z(s,E)
\end{equation*}
with initial condition $Z(0,E)=0=\dot Z(0,E)$.
Denoting the corresponding semigroup by $e^{s A}$, $A = \begin{pmatrix} 0 & 1 \\ -\alpha^2 & 0 \end{pmatrix}$, it follows that
\begin{align}\label{E:ZDUHAMEL}
\begin{pmatrix} Z(s) \\ \dot Z(s) \end{pmatrix} = \int_0^s e^{A(s-\sigma)} \begin{pmatrix} 0 \\ \left(\alpha^2 - \Psi''(R(\sigma,E))\right) z(\sigma,E) \end{pmatrix} \diff \sigma,\qquad s\in\R. 
\end{align}
Since $|z(s,E)|\leq C$ for $s\in[0,\Tmax]$ and $E\in[\Emin,E_0]$ and $\lv \alpha^2 - \Psi''(R(s,E))\rv \le C \sqrt{E-\Emin}$, it follows 
from~\eqref{E:ZDUHAMEL} that 
\begin{align}\label{E:BIGZBOUND}
\lv Z(s,E)\rv + \lv \dot Z(s,E)\rv \le C \sqrt{E-\Emin},\qquad E\in[\Emin,E_0],\;s\in[0,\Tmax].
\end{align}

For any $E,\tilde E\in I$ we then have via the mean value theorem
\begin{align}
r(\theta, E) - r(\theta,\tilde E) & = \int_{\tilde E}^E \pa_e r(\theta, e) \diff e \notag\\
& =  \int_{\tilde E}^E  W(\theta\,T(e),e)\,\theta\,T'(e)+\partial_eR(\theta\,T(e),e) \diff e \notag\\
& = \theta \int_{\tilde E}^E  W(\theta\,T(e),e)\,T'(e) \diff e+ \int_{\tilde E}^E \partial_er_-(e) z(\theta\,T(e),e) \diff e. \label{E:RDIFFQUOTIENT}
\end{align}
Note that 
\begin{align}
z(\theta T(e),e) & = z(\theta \Tmin,\Emin) + z(\theta T(e),e)-z(\theta \Tmin,e) + Z(\theta \Tmin, e) \notag\\
& = \cos(2\pi\theta) + z(\theta T(e),e)-z(\theta \Tmin,e) + Z(\theta \Tmin, e) \label{E:LITTLEZDECOMPOSITION}
\end{align}
since $\Tmin=T(\Emin)=\frac{2\pi}{\alpha}$ by~\eqref{eq:Tminexplicit}.
Since $\dot z$ is uniformly bounded and $T\in C^1(\bar I)$, it follows that $\lv z(\theta T(e),e)-z(\theta \Tmin,e)\rv\le C (e-\Emin)$.
Together with~\eqref{E:BIGZBOUND} and~\eqref{E:LITTLEZDECOMPOSITION} we hence deduce
\begin{equation}\label{E:ZLEADINGORDER}
z(\theta T(e),e) = \cos(2\pi\theta) + \mathcal O(\sqrt{e-\Emin}),
\end{equation}
and the bound is uniform in~$\theta\in\mathbb S^1$.
We also recall that 
\begin{align}
\pa_e r_-(e) & = \frac1{\Psi'(r_-(e))} = \frac1{\Psi''(r_\ast)(r_-(e)-r_\ast)} \left(1+\mathcal O_{e\to\Emin}(\sqrt{e-\Emin})\right) \notag\\
& = -\frac{1}{\sqrt 2\alpha (e-\Emin)^{\frac12}}\left(1+\mathcal O_{e\to\Emin}(\sqrt{e-\Emin})\right) .\label{E:PARTIALRMINUSEXPANSION}
\end{align}
From the uniform boundedness of $W$,~\eqref{E:RDIFFQUOTIENT},~\eqref{E:ZLEADINGORDER}, and~\eqref{E:PARTIALRMINUSEXPANSION} we conclude
\begin{align}
r(\theta, E) - r(\theta,\tilde E) =&  \int_{\tilde E}^E  \mathcal O_{e\to\Emin}(1) \diff e + \cos(2\pi\theta) \int_{\tilde E}^E \pa_e r_-(e)\diff e \nonumber\\
&\qquad\qquad\qquad+\int_{\tilde E}^E   \pa_e r_-(e) \mathcal O(\sqrt{e-\Emin}) \diff e \nonumber\\
\label{E:RTHETAEXPANSIONINT} &= -\frac{ \cos(2\pi\theta)}{\sqrt 2\alpha} \int_{\tilde E}^E (e-\Emin)^{-\frac12}\diff e + \mathcal O(E-\tilde E).
\end{align}
We now let $\tilde E\searrow\Emin$ and upon dividing by $(E-\Emin)^{\frac12}$ we obtain the first expansion term in~\eqref{E:HOLDERONEHALFR}.
In order to obtain~\eqref{E:HOLDERONEHALFW}, first note $\partial_EW(\theta T(E),E)=\partial_Er_-(E)\,\dot z(\theta T(E),E)$. We hence proceed as above and write
\begin{equation}\label{E:WDIFFQUOTIENT}
	w(\theta,E)-w(\theta,\tilde E)=-\theta\int_{\tilde E}^E  \Psi'(r(\theta,e))\,T'(e) \diff e+\int_{\tilde E}^E \partial_er_-(e)\,\dot z(\theta\,T(e),e) \diff e
\end{equation}
for $E,\tilde E\in I$. 
After letting $\tilde E\searrow\Emin$, the first term on the right-hand side is $\mathcal O(E-\Emin)$. For the second one, differentiating~\eqref{E:LITTLEZDECOMPOSITION} with respect to $\th$ gives
\begin{equation*}
\Tmin\dot z(\theta T(e),e) = -\al \Tmin\sin(2\pi\theta) + \Tmin\big(\dot z(\theta T(e),e)-\dot z(\theta \Tmin,e)\big) + \Tmin\dot Z(\theta \Tmin, e) .
\end{equation*}
The boundedness of $\ddot z$, $T\in C^1(\bar I)$, and~\eqref{E:ZDUHAMEL} yield 
$$\big|\dot z(\theta T(e),e)-\dot z(\theta \Tmin,e)\big| + |\dot Z(\theta \Tmin, e)|\leq C\big(e-E_{\min})^{\frac12},$$
and so
\begin{equation*}
	\dot z(\theta T(e),e)=-\alpha\sin(2\pi\theta)+\mathcal O_{e\to\Emin}(\sqrt{e-\Emin}).
\end{equation*}
Inserting this expansion and~\eqref{E:PARTIALRMINUSEXPANSION} into~\eqref{E:WDIFFQUOTIENT} yields the first expansion term in~\eqref{E:HOLDERONEHALFW} upon letting $\tilde E\searrow\Emin$.
\end{proof}

%%%%%%%%%%%%%%%%%%%%%%%%%%%%%%%%%%%%%%%
%%%%%%%%%%%%%%%%%%%%%%%%%%%%%%%%%%%%%%%

\begin{remark}\label{R:HOLDERONEHALFDERIVATIVES}
	Notice that we have proven~\eqref{E:HOLDERONEHALFR} and~\eqref{E:HOLDERONEHALFW} by integrating the expansion formulae for the $E$-derivatives of~$r$ and $w$, respectively.
	This means that na\"ively differentiating~\eqref{E:HOLDERONEHALFR} and~\eqref{E:HOLDERONEHALFW} w.r.t.~$E$ indeed yields the correct expansions of $\partial_E^kr$ and $\partial_E^kw$, respectively, for $k\in\{1,2\}$.
\end{remark}

The previous lemma allows us to capture effective cancellations that will be important in the proof of the main theorem.

\begin{lemma}\label{L:IMPROVED0}
As $E\searrow\Emin$, the following bounds hold uniformly in $\theta,\tilde\theta\in\mathbb S^1$:
\begin{align}
w(\theta,E)w(\tilde\theta,E) &= 2\sin (2\pi\theta)\sin(2\pi\tilde\theta)\,(E-\Emin)\nonumber\\
&\quad-\frac{\sqrt2\Psi'''(r_\ast)}{3\alpha^3}\left(\sin(2\pi\theta)\sin(4\pi\tilde\theta)+\sin(4\pi\theta)\sin(2\pi\tilde\theta)\right)(E-\Emin)^{\frac32} \nonumber\\&\quad+ \mathcal O((E-\Emin)^2), \label{E:WWEXPANSION}\\
w(\theta,E)\Psi'(r(\tilde\theta,E)) &=  -2 \alpha \sin(2\pi\theta)\cos(2\pi\tilde\theta)\,(E-\Emin)\nonumber\\&\quad+\frac{\sqrt2\Psi'''(r_\ast)}{3\alpha^2}\left(2\sin(2\pi\theta)\cos(4\pi\tilde\theta)+\sin(4\pi\theta)\cos(2\pi\tilde\theta)\right)(E-\Emin)^{\frac32} \nonumber\\&\quad+ \mathcal O((E-\Emin)^2);\label{E:PSIPRIMEWEXPANSION}
\end{align} 
recall that $\alpha=\sqrt{\Psi''(r_\ast)}$.
Similar to Remark~\ref{R:HOLDERONEHALFDERIVATIVES}, na\"ively differentiating both sides of~\eqref{E:WWEXPANSION} and~\eqref{E:PSIPRIMEWEXPANSION} once or twice w.r.t.~$E$ results in the correct expansion formulae of $\partial_E^k(w(\theta,E)w(\tilde\theta,E))$ and $\partial_E^k(w(\theta,E)\Psi'(r(\tilde\theta,E)))$ for $k\in\{1,2\}$.
In particular, the functions $I\ni E\mapsto w(\theta,E)w(\tilde\theta,E)$ and $I\ni E\mapsto w(\theta,E)\Psi'(r(\tilde\theta,E))$ extend to continuously differentiable functions at $E=\Emin$ and their derivatives are  $C^{0,\frac12}$ on $\bar I$.
\end{lemma}

%%%%%%%%%%%%%%%%%%%%%%%%%%%%%%%%%%%%%%%

\begin{proof}
Both expansions follow directly from Lemma~\ref{L:HOLDERONEHALF}; the second one further relies on an expansion of $\Psi'$ around~$r_\ast$.
\end{proof}
%%%%%%%%%%%%%%%%%%%%%%%%%%%%%%%%%%%%%%%
%%%%%%%%%%%%%%%%%%%%%%%%%%%%%%%%%%%%%%%

\begin{lemma}\label{L:IMPROVED}
Let $f_0\in W^{3,\infty}(\Omega)$.
There exists a constant $C>0$ such that for any $m\in\Z^\ast$, $E\in I$, and $\theta\in\mathbb S^1$ the following bounds hold:
\begin{align}
\lv\partial_E\left(\fhat_0(m,E)w(\theta,E)\right)\rv &\le \frac {C\|f_0\|_{W^{2,\infty}(\Om)}}{|m|} , \label{E:IMPROVED1}\\
\lv\partial_E^2\left(\fhat_0(m,E)w(\theta,E)\right)\rv &\le \frac {C\|f_0\|_{W^{3,\infty}(\Om)}}{|m|} (E-\Emin)^{-\frac12}, \label{E:IMPROVED2}\\
\lv\partial_E\left(\fhat_0(m,E)\pa_Er(\theta,E)\right)\rv &\le \frac {C\|f_0\|_{W^{2,\infty}(\Om)}}{|m|} (E-\Emin)^{-\frac12}. \label{E:IMPROVED3}
\end{align}
\end{lemma}

%%%%%%%%%%%%%%%%%%%%%%%%%%%%%%%%%%%%%%%

\begin{proof}
We first integrate by parts to obtain
\begin{align}
	\fhat_0&(m,E)\,w(\theta,E)=w(\theta,E)\int_{\mathbb S^1}f_0(\tilde\theta,E)e^{-2\pi im\tilde\theta}\diff\tilde\theta=\frac{w(\theta,E)}{2\pi im}\int_{\mathbb S^1}\partial_\theta f_0(\tilde\theta,E)e^{-2\pi im\tilde\theta}\diff\tilde\theta\nonumber\\
	&=\frac{T(E)}{2\pi im}\int_{\mathbb S^1}\left(\partial_rf_0(\tilde\theta,E)\,w(\tilde\theta,E)-\partial_wf_0(\tilde\theta,E)\,\Psi'(r(\tilde\theta,E))\right)w(\theta,E)\,e^{-2\pi im\tilde\theta}\diff\tilde\theta.\label{eq:fhatw}
\end{align}
Our aim is study the rates of degeneracy as $E\searrow\Emin$ of all terms in the latter expression.

In order to derive the expansions of $\partial_rf_0$ and~$\partial_wf_0$, first note that for any $h\in W^{1,\infty}(\Omega)$,
\begin{multline*}
	\partial_Eh(\theta,E)=\partial_rh(\theta,E)\,\partial_Er(\theta,E)+\partial_wh(\theta,E)\,\partial_Ew(\theta,E)=\\
	=-\partial_rh(\theta,E)\frac{\cos(2\pi\theta)}{\sqrt2\alpha}(E-\Emin)^{-\frac12}+\partial_wh(\theta,E)\frac{\sin(2\pi\theta)}{\sqrt2}(E-\Emin)^{-\frac12}+\mathcal O(1)
\end{multline*}
uniformly in $\theta\in\mathbb S^1$ as $E\searrow\Emin$; recall Remark~\ref{R:HOLDERONEHALFDERIVATIVES}.
Integrating this bound w.r.t.~$E$ yields
\begin{equation*}
	h(\theta,E)=h(r_\ast,0)+\mathcal O(\sqrt{E-\Emin});
\end{equation*}
note that $h(\theta,\Emin)=h(r_\ast,0)$ for any $\theta\in\mathbb S^1$.
Iterating this argument once more leads to
\begin{multline*}
	g(\theta,E)=g(r_\ast,0)+\left(-\partial_rg(r_\ast,0)\frac{\sqrt2}\alpha\cos(2\pi\theta)+\partial_wg(r_\ast,0)\sqrt2\sin(2\pi\theta)\right)\sqrt{E-\Emin}\\
	+\mathcal O(E-\Emin)
\end{multline*}
for $g\in W^{2,\infty}(\Omega)$ uniformly in $\theta\in\mathbb S^1$ as $E\searrow\Emin$.
Inserting the latter bound with $g=\partial_rf_0$ and $g=\partial_wf_0$, Lemma~\ref{L:IMPROVED0}, and the obvious expansion $T(E)=\Tmin+\mathcal O(E-\Emin)$ into~\eqref{eq:fhatw} results in
\begin{equation}\label{eq:hatf0wexpansion}
	\fhat_0(m,E)\,w(\theta,E)=\frac{\Tmin}{2\pi im}(E-\Emin)\beta_1(\theta,m)+\frac{\Tmin}{2\pi im}(E-\Emin)^{\frac32}\beta_2(\theta,m)+\mathcal O(\frac{(E-\Emin)^2}{|m|})
\end{equation}
uniformly in $\theta\in\mathbb S^1$ as $E\searrow\Emin$.
Here, $\beta_1(\cdot,m),\beta_2(\cdot,m)\in C^\infty(\mathbb S^1;\mathbb C)$ are suitable combinations of trigonometric functions which are bounded by $C\|f_0\|_{W^{1,\infty}}$ and $C\|f_0\|_{W^{2,\infty}}$, respectively.

Due to Remark~\ref{R:HOLDERONEHALFDERIVATIVES}, na\"ively differentiating~\eqref{eq:hatf0wexpansion} yields the correct expansions of $\partial_E(\fhat_0(m,E)\,w(\theta,E))$ and $\partial_E^2(\fhat_0(m,E)\,w(\theta,E))$.
We thus conclude~\eqref{E:IMPROVED1} and~\eqref{E:IMPROVED2}.

The remaining bound~\eqref{E:IMPROVED3} follows similarly by employing the expansion of $\partial_Er(\theta,E)$ established in Remark~\ref{R:HOLDERONEHALFDERIVATIVES}. 
More precisely,
\begin{equation}\label{eq:hatf0parexpansion}
	\fhat_0(m,E)\,\partial_Er(\theta,E)=\frac{\Tmin}{2\pi im}\gamma_1(\theta,m)+\frac{\Tmin}{2\pi im}\sqrt{E-\Emin}\gamma_2(\theta,m)+\mathcal O(\frac{E-\Emin}{|m|})
\end{equation}
uniformly in $\theta\in\mathbb S^1$ as $E\searrow\Emin$, where $\gamma_1$ and $\gamma_2$ are functions with similar properties to~$\beta_1$ and~$\beta_2$.
Furthermore, the expansion of $\partial_E\left(\fhat_0(m,E)\pa_Er(\theta,E)\right)$ can again be obtained by na\"ively differentiating~\eqref{eq:hatf0parexpansion}, and we conclude~\eqref{E:IMPROVED3}. 
\end{proof}

We recall the definitions~\eqref{E:GRDEF} and~\eqref{E:HRDEF} of $g_R$ and $h_R$, respectively. 
The Fourier description of these functions in the action-angle variables plays a key role in the
understanding of the decay-in-$t$ of the gravitational potential.

\begin{lemma}[Action-angle descriptions of the supports of $g_R$ and $h_R$]\label{L:SUPPORTS}
The supports $\Omega_{g_R}$ and $\Omega_{h_R}$ of $g_R$ and $h_R$ in action-angle variables are respectively given by
\begin{align*}
  \Omega_{g_R}
& =
\begin{cases}
 \cup_{E\in [\Psi(R),E_0]}\big( [\theta(R,E),1-\theta(R,E)] \times \{E\} \big), & R>r_\ast, \\
 \cup_{E\in (\Psi(R),E_0]}\big( [\theta(R,E),1-\theta(R,E)] \times \{E\} \big)\cup \big(\mathbb S^1 \times [\Emin,\Psi(R)]\big), & R\leq r_\ast,
\end{cases} \\
 \Omega_{h_R}
&=
\begin{cases}
 \cup_{E\in (\Psi(R),E_0]}\big( [-\theta(R,E),\theta(R,E)] \times \{E\} \big)\cup \big(\mathbb S^1 \times [\Emin,\Psi(R)]\big), & R\ge r_\ast, \\
 \cup_{E\in [\Psi(R),E_0]}\big( [-\theta(R,E),\theta(R,E)] \times \{E\} \big), & R< r_\ast.
\end{cases}
\end{align*}
Moreover, $g_R(r(\theta,E)) =  \frac1{r(\theta,E)^2}\chi_{\Omega_{g_R}}(\theta,E)$ and $h_R(r(\theta,E)) =  \chi_{\Omega_{h_R}}(\theta,E)$.
\end{lemma}
\begin{proof}
	The characterisations of $\Omega_{g_R}$ and $\Omega_{h_R}$ are due to~\eqref{E:rthetaEgeqR_Rlowerrast} and~\eqref{E:rthetaEgeqR_Rlargerrast}.
\end{proof}

%%%%%%%%%%%%%%%%%%%%%%%%%%%%%%%%%
%%%%%%%%%%%%%%%%%%%%%%%%%%%%%%%%%

\begin{lemma}[Fourier descriptions of $wg_R$ and $h_R$]\label{L:GHFOURIER}
Let $R\in[\Rmin,\Rmax]$ be given. Then for any $m\in\mathbb Z_\ast$ we have the following formulae for the Fourier coefficients of $wg_R$ and $h_R$:
\begin{align}
\widehat{w g_R}(m,E) & = \begin{cases} \int_{\theta(R,E) }^{1-\theta(R,E) } \frac{w}{r^2} e^{-2\pi i m \theta}\diff\theta,&\text{if }E\in(\Psi(R),E_0],\\
\int_{\mathbb S^1} \frac{w}{r^2} e^{-2\pi i m \theta}\diff\theta,&\text{if }\Emin<E\leq\Psi(R)\,\land\,R<r_\ast,\\
0,&\text{if }\Emin<E\leq\Psi(R)\,\land\,R>r_\ast,
\end{cases}\label{E:WGREXPLICIT} \\
\widehat{h_R}(m,E) & = \begin{cases} \frac{1}{m\pi}\sin(2\pi m\theta(R,E)),  
&\text{if }E\in(\Psi(R),E_0],\\    
0 & \text{if }E\in[\Emin,\Psi(R)].
\end{cases}
\label{E:HREXPLICIT}
\end{align}
\end{lemma}

%%%%%%%%%%%%%%%%%%%%%%%%%%%%%%%%%

\begin{proof}
We use the formula
\begin{equation*}
\widehat{w g_R}(m,E)  = \int_{\mathbb S^1} \frac{w}{r^2}\chi_{\{r(\theta,E)\ge R\}} e^{-2\pi i m \theta}\diff\theta 
\end{equation*}
and Lemma~\ref{L:SUPPORTS} to obtain~\eqref{E:WGREXPLICIT}.
To show~\eqref{E:HREXPLICIT} we use the formula
\begin{equation*}
\widehat{h_R}(m,E) = \int_{-\theta(R,E)}^{\theta(R,E)}  e^{-2\pi i m\theta}\diff\theta,\quad m\in\Z_*,
\end{equation*}
where we note that in the region $r\ge R_\ast$ for any $E\in[\Emin,\Psi(R)]$, we have 
$\int_{\mathbb S^1} h_R(\theta, E)e^{-2\pi i m}\diff\theta=\int_{\mathbb S^1} e^{-2\pi i m}\diff\theta=0$.
\end{proof}

To prove quantitative decay of the solution, we need to carefully describe the regularity of the maps, defined for a given $f_0\colon\overline{\Om}\to\R$,
\begin{align}
G_{f_0}(m,\cdot)\colon I\to\mathbb C,\;G_{f_0}(m,E) &\coloneqq\fhat_0(m,E)\,\overline{\widehat{w g_R}(m,E)}, \label{E:HDEF} \\
H_{f_0}(m,\cdot)\colon I\to\mathbb C,\;H_{f_0}(m,E)& \coloneqq\fhat_0(m,E)\,\overline{\widehat{h_R}(m,E)}, \label{E:HDEFNEW}
\end{align}
for $m\in\Z^\ast$
near the elliptic point $E=\Emin$.
To that end we first write down the expressions for the first derivatives of $G_{f_0}(m,\cdot)$ and $H_{f_0}(m,\cdot)$ as well as the second derivative of $G_{f_0}(m,\cdot)$.

%%%%%%%%%%%%%%%%%%%%%%%%%%%%%%%
%%%%%%%%%%%%%%%%%%%%%%%%%%%%%%%

\begin{lemma}
Let $R\in [\Rmin,\Rmax]$ and $f_0\in W^{2,\infty}(\Omega)$ be given. Then for any $m\in\Z_\ast$ and $E\in(\Psi(R),E_0]$ the following identities hold:
\begin{align}
&\frac d{dE}G_{f_0}(m,E) \notag\\
& = \int_{\theta(R,E) }^{1-\theta(R,E) } \frac{d}{dE}\left(\frac{\fhat_0(m,E)w}{r^2}\right) e^{2\pi i m \theta}\diff\theta 
+ \frac{2 i\fhat_0(m,E)\pa_E r(\theta(R,E)  ,E)}{T(E)\,R^2} \sin(2\pi m \theta(R,E)  ), \label{E:EXTERIOR0} \\
&\frac d{dE}H_{f_0}(m,E) \notag\\
& = \frac1{\pi m} \pa_E \fhat_0(m,E)\sin(2\pi m \theta(R,E))+ 2 \fhat_0(m,E)\cos(2\pi m \theta(R,E))\partial_E\theta(R,E), \label{E:HDERIVATIVE}
\end{align}
\begin{align}
&\frac {d^2}{dE^2}G_{f_0}(m,E) \notag \\
&= \int_{\theta(R,E) }^{1-\theta(R,E) } \frac{d^2}{dE^2}\left(\frac{\fhat_0(m,E)w}{r^2}\right) e^{2\pi i m \theta}\diff\theta \notag\\
& \ \ \ \ -\pa_E\theta(R,E)   \left[\pa_E\left(\frac{\fhat_0(m,E)w}{r^2}\right)\Big|_{\theta= 1-\theta(R,E)  } e^{-2\pi i m \theta(R,E)  } 
+\pa_E\left(\frac{\fhat_0(m,E)w}{r^2}\right)\Big|_{\theta= \theta(R,E)  } e^{2\pi i m \theta(R,E)  }\right] \notag\\
& \ \ \ \ +\frac{2 i\pa_E\theta(R,E)  \fhat_0(m,E)}{R^2} \big(\pa_{E\theta}r(\theta(R,E)  ,E)\sin(2\pi m \theta(R,E)  )+2\pi m \pa_Er(\theta(R,E)  ,E)\cos(2\pi m \theta(R,E))\big) \notag\\
& \ \ \ \ + \frac{2 i}{R^2}\pa_{E} \left(\frac{\fhat_0(m,E) \pa_Er(s,E)}{T(E)}\right)\Big|_{s=\theta(R,E)  } \sin(2\pi m \theta(R,E)  ). \label{E:EXTERIOR1}
\end{align}
\end{lemma}

%%%%%%%%%%%%%%%%%%%%%%%%%%%%%%%

\begin{proof}
Differentiating $G_{f_0}(m,E)$ with respect to $E$ we get
\begin{align}
\frac d{dE}G_{f_0}(m,E)
& = \int_{\theta(R,E) }^{1-\theta(R,E) } \frac{d}{dE}\left(\frac{\fhat_0(m,E)w}{r^2}\right) e^{2\pi i m \theta}\diff\theta \notag\\
& \ \ \ \ - \frac{2i\fhat_0(m,E)w(\theta(R,E)  ,E)}{R^2}\pa_E\theta(R,E)  \sin(2\pi m \theta(R,E)  ). \label{E:EXTERIOR-1}
\end{align}
Upon plugging~\eqref{E:THETADEF2} into~\eqref{E:EXTERIOR-1}, we obtain the claimed identity~\eqref{E:EXTERIOR0}. Identity~\eqref{E:HDERIVATIVE}
is trivial; it follows from the product rule and the formula~\eqref{E:HREXPLICIT}.
Identity~\eqref{E:EXTERIOR1} follows by taking a further derivative.
\end{proof}

%%%%%%%%%%%%%%%%%%%%%%%%%%%%%%%
%%%%%%%%%%%%%%%%%%%%%%%%%%%%%%%

\begin{lemma}\label{L:HINTEGRABLE}
Let $R\in [\Rmin,\Rmax]$ and $f_0\in W^{3,\infty}(\Omega)$ be given.
Then there exists a constant $C>0$ such that 
\begin{align}
\lv \frac d{dE}G_{f_0}(m,E) \rv &\le \frac {C\|f_0\|_{W^{2,\infty}(\Om)}}{|m|},  \ \ &&E\in(\Psi(R),E_0], \label{E:PARTIALEHBOUND}\\
\lv \frac{d^2}{dE^2}G_{f_0}(m,E) \rv &\le \frac {C\|f_0\|_{W^{3,\infty}(\Om)}}{|m|} \frac1{(E-\Emin)^{\frac12}(E-\Psi(R))^{\frac12}}, \ \  &&E\in(\Psi(R),E_0]. \label{E:PARTIALEEHBOUND}
\end{align}
If $R\in[\Rmin,r_\ast)$ then 
\begin{align}
\lv \frac d{dE}G_{f_0}(m,E) \rv &\le \frac {C\|f_0\|_{W^{2,\infty}(\Om)}}{|m|},  \ \ &&E\in(\Emin,\Psi(R)], \label{E:PARTIALEHBOUNDLOW}\\
\lv \frac{d^2}{dE^2}G_{f_0}(m,E) \rv &\le  \frac {C\|f_0\|_{W^{3,\infty}(\Om)}}{|m|} (E-\Emin)^{-\frac12}, \ \  &&E\in(\Emin,\Psi(R)]. \label{E:PARTIALEEHBOUNDLOW}
\end{align}
Moreover, for any $R\in [\Rmin,\Rmax]$ we have the bound
\begin{align}\label{E:PARTIALEHNEWBOUND}
\lv \frac d{dE}H_{f_0}(m,E) \rv \le \frac{C\|f_0\|_{W^{2,\infty}(\Om)}}{|m|} (E-\Psi(R))^{-\frac12},\qquad E\in(\Psi(R),E_0]
\end{align}
\end{lemma}

%%%%%%%%%%%%%%%%%%%%%%%%%%%%%%%

\begin{proof}
The bound~\eqref{E:PARTIALEHBOUND} follows by plugging the expansions~\eqref{eq:hatf0wexpansion} and~\eqref{eq:hatf0parexpansion} established in the proof of Lemma~\ref{L:IMPROVED} as well as the expansion of~$\partial_Er$ given by Remark~\ref{R:HOLDERONEHALFDERIVATIVES} into~\eqref{E:EXTERIOR0}. 

In the same way, 
we can bound the first line on the right-hand side of~\eqref{E:EXTERIOR1} by 
\begin{align}
	\lv \int_{\theta(R,E) }^{1-\theta(R,E) } \frac{d^2}{dE^2}\left(\frac{\fhat_0(m,E)w}{r^2}\right) e^{2\pi i m \theta}\diff\theta  \rv
	\le \frac {C\|f_0\|_{W^{3,\infty}(\Om)}}{|m|} (E-\Emin)^{-\frac12}. \label{E:B1}
\end{align}
By additionally using Lemma~\ref{L:THETAREG},
we can bound the second line on the right-hand side of~\eqref{E:EXTERIOR1} by
\begin{align*}
&\lv\pa_E\theta(R,E)   \left[\frac{d}{dE}\left(\frac{\fhat_0(m,E)w}{r^2}\right)\Big|_{\theta= 1-\theta(R,E)  } e^{-2\pi i m \theta(R,E)} 
+\pa_E\left(\frac{\fhat_0(m,E)w}{r^2}\right)\Big|_{\theta= \theta(R,E)} e^{2\pi i m \theta(R,E)}\right]  \rv \notag\\
& \le \frac {C\|f_0\|_{W^{2,\infty}(\Om)}}{|m|} (E-\Psi(R))^{-\frac12}  (E-\Emin)^{-\frac12} .
\end{align*}
Similarly, by using Lemmas~\ref{lemma:Eminvanishing} and~\ref{L:THETAREG} as well as~\eqref{E:PARTIALER}, 
we conclude that the third line on the right-hand side of~\eqref{E:EXTERIOR1} is bounded by
\begin{align*}
& \lv\frac{2\pa_E\theta(R,E)\fhat_0(m,E)}{R^2} \left(\pa_{E\theta}r(\theta(R,E),E)\sin(2\pi m \theta(R,E))+2\pi m \partial_Er(\theta(R,E),E)\cos(2\pi m \theta(R,E))\right) \rv \notag\\
& \le \frac {C\|f_0\|_{W^{2,\infty}(\Om)}}{|m|} (E-\Psi(R))^{-\frac12}(E-\Emin)^{-\frac12};
\end{align*}
to estimate the second term, we used the case $k=2$ in Lemma~\ref{lemma:Eminvanishing}.

Finally, using~\eqref{E:IMPROVED3} we have the bound
\begin{align}
\lv \frac{2}{R^2}\pa_{E} \left(\frac{\fhat_0(m,E) \pa_Er(s,E)}{T(E)}\right)\Big|_{s=\theta(R,E)} \sin(2\pi m \theta(R,E)) \rv
\le  \frac {C\|f_0\|_{W^{2,\infty}(\Om)}}{|m|} (E-\Emin)^{-\frac12}. \label{E:B4}
\end{align}
Summing~\eqref{E:B1}--\eqref{E:B4} and using~\eqref{E:EXTERIOR1}, we conclude~\eqref{E:PARTIALEEHBOUND}.

For $R\in[\Rmin,r_\ast)$ we have
\begin{equation*}
G_{f_0}(m,E) =  \int_{\mathbb S^1} \frac{\fhat_0(m,E)\,w}{r^2} e^{2\pi i m \theta}\diff\theta,\quad  E\in(\Emin,\Psi(R)],
\end{equation*}
by~\eqref{E:WGREXPLICIT},
so~\eqref{E:PARTIALEHBOUNDLOW} and~\eqref{E:PARTIALEEHBOUNDLOW} follow in the same way as~\eqref{E:PARTIALEHBOUND} and~\eqref{E:B1}, respectively. 
It remains to show~\eqref{E:PARTIALEHNEWBOUND}. This follows from Lemmas~\ref{lemma:Eminvanishing} and~\ref{L:THETAREG} together with the bound
\begin{equation}\label{eq:paEhatf}
\lv\pa_E\fhat_0(m,E)\rv \le\frac{C\|f_0\|_{W^{2,\infty}}}{|m|\sqrt{E-\Emin}},
\end{equation}
which can be obtained similarly to Lemma~\ref{lemma:Eminvanishing} using~\eqref{E:PARTIALER}.
\end{proof}

In order to obtain the fractional gain in decay as stated in~\eqref{E:PARTIALTUFDECAYEXTERIORINTRO} and~\eqref{E:PARTIALRUFDECAYEXTRAMAIN}, we shall need
a precise control on even higher-order $E$-derivatives of functions $G(m,\cdot)$ and $H(m,\cdot)$.

\begin{lemma}[Higher-order regularity]\label{L:HIGHORDER}
Let $f_0\in W^{3,\infty}(\Omega)$.
Then there exists a constant $C>0$ such that for any $m\in\Z_\ast$, $R\in[\Rmin,\Rmax]$, and $E\in(\Psi(R),E_0]$ the following bounds hold:
\begin{align}
\lv \pa_E\left((E-\Psi(R))^{\frac12} \pa_{E}^2 G_{f_0}(m,E)\right)\rv&\leq C\|f_0\|_{W^{3,\infty}(\Om)} (E-\Emin)^{-1} (E-\Psi(R))^{-\frac12},  \label{E:PARTIALEEEG}\\
\lv \pa_E\left((E-\Psi(R))^{\frac12} \pa_E H_{f_0}(m,E)\right)\rv&\leq C\|f_0\|_{W^{2,\infty}(\Om)} (E-\Emin)^{-\frac12} (E-\Psi(R))^{-\frac12}.\label{E:PARTIALEEH}
\end{align}
\end{lemma}
\begin{proof}
	Let us start by proving the easier bound~\eqref{E:PARTIALEEH}.
	The relations $\sqrt{E-\Psi(R)}=\frac1{\sqrt2}w(\theta(R,E),E)$ and~\eqref{E:THETADEF2} as well as~\eqref{E:HDERIVATIVE} yield
	\begin{multline*}
		(E-\Psi(R))^{\frac12} \partial_E H_{f_0}(m,E)=\frac1{\sqrt2\pi m} w(\theta(R,E),E)\pa_E \fhat_0(m,E)\sin(2\pi m \theta(R,E))\\- \frac{\sqrt2}{T(E)} \partial_Er(\theta(R,E),E)\fhat_0(m,E)\cos(2\pi m \theta(R,E)).
	\end{multline*}
	Estimating the $E$-derivative of this expression is based on the boundedness of~$T$ and~$T'$, 
	the bound on~$\partial_E\theta(R,E)$ from Lemma~\ref{L:THETAREG},
	the bounds on $w(\theta,E)$, $\Psi'(R)=\Psi'(r(\theta(R,E),E))$, $\partial_Er(\theta,E)$, and $\partial_Ew(\theta,E)$ from Lemma~\ref{L:HOLDERONEHALF} and Remark~\ref{R:HOLDERONEHALFDERIVATIVES}, 
	the bound on $\partial_E(\fhat_0(m,E)\pa_Er(\theta,E))$ from Lemma~\ref{L:IMPROVED},
	as well as the bounds on $\fhat_0(m,E)$ and $\partial_E\fhat_0(m,E)$ from Lemma~\ref{lemma:Eminvanishing} and~\eqref{eq:paEhatf}.
	Together with the additional estimate
	\begin{equation*}
		\lv\partial_E\left(\partial_E\fhat_0(m,E)w(\theta,E)\right)\rv\leq C\|f_0\|_{W^{2,\infty}(\Om)} (E-\Emin)^{-\frac12},\qquad\theta\in\mathbb S^1,
	\end{equation*}
	which can be obtained similarly to Lemma~\ref{L:IMPROVED}, we conclude~\eqref{E:PARTIALEEH}.
	
	In order to prove~\eqref{E:PARTIALEEEG}, we first rewrite the factor $(E-\Psi(R))^{\frac12}$ in the same way as above and use~\eqref{E:EXTERIOR1} to arrive at
	\begin{align*}
		(&E-\Psi(R))^{\frac12} \pa_{E}^2 G_{f_0}(m,E)=\frac1{\sqrt2} \int_{\theta(R,E) }^{1-\theta(R,E)} \partial_E^2\left(\frac{\fhat_0(m,E)w(\theta,E)}{r(\theta,E)^2}\right)w(\theta(R,E),E)e^{2\pi i m \theta}\diff\theta \nonumber\\
		&\qquad +\frac{\pa_Er(\theta(R,E),E)}{\sqrt2T(E)} \Bigg[\partial_E\left(\frac{\fhat_0(m,E)w}{r^2}\right)\Big|_{\theta= 1-\theta(R,E)  } e^{-2\pi i m \theta(R,E)  } 
		\nonumber\\&\qquad\qquad\qquad\qquad\qquad\qquad+\pa_E\left(\frac{\fhat_0(m,E)w}{r^2}\right)\Big|_{\theta= \theta(R,E)  } e^{2\pi i m \theta(R,E)  }\Bigg] \nonumber\\
		&\qquad-\frac{\sqrt2 i\pa_Er(\theta(R,E),E)\fhat_0(m,E)}{T(E)\,R^2} \Big(\pa_{E\theta}r(\theta(R,E)  ,E)\sin(2\pi m \theta(R,E)  )\nonumber\\&\qquad\qquad\qquad\qquad\qquad\qquad\qquad\qquad\quad+2\pi m \pa_Er(\theta(R,E)  ,E)\cos(2\pi m \theta(R,E))\Big) \notag\\
		&\qquad + \frac{\sqrt2 i}{R^2}w(\theta(R,E),E)\pa_{E}\left(\frac{\fhat_0(m,E) \pa_Er(\theta,E)}{T(E)}\right)\Big|_{\theta=\theta(R,E)  } \sin(2\pi m \theta(R,E)  )\nonumber\\
		&\eqqcolon I+II+III+IV.
	\end{align*}
	Estimating the $E$-derivative of this expression requires some additional bounds compared to the above.
	
	Concretely, to estimate $\partial_EI$, first observe that Remark~\ref{R:HOLDERONEHALFDERIVATIVES} and an expansion formula similar to~\eqref{eq:hatf0wexpansion} from the proof of Lemma~\ref{L:IMPROVED} yield
	\begin{equation}\label{eq:pa2hatfwoverr}
		\partial_E^2\left(\frac{\fhat_0(m,E)w(\theta,E)}{r(\theta,E)^2}\right)=\tilde\beta_1(\theta,m)(E-\Emin)^{-\frac12}+\tilde\beta_2(\theta,m)+\mathcal O({\sqrt{E-\Emin}})
	\end{equation}
	as $E\searrow\Emin$ uniformly in $\theta\in\mathbb S^1$, where $\tilde\beta_1(\cdot,m),\tilde\beta_2(\cdot,m)\in C^{\infty}(\mathbb S^1)$ are suitable combinations of trigonometric functions which can be bounded uniformly in~$\|f_0\|_{W^{2,\infty}}$ and $\|f_0\|_{W^{3,\infty}}$, respectively.
	This expansion obviously yields the bound
	\begin{equation*}
		\lv\partial_E^2\left(\frac{\fhat_0(m,E)w(\theta,E)}{r(\theta,E)^2}\right)\rv\leq C\,\|f_0\|_{W^{2,\infty}}(E-\Emin)^{-\frac12}.
	\end{equation*}
	Furthermore, combined with the expansion of~$w$ from Lemma~\ref{L:HOLDERONEHALF}, \eqref{eq:pa2hatfwoverr} indicates
	\begin{equation}\label{eq:pa3hatfwoverr}
		\lv\partial_E\left[w(\tilde\theta,E)\,\partial_E^2\left(\frac{\fhat_0(m,E)w(\theta,E)}{r(\theta,E)^2}\right)\right]\rv\leq C\,\|f_0\|_{W^{3,\infty}}(E-\Emin)^{-\frac12}.
	\end{equation}
	The estimate~\eqref{eq:pa3hatfwoverr} can indeed be proven by extending the results from Remark~\ref{R:HOLDERONEHALFDERIVATIVES} from second-order to third-order $E$-derivatives by iterating the arguments from the proof of Lemma~\ref{L:HOLDERONEHALF} once more. 
	In the same way as shown in Lemma~\ref{L:IMPROVED}, this allows to derive an analogue of~\eqref{eq:pa2hatfwoverr} for the third-order $E$-derivative, which then leads to~\eqref{eq:pa3hatfwoverr}.
	Combining all these bounds yields
	\begin{equation}\label{eq:pa3estI}
		\lv\partial_EI\rv\leq C\,\|f_0\|_{W^{3,\infty}} (E-\Emin)^{-\frac12}(E-\Psi(R))^{-\frac12}.
	\end{equation}
	Similarly to~\eqref{E:IMPROVED1} we derive the following estimate for $\partial_EII$:
	\begin{equation}\label{eq:pa3estII}
		\lv\partial_EII\rv\leq C\,\|f_0\|_{W^{2,\infty}} (E-\Emin)^{-1}(E-\Psi(R))^{-\frac12}.
	\end{equation}
	Similarly, by additionally applying~\eqref{E:IMPROVED3} and the case $k=2$ from Lemma~\ref{lemma:Eminvanishing}, it follows that 
	\begin{equation}\label{eq:pa3estIII}
		\lv\partial_EIII\rv\leq C\,\|f_0\|_{W^{3,\infty}} (E-\Emin)^{-1}(E-\Psi(R))^{-\frac12}.
	\end{equation}
	In order to estimate $\partial_EIV$, observe that generalizing the arguments from Lemma~\ref{L:IMPROVED} similarly to~\eqref{eq:pa3hatfwoverr} yields
	\begin{equation*}
		\lv\partial_E\left[w(\tilde\theta,E)\,\pa_{E}\left(\frac{\fhat_0(m,E) \pa_Er(\theta,E)}{T(E)}\right)\right]\rv\leq\frac C{|m|}\,\|f_0\|_{W^{3,\infty}}(E-\Emin)^{-\frac12},
	\end{equation*}
	which leads to
	\begin{equation}\label{eq:pa3estIV}
		\lv\partial_EIV\rv\leq C\,\|f_0\|_{W^{3,\infty}}(E-\Emin)^{-1}(E-\Psi(R))^{-\frac12}.
	\end{equation}
	We thus conclude~\eqref{E:PARTIALEEEG} by combining~\eqref{eq:pa3estI}--\eqref{eq:pa3estIV}.
\end{proof}

%%%%%%%%%%%%%%%%%%%%%%%%%%%%%%%
%%%%%%%%%%%%%%%%%%%%%%%%%%%%%%%

\subsection{The radially symmetric case}\label{S:AAREG3D}

We now analyse the regularity of the action-angle variables in the general radially symmetric case.
To that end, we suppose henceforth that we are working with a steady state 
of the form~\eqref{E:POLY2}  with polytropic exponents $k$, $\ell$ and smallness parameter $\varepsilon>0$ and we introduce some further notation. For a sufficiently smooth function $F:J\to \mathbb C$ and for any $j\in\mathbb N$ we denote 
by 
\[
\pa^{(j)}F
\]
the collection of all derivatives of the form $\pa_E^n\pa_L^{j-n}F$,
$n\in\{0,\dots, j\}$.

%%%%%%%%%%%%%%%%%%%%%%%%%%%%%%%
%%%%%%%%%%%%%%%%%%%%%%%%%%%%%%%

\begin{lemma}\label{L:RELCONTROL}
	Let $k>1$, $\ell>0$. Then for any $j\in\N\cap(0,k+\ell+\frac52)$ there exists a constant $C>0$ such that 
	\begin{equation}\label{E:RCONTROL}
		\Big| \pa^{(j)}r(\theta,E,L)\Big| \le C (E-\EminL)^{\frac12-j}, \ \ \ \theta\in\mathbb S^1, \ (E,L)\in J.
	\end{equation}
	Moreover, $r(\theta,E,L)$ admits the following expansion. Let $\al_L:=\sqrt{\Psi_L''(r_L)}$ for $L>0$. Then 
	\beqa\label{E:RELEXP}
	r(\th,E,L)-r_L=-\frac{\sqrt{2}\cos(2\pi\th)}{\al_L}\sqrt{E-\EminL }+\mathcal O(E-\EminL )
	\eeqa
	as $E\to\EminL$ locally uniformly in $L>0$.
	Moreover, if we set
	\beqa\label{eq:paELr0}
	(\pa_E r)_0=-\frac{\cos(2\pi\th)}{\sqrt{2}\al_L\sqrt{E-\EminL }},\qquad (\pa_Lr)_0=\frac{\cos(2\pi\th)}{2\sqrt{2}r_L^2\al_L\sqrt{E-\EminL }},
	\eeqa 
	the derivatives satisfy
	\begin{align}
		\pa_E r(\th,E,L)=(\pa_E r)_0+\widetilde{\pa_E r},\label{ineq:paEr}\\
		\pa_L r(\th,E,L)=(\pa_L r)_0+\widetilde{\pa_L r},\label{ineq:paLr}
	\end{align}
	where
	\beq\label{ineq:paELrtildebounds}
	\Big| \pa^{(j)}\widetilde{\pa_E r}(\theta,E,L)\Big| +\Big| \pa^{(j)}\widetilde{\pa_L r}(\theta,E,L)\Big|  \le C (E-\EminL)^{\frac12-j}, \ \ \ j\in\N_0\text{ with }j<k+\ell+\frac32
	\eeq
	for $\theta\in\mathbb S^1$ and $(E,L)\in J$.
	
	Finally, $T,\om\in C^\infty(J)\cap C^j(\bar J)$  for all  $j\in\N\cap(0,k+\ell+\frac52)$.
\end{lemma}
\begin{proof}
	It has already been shown in~\cite[Rem.~A.3.5]{St23} that $r(\theta,\cdot,\cdot)\in C^j(\bar J)$ for any $j\in\N\cap(0,k+\ell+\frac52)$. Similarly, the claim~\eqref{E:RCONTROL}
	in the case $j=1$ has been shown in~\cite[Lemma 6.3.2]{St23}. It follows by differentiating the characteristic flow with respect to $E$ and $L$, which promptly implies that the leading order singular behaviour near the elliptic point $\EminL$ is described by the $\pa_Er_-(E,L)$ and $\pa_Lr_-(E,L)$ near $(\EminL,L)$. This is turn is easily seen to be proportional to the inverse of $\Psi_L'(r)$ near $r=r_L$, which gives~\eqref{E:RCONTROL} with $j=1$. The bounds on higher derivatives now follow by similar arguments, differentiating the characteristic flow. The limitation $j<k+\ell+\frac52$ comes simply from the regularity 
	of the potential $\Psi_L$, which at the boundary of $[\Rmin,\Rmax]$ is in $C^j$ for any $j\in \N\cap(0,k+\ell+\frac72)$, see~\cite{St23}.
	
	The expansions~\eqref{E:RELEXP}--\eqref{ineq:paLr} follow similarly to Lemma~\ref{L:HOLDERONEHALF}; see~\cite[Sc.~A.3]{Ku21} for a detailed discussion of the properties of~$\partial_LR$. We observe that, close to a trapping point $r_L$, the potential $\Psi_L$ is actually smooth, and hence the limitation on the regularity occurs only at the boundary points, $\{R_{\min},R_{\max}\}$, leading to the claimed restriction on the number of derivatives. The regularity statement for the period function $T$ (and therefore the frequency function $\om$) now follow by arguments analogous to~\cite[Ch.~3]{Ku21} and~\cite[App.~A]{St23}.
\end{proof}

%%%%%%%%%%%%%%%%%%%%%%%
%%%%%%%%%%%%%%%%%%%%%%%

We introduce the change of variables
\beq
(E,L)\mapsto(q(E,L),L),\qquad q(E,L)=\om(E,L).\nonumber
\eeq
Then
\beq
\frac{\pa (q,L)}{\pa (E,L)}=\begin{pmatrix}
	\om' & \pa_L\om \\ 0 & 1
\end{pmatrix}, \qquad \frac{\pa (E,L)}{\pa (q,L)}=\begin{pmatrix}
	(\om')^{-1} & -\frac{\pa_L\om}{\om'} \\ 0 & 1
\end{pmatrix},\nonumber
\eeq
where we recall the monotonicity assumption~\eqref{E:MONOTONEPERIODL}.
As in~\eqref{E:ASC}, the symbol~$'$ always denotes the $E$-derivative of a function depending on~$(E,L)$. 
Further note that $\omega'=O(1)$ and $\partial_L\omega=O(\varepsilon)$ as $\varepsilon\to0$ by~\cite[Sc.~6.2.1]{St23}. 
Given $R\in[R_{\min},R_{\max}]$, we introduce the notations
\beq
\om_R(L)=\om(\Psi_L(R),L),\qquad \om_*(L)=\om(E^{\min}_L,L),\qquad \om_0(L)=\om(E_0,L),\nonumber
\eeq
as well as ${\om_*^{\min}}\coloneqq\inf_J\omega$ and ${\om_0^{\max}}\coloneqq\sup_J\omega$. 

%%%%%%%%%%%%%%%%%%%%%%%
%%%%%%%%%%%%%%%%%%%%%%%

The next lemma enables us to exchange $q$ derivatives of $\widehat{h_R}(m,E(q,L),L)$ for an $L$ derivative at the cost of a less singular remainder term. This will be crucial in using $L$ regularity to obtain enhanced decay.

%%%%%%%%%%%%%%%%%%%%%%%
%%%%%%%%%%%%%%%%%%%%%%%

\begin{lemma}\label{L:PAQHR}
	Let $R\in [R_{\min},R_{\max}]$. Then there exist functions $Q_i(E,L)$, $Q_i\in C^2(J)$, $i=1,2$, bounded on $\bar J$ and satisfying 
	\beq\label{ineq:dQj}
	\big|\pa^{(j)}Q_1(E,L)\big| + (E-\EminL)^{j-\frac12}\big|\pa^{(j)} Q_2(E,L)\big|\leq C, \ \  (E,L)\in J,
	\eeq
	for $1\leq j<k+\ell+\frac12$ such that
	\begin{multline*}
		\pa_q\big(\widehat{h_R}(m,E(q,L),L)\big)=\,Q_1(E(q,L),L)\pa_L\big(\widehat{h_R}(m,E(q,L),L)\big)\\
		+Q_2(E(q,L),L)\cos(2\pi m\theta(R,E(q,L),L))\pa_L\big(w(R,E(q,L),L)\big)
	\end{multline*}
	for $(q,L)$ with $(E(q,L),L)\in J$ and $E(q,L)\in(\Psi_L(R),E_0]$, i.e., $q<\omega_R(L)$.
\end{lemma}

%%%%%%%%%%%%%%%%%%%%%%%
%%%%%%%%%%%%%%%%%%%%%%%

\begin{proof}
	We begin by expanding the $E$ and $L$ derivatives of $\theta(R,E,L)$. Noting the simple identities
	\begin{align}
		\pa_E\th(R,E,L)&=-\frac{\pa_E r(\th(R,E,L),E,L)}{T(E,L) w(R,E,L)},\label{E:PAER}\\
		\pa_L\th(R,E,L)&=-\frac{\pa_L r(\th(R,E,L),E,L)}{T(E,L) w(R,E,L)},\label{E:PALR}
	\end{align}
	we then apply~\eqref{eq:paELr0}--\eqref{ineq:paLr} to decompose the derivatives $\pa_E r$ and $\pa_Lr$ and observe
	\beqa\label{ineq:dELrexpansions}
	{}&\Big|\pa^{(j)}\Big(\frac{(\pa_E r)_0}{(\pa_L r)_0}\Big)\Big|\leq C, \ \ \ \ j=0,1,2.
	\eeqa
	Thus
	\beq
	\pa_E\th(R,E,L)=-\frac{\pa_E r}{(\pa_Lr)_0}\frac{(\pa_Lr)_0}{Tw}=\frac{(\pa_E r)_0}{(\pa_Lr)_0}\big(\pa_L\th + \frac{\widetilde{\pa_L r}}{Tw}\big)-\frac{\widetilde{\pa_E r}}{Tw}.\nonumber
	\eeq
	Now observe that
	\begin{multline*}
		\pa_L\big(\widehat{h_R}(m,E(q,L),L)\big)=2\cos(2\pi m\theta(R,E(q,L),L))\,\partial_L\left(\theta(R,E(q,L),L)\right)\\
		=2\cos(2\pi m \theta(R,E(q,L),L))\Big((\pa_E\th)(R,E(q,L),L)\big(-\frac{\pa_L\om}{\om'}\big)+(\pa_L\th)(R,E(q,L),L)\Big).
	\end{multline*}
	On the other hand, 
	\begin{equation*}
		\pa_q\big(\widehat{h_R}(m,E(q,L),L)\big)=2\cos(2\pi m \theta(R,E(q,L),L))(\pa_E\th)(R,E(q,L),L)\frac{1}{\om'}.
	\end{equation*}
	Now, setting
	$$\mu(E,L):=-\frac{(\pa_E r)_0}{(\pa_Lr)_0}\frac{\pa_L\om}{\om'}\Big(1-\frac{(\pa_E r)_0}{(\pa_Lr)_0}\frac{\pa_L\om}{\om'}\Big)^{-1}=O(\k),$$
	where we note that $\mu$ is bounded in $C^2(\bar J)$ by~\eqref{ineq:dELrexpansions}, we expand
	\begin{multline*}
		(\pa_E\th)(R,E(q,L),L)\frac{1}{\om'}=\mu\,\pa_E\th\,\frac{1}{\om'}+\frac{1-\mu}{\omega'}\left[\frac{(\pa_E r)_0}{(\pa_Lr)_0}\big(\pa_L\th + \frac{\widetilde{\pa_L r}}{Tw}\big)-\frac{\widetilde{\pa_E r}}{Tw}\right]\\
		=\frac{(1-\mu)}{\om'}\frac{(\pa_E r)_0}{(\pa_Lr)_0}\Big((\pa_E\th)\big(-\frac{\pa_L\om}{\om'}\big)+(\pa_L\th)\Big)+\frac{1-\mu}{\omega'}\left[\frac{(\pa_E r)_0}{(\pa_Lr)_0}\,\frac{\widetilde{\pa_L r}}{Tw}-\frac{\widetilde{\pa_E r}}{Tw}\right].
	\end{multline*}
	Next, we note that $w(R,E(q,L),L)=\sqrt{2(E(q,L)-\Psi_L(R))}$, so that
	\beq\label{eq:dLwidentity}
	\pa_L\big(w(R,E(q,L),L)\big)=\frac{\pa_L(E(q,L))-\pa_L(\Psi_L(R))}{w}=-\frac{1}{w}\Big(\frac{\pa_L\om}{\om'}+\frac{1}{2R^2}\Big).
	\eeq
	Again using that $\om_L=O(\eps)$, the factor on the right hand side is uniformly bounded away from zero. We therefore re-write 
	\begin{multline*}
		\frac{1-\mu}{\omega'}\left[\frac{(\pa_E r)_0}{(\pa_Lr)_0}\,\frac{\widetilde{\pa_L r}}{T}-\frac{\widetilde{\pa_E r}}{T}\right]\frac1w\\
		=-\Big(\frac{\pa_L\om}{\om'}+\frac{1}{2R^2}\Big)^{-1}\frac{1-\mu}{\omega'}\left[\frac{(\pa_E r)_0}{(\pa_Lr)_0}\,\frac{\widetilde{\pa_L r}}{T}-\frac{\widetilde{\pa_E r}}{T}\right]\pa_L\big(w(R,E(q,L),L)\big).
	\end{multline*}
	Now setting
	\beqa
	Q_1(E,L)=&\,\frac{(1-\mu)}{\om'}\frac{(\pa_E r)_0}{(\pa_Lr)_0},\\
	Q_2(E,L)=&\,\Big(\frac{\pa_L\om}{\om'}+\frac{1}{2R^2}\Big)^{-1}\frac{1-\mu}{\omega'}\left[\frac{(\pa_E r)_0}{(\pa_Lr)_0}\,\frac{\widetilde{\pa_L r}}{T}-\frac{\widetilde{\pa_E r}}{T}\right],\nonumber
	\eeqa
	we have obtained
	\begin{align*}
		\pa_q&\big(\widehat{h_R}(m,E(q,L),L)\big)
		=2\cos(2\pi m \theta(R,E(q,L),L))(\pa_E\th)(R,E(q,L),L)\frac{1}{\om'}\\
		=&\,2\cos(2\pi m \theta(R,E(q,L),L))\Big(Q_1\pa_L\big(\theta(R,E(q,L),L)\big)+ Q_2\pa_L\big(w(R,E(q,L),L)\big)\Big)\\
		=&\,Q_1\pa_L\big(\widehat{h_R}(m,E(q,L),L)\big)+2\cos(2\pi m \theta(R,E(q,L),L)) Q_2\pa_L\big(w(R,E(q,L),L)\big),
	\end{align*}
	as claimed, where the bounds on $Q_1$ and $Q_2$ come from~\eqref{ineq:dELrexpansions},~\eqref{ineq:paELrtildebounds}, and the smoothness of $\om$.
\end{proof}

%%%%%%%%%%%%%%%%%%%%%%%
%%%%%%%%%%%%%%%%%%%%%%%

\begin{lemma}[Integral identities]\label{L:RECURSION1}
	Let $F:J\to\mathbb C$ be a $C^1$-function that vanishes on~$\partial J_{{\text{vac}}}$, see~\eqref{E:JVACDEF}.
	Then the following identities hold
	\begin{align}
		&\int_{L_0}^{L_{\max}} \int_{\om_R(L)}^{\om_0(L)} e^{-2\pi mi q t}  \widehat{h_R}(m,q,L) F \diff q\,\diff L \notag \\
		& = -\frac{1}{2\pi i m t}\int_{L_0}^{L_{\max}} \int_{\om_R(L)}^{\om_0(L)}e^{-2\pi mi q t} 
		\Bigg[ \widehat h_R\A_1 F+ w\cos(2\pi m\theta) \A_2 F
		\Bigg] \diff q\diff L \label{E:INTID1}
	\end{align}
	where
	\beqa\label{DEF:AJOPS}
	\A_1 F & : = \pa_L\Big(Q_1 F\Big)  +2\pi^2 m^2\big(\frac{\pa_E r}{\om'}+\pa_Lr\big)\frac{Q_2F}{T} -\pa_qF, \\
	\A_2 F & := \pa_L\Big(Q_2F\Big).
	\eeqa
	Moreover, 
	\begin{align}
		&\int_{L_0}^{L_{\max}} \int_{\om_R(L)}^{\om_0(L)} e^{-2\pi mi q t} w\cos(2\pi m\th) F \diff q\diff L \notag\\
		&= - \frac1{2\pi i mt} \int_{L_0}^{L_{\max}} \int_{\om_R(L)}^{\om_0(L)} e^{-2\pi mi q t}\Bigg[
		\widehat{h_R} \B_1 F +w\cos(2\pi m\theta) \B_2 F\Bigg] \diff q\diff L,\label{E:INTID2}
	\end{align}
	where
	\beqa\label{DEF:BJOPS}
	\B_1 F & : =-2\pi^2 m^2\frac{\pa_E r}{T\om'}F-2\pi^2m^2\big(\frac{\pa_E r}{\om'}+\pa_Lr\big)\frac F{\om' \big(\frac{\pa_L\om}{\om'}+\frac{1}{2R^2}\big) T}, \\
	\B_2 F & := -\pa_q F-\pa_L\left(\frac{F}{\om'}\big(\frac{\pa_L\om}{\om'}+\frac{1}{2R^2}\big)^{-1}\right).
	\eeqa
\end{lemma}

%%%%%%%%%%%%%%%%%%%%%%%

\begin{proof}
	Integrating-by-parts we get
	\begin{equation*}
		\int_{L_0}^{L_{\max}} \int_{\om_R(L)}^{\om_0(L)}e^{-2\pi mi q t}\widehat{h_R} F\,\diff q\,\diff L
		=\frac{1}{2\pi i m t}\int_{L_0}^{L_{\max}} \int_{\om_R(L)}^{\om_0(L)}e^{-2\pi mi q t}\pa_q\Big(\widehat{h_R} F\Big)\,\diff q\,\diff L,
	\end{equation*}
	where we have used that~$F$ vanishes on $\partial J_{{\text{vac}}}$
	to eliminate the boundary terms in the integration by parts due to $F(E(\om_0(L),L),L)=F(E_0,L)=0$.
	We also used 
	\beq\label{E:TPVANISH}
	\sin(2\pi m\theta(R,E(\om_R(L),L),L))=\sin(2\pi m\theta(R,\Psi_L(R),L))=0.
	\eeq
	Distributing derivatives, recalling that $\widehat{h_R}(m,E,L)$ is real-valued, and applying Lemma~\ref{L:PAQHR} gives
	\beqa
	\pa_q\Big(\widehat{h_R} F\Big)=\widehat{h_R} \pa_qF+\pa_q\big(\widehat{h_R}\big)F=\widehat{h_R} \pa_qF+Q_1\pa_L\big(\widehat{h_R}\big) F+Q_2\cos(2\pi m\theta)\pa_L\big(w\big)F .\nonumber
	\eeqa
	
	We next wish to integrate-by-parts with respect to $L$. To do that we use Fubini's theorem; to change the domain of integration, we observe that the boundary segment $\{(q,L)=(\om_R(L),L)\}$ corresponds to the curve $E=\Psi_L(R)$, at which $w$ and $\sin(2\pi m\th(R,E,L)$ both vanish. Moreover, as $|\om'|\geq c_0>0$ and $\pa_L\om=O(\varepsilon)$, we find that $\om_R(L)$ is an invertible function of $L$  in the regime $0<\varepsilon\ll1$, so that 
	\begin{align*}
		{}&\frac{1}{2\pi i m t}\int_{L_0}^{L_{\max}} \int_{\om_R(L)}^{\om_0(L)}e^{-2\pi mi q t}\Big(Q_1\pa_L\big(\widehat{h_R}\big)F +Q_2\cos(2\pi m\theta)\pa_L\big(w\big)F  \Big)\,\diff q\,\diff L\\
		&=\frac{1}{2\pi i m t}\int_{\om_*^{\min}}^{\om_0^{\max}}\int_{L_0}^{\om_R^{-1}(q)}e^{-2\pi mi q t}\Big(\widehat{h_R}\pa_L\Big(Q_1F \Big)+w\pa_L\Big(Q_2\cos(2\pi m\theta)F \Big) \Big)\,\diff L\,\diff q\\
		&\qquad\qquad-\frac{1}{2\pi i m t}\int_{\om_*^{\min}}^{\om_0^{\max}}\bigg[e^{-2\pi mi q t}\Big(Q_1\widehat{h_R}F +wQ_2\cos(2\pi m\theta)F  \Big)\bigg]\bigg|^{L=\om_R^{-1}(q)}_{L=L_0}\diff q\\
		&=-\frac{1}{2\pi i m t}\int_{L_0}^{L_{\max}} \int_{\om_R(L)}^{\om_0(L)}e^{-2\pi mi q t}\Big(\widehat{h_R}\pa_L\Big(Q_1F \Big)+w\pa_L\Big(Q_2\cos(2\pi m\theta)F \Big) \Big)\,\diff q\,\diff L,
	\end{align*}
	where we have used the above remark to eliminate the upper boundary term, and the lower boundary term vanishes because $F$ vanishes on $\partial J_{{\text{vac}}}$.
	We have again applied Fubini's theorem in the last line.
	To further simplify one of the resulting terms, we observe that for any sufficiently smooth function $j:J\to\mathbb C$ we have
	\begin{align}
		w\pa_L\Big(\cos(2\pi m\theta)j\Big)&=-2\pi m\sin(2\pi m\theta) w\big(\pa_E\th\frac{1}{\om'}+\pa_L\th\big)j+w\cos(2\pi m\theta)\pa_Lj\nonumber\\
		&=2\pi^2 m^2\widehat{h_R}\big(\frac{\pa_E r}{\om'}+\pa_Lr\big)\frac j T+w\cos(2\pi m\theta)\pa_L,
		\label{E:AUXLJ}
	\end{align}
	where we have again used~\eqref{E:PAER} and~\eqref{E:PALR}.
	
	Altogether, we have obtained
	\begin{align*}
		&\int_{L_0}^{L_{\max}} \int_{\om_R(L)}^{\om_0(L)}e^{-2\pi mi q t}\widehat{h_R} F\,\diff q\,\diff L\\
		&=-\frac{1}{2\pi i m t}\int_{L_0}^{L_{\max}} \int_{\om_R(L)}^{\om_0(L)} e^{-2\pi mi q t} \bigg[\widehat{h_R}\pa_L\Big(Q_1F \Big)+w\cos(2\pi m\theta)\pa_L\Big(Q_2F \Big)\\
		&\hspace{48mm} +2\pi^2 m^2\widehat{h_R}\big(\frac{\pa_E r}{\om'}+\pa_Lr\big)\frac{Q_2}{T}F - \widehat{h_R}\pa_q F\bigg]\diff q\,\diff L\\
		& = -\frac{1}{2\pi i m t}\int_{L_0}^{L_{\max}} \int_{\om_R(L)}^{\om_0(L)} e^{-2\pi mi q t}
		\Bigg[\widehat{h_R}\Big[\pa_L\Big(Q_1 F\Big) +2\pi^2 m^2\big(\frac{\pa_E r}{\om'}+\pa_Lr\big)\frac{Q_2}{T}F - \pa_q F\Big] \\
		&\hspace{48mm} + w\cos(2\pi m\theta)\pa_L\Big(Q_2 F\Big)\,  \Bigg]\diff q\,\diff L ,
	\end{align*}
	where we have also used~\eqref{E:AUXLJ} with $j = Q_2F $. This is precisely~\eqref{E:INTID1}.

	We note that, by~\eqref{eq:dLwidentity},
	$$\pa_q \big(w(R,E(q,L),L)\big)=-\frac{1}{\om'}\big(\frac{\pa_L\om}{\om'}+\frac{1}{2R^2}\big)^{-1}\pa_L\big(w(R,E(q,L),L)\big),$$
	so that
	\begin{align}
		{}\pa_q\big(w\cos(2\pi& m\th(R,E(q,L),L))\big) \notag \\
		&=-\frac{1}{\om'}\big(\frac{\pa_L\om}{\om'}+\frac{1}{2R^2}\big)^{-1}\cos(2\pi m \th)\pa_L(w)+2\pi m\sin(2\pi m\th)\frac{\pa_E r}{T\om'} \notag\\
		& =-\frac{1}{\om'}\big(\frac{\pa_L\om}{\om'}+\frac{1}{2R^2}\big)^{-1}\cos(2\pi m \th)\pa_L(w)+2\pi^2 m^2\widehat{h_R}\frac{\pa_E r}{T\om'},\notag
	\end{align}
	where we have used the chain rule and re-written $\pa_q \th=\pa_E\th \frac{1}{\om'}$ and $\pa_E\th=-\frac{\pa_E r}{Tw}$. 
	Therefore, with the test function $F$ vanishing at $\partial J_{{\mathrm{vac}}}$
	and~$w$ vanishing at $q=\omega_R(L)$ 
	we have
	\begin{align}
		&\int_{L_0}^{L_{\max}} \int_{\om_R(L)}^{\om_0(L)} e^{-2\pi mi q t} w\cos(2\pi m\th) F \diff q\diff L \notag\\
		& =  \frac1{2\pi i mt} \int_{L_0}^{L_{\max}} \int_{\om_R(L)}^{\om_0(L)}e^{-2\pi mi q t} \Big[\pa_q \Big(w\cos(2\pi m\th)\Big) F + w\cos(2\pi m\th) \pa_q F\Big]\diff q \diff L  \notag\\
		& =   \frac1{2\pi i mt} \int_{L_0}^{L_{\max}} \int_{\om_R(L)}^{\om_0(L)} e^{-2\pi mi q t}\bigg[-
		\frac{1}{\om'}\big(\frac{\pa_L\om}{\om'}+\frac{1}{2R^2}\big)^{-1}\cos(2\pi m \th) F \pa_L(w)  \notag\\
		&\hspace{48mm} + \widehat{h_R} 2\pi^2 m^2\frac{\pa_E r}{T\om'} F + w\cos(2\pi m\th) \pa_q F\bigg]\dif q\,\dif L. \label{E:FAUX1}
	\end{align}
	
	We now apply Fubini and integrate by parts with respect to $L$ in the first term above to get
	\begin{align}
		&  \frac1{2\pi i mt} \int_{L_0}^{L_{\max}} \int_{\om_R(L)}^{\om_0(L)} e^{-2\pi mi q t}
		\frac{1}{\om'}\big(\frac{\pa_L\om}{\om'}+\frac{1}{2R^2}\big)^{-1}\cos(2\pi m \th) F \pa_L(w) \diff q \diff L \notag \\
		& = -\frac1{2\pi i mt} \int_{\om_*^{\min}}^{\om_0^{\max}}\int_{L_0}^{\om_R^{-1}(q)} e^{-2\pi mi q t}
		\frac{1}{\om'}\big(\frac{\pa_L\om}{\om'}+\frac{1}{2R^2}\big)^{-1}\cos(2\pi m \th) F \pa_L(w) \diff L \diff q \notag\\
		& = - \frac1{2\pi i mt} \int_{L_0}^{L_{\max}} \int_{\om_R(L)}^{\om_0(L)} e^{-2\pi mi q t}\Bigg[
		2\pi^2 m^2\widehat{h_R}\big(\frac{\pa_E r}{\om'}+\pa_Lr\big)\frac F{\om' \big(\frac{\pa_L\om}{\om'}+\frac{1}{2R^2}\big) T} \notag\\
		&\hspace{48mm} +w\cos(2\pi m\theta)\pa_L\left(\frac{F}{\om'}\big(\frac{\pa_L\om}{\om'}+\frac{1}{2R^2}\big)^{-1}\right)\Bigg] \diff q \diff L, \label{E:FAUX2}
	\end{align}
	where no boundary terms are present due to $F=0$ on $\partial J_{{\mathrm{vac}}}$, and the fact that both $\widehat{h_R}$ and $w$ vanish at $E=\Psi_L(R)$. In the last line above we used~\eqref{E:AUXLJ} with $j =\frac{F}{\om'}\big(\frac{\pa_L}{\om'}+\frac{1}{2R^2}\big)^{-1}$. We now combine~\eqref{E:FAUX1}--\eqref{E:FAUX2} to obtain~\eqref{E:INTID2}.
\end{proof}

%%%%%%%%%%%%%%%%%%%%%%%
%%%%%%%%%%%%%%%%%%%%%%%

We now iterate the previous lemma to obtain a suitable high-order version of it, which will be used to obtain the time-decay of the gravitational force field. 

%%%%%%%%%%%%%%%%%%%%%%%
%%%%%%%%%%%%%%%%%%%%%%%

\begin{lemma}\label{L:PQ}
	Let $p\in\mathbb N$, $m\in\Z_\ast$,  and $R\in[\Rmin,\Rmax]$ be given. Assume that for any $\theta\in\mathbb S^1$, the map $(E,L)\mapsto r(\theta,E,L)\in C^{p+1}(\bar J)$, $T\in C^{p+1}(\bar J)$. Moreover,  let $F\in C^p(J)$ and  assume that for all $j\in\{0,\dots,p-1\}$,
	$\pa^{(j)}F$ vanishes on~$\partial J_{\mathrm{vac}}$ (see~\eqref{E:JVACDEF}).
	Then the following formula holds
	\begin{align}
		&\int_{L_0}^{L_{\max}} \int_{\Psi_L(R)}^{E_0} \widehat{h_R}(m,E,L) \fhat_0 T \diff E\diff L \notag\\
		& = \frac{1}{(-2\pi i m t)^p} \int_{L_0}^{L_{\max}} \int_{\om_R(L)}^{\om_0(L)} e^{-2\pi mi q t}\Bigg[ \widehat{h_R}  \P_p(\frac{\fhat_0 T}{\om'})
		+ w\cos(2\pi m \theta) \Q_p(\frac{\fhat_0 T}{\om'})\Bigg] \diff q \diff L,\nonumber
	\end{align}
	as long as the integrand of the right-hand side belongs to $L^1(J)$. Here
	the operators $\P_p, \Q_p$ are given through the following formula
	\begin{align}\label{DEF:PJQJ}
		\begin{pmatrix}
			\P_{p}\\ \Q_{p}
		\end{pmatrix}
		=
		\begin{pmatrix}
			\A_1 & \B_1 \\ \A_2 & \B_2
		\end{pmatrix}^{p-1} 
		\begin{pmatrix}
			\A_1\\ \A_2
		\end{pmatrix}.
	\end{align}
\end{lemma}

%%%%%%%%%%%%%%%%%%%%%%%
\begin{proof}
	Changing variables to $(q,L)$ we obtain
	\begin{align}
		\int_{L_0}^{L_{\max}} \int_{\Psi_L(R)}^{E_0} e^{-2\pi mi \om(E,L) t} \widehat{h_R} \fhat_0 T\diff E\,\diff L 
		=\int_{L_0}^{L_{\max}} \int_{\om_R(L)}^{\om_0(L)}e^{-2\pi mi q t}\widehat{h_R} \frac{\fhat_0 T}{\om'}\,\diff q\,\diff L.\nonumber
	\end{align}
	We now apply Lemma~\ref{L:RECURSION1} successively $p$-times and the above formula follows. The vanishing at vacuum assumptions are used to eliminate the boundary integrals arising at $E=E_0$ and $L=L_0$. At turning points of the flow, which represent another boundary contribution from the integration, the radial velocity $w$ vanishes and we also use the identity~\eqref{E:TPVANISH}, which implies that $\widehat{h_R}$ vanishes at those points.
\end{proof}

%%%%%%%%%%%%%%%%%%%%%%%
%%%%%%%%%%%%%%%%%%%%%%%

\begin{lemma}
	Let $k>1$ and $\ell>0$ be the given polytropic indices of the steady state~\eqref{E:POLY2}. 
	\begin{enumerate}
		\item[(a)]
		Given $f_0=|\varphi'(E,L)|g_0$ for some $g_0\in W^{N,\infty}(\Om)$, $N\in\N$, 
		there exists a constant $C>0$ such that for any $m\in\Z_\ast$ we have, for $(E,L)\in J$, $j\leq N$,
		\begin{multline}
			\Big| \pa^{(j)}\fhat_0(m,E,L)\Big|\\\le C\sum_{n=0}^j\frac{1}{|m|^{N-n}}(E-\EminL)^{\frac12-n}\Big(\sum_{s=0}^{j-n}(E_0-E)^{k-1-s}(L-L_0)^{\ell-(j-n-s)}\Big).\label{E:FZEROCONTROL}
		\end{multline} 
		\item[(b)]
		Under the assumptions of part (a) with $g_0\in W^{2,\infty}(\Om)$, there exists a constant $C>0$ such that 
		for any $m\in\Z_\ast$ we have, for $(E,L)\in J$, 
		\begin{align}
			&\Big|\P_1\big(\frac{\fhat_0(m,\cdot,\cdot)T}{\om'}\big)\Big|  \leq C|\varphi'|\bigg(\sum_{n=0}^1\frac{1}{|m|^{2-n}}(E-\EminL)^{\frac12-n}\Big(\sum_{s=0}^{1-n}(E_0-E)^{-s}(L-L_0)^{-(1-n-s)}\Big)+1\bigg)\label{E:P1FBOUNDS}\\
			&\Big|\Q_1\big(\frac{\fhat_0(m,\cdot,\cdot)T}{\om'}\big)\Big| \leq C|\varphi'|\bigg(\sum_{n=0}^1\frac{1}{|m|^{2-n}}(E-\EminL)^{\frac12-n}\Big(\sum_{s=0}^{1-n}(E_0-E)^{-s}(L-L_0)^{-(1-n-s)}\Big)\label{E:Q1FBOUNDS}
		\end{align}
		and
		\begin{align}
			&\Big|\P_2\big(\frac{\fhat_0(m,\cdot,\cdot)T}{\om'}\big)\Big| \leq C|\varphi'|\bigg(\sum_{n=0}^2\frac{1}{|m|^{2-n}}(E-\EminL)^{\frac12-n}\Big(\sum_{s=0}^{2-n}(E_0-E)^{-s}(L-L_0)^{-(2-n-s)}\Big)\notag\\
			&\ \ \ +|m|(E-\EminL)^{-\frac12}+(E_0-E)^{-1}+(L-L_0)^{-1}+m^2(E-\EminL)^{-\frac12}\Big(\frac{|m^2\widehat{g_0}(m,E,L)|}{(E-\EminL)^{\frac12}}\Big)\bigg)\label{E:P2FBOUNDS}\\
			&\Big|\Q_2\big(\frac{\fhat_0(m,\cdot,\cdot)T}{\om'}\big)\Big| \leq C|\varphi'|\bigg(\sum_{n=0}^2\frac{1}{|m|^{2-n}}(E-\EminL)^{\frac12-n}\Big(\sum_{s=0}^{2-n}(E_0-E)^{-s}(L-L_0)^{-(2-n-s)}\Big)\notag\\
			&\ \ \ +\Big|m\pa^{(1)}\Big(\frac{m\widehat{g_0}(m,E,L)}{(E-\EminL)^{\frac12}}\Big)\Big|+(E_0-E)^{-1}+(L-L_0)^{-1}\bigg).\label{E:Q2FBOUNDS}
		\end{align}
	\end{enumerate}
\end{lemma}

%%%%%%%%%%%%%%%%%%%%%%%

\begin{proof}
	Claim~\eqref{E:FZEROCONTROL} is a trivial consequence of  the explicit structures of~$f_0$ and~$\varphi$ and Lemma~\ref{lemma:Eminvanishing} which with trivial modifications applies to functions on $\mathbb S^1\times J$ and, by Lemma~\ref{L:RELCONTROL}, also to their derivatives with respect to~$E$ and~$L$.
	
	The estimate~\eqref{E:P1FBOUNDS} follows now directly from~\eqref{E:FZEROCONTROL} and the formulae~\eqref{DEF:AJOPS},~\eqref{DEF:BJOPS}, and~\eqref{DEF:PJQJ}. More precisely, we observe that
	\begin{multline*}
		\P_1(\frac{\fhat_0(m,\cdot,\cdot)T}{\om'}) =\A_1(\frac{\fhat_0(m,\cdot,\cdot)T}{\om'}) \\
		= \pa_L\Big(Q_1 \frac{\fhat_0(m,\cdot,\cdot)T}{\om'}\Big)  +2\pi^2 m^2\big(\frac{\pa_E r}{\om'}+\pa_Lr\big)\frac{Q_2\fhat_0(m,\cdot,\cdot)}{\om'} -\pa_q(\frac{\fhat_0(m,\cdot,\cdot)T}{\om'}).
	\end{multline*}
	We note that $|\pa_Er|,|\pa_Lr|\lesssim (E-\EminL)^{-\frac12}$ and since $g_0\in W^{2,\infty}(\Om)$ we also have
	$|\ghat_0(m,\cdot,\cdot)|\le   \frac{C\sqrt{E-\EminL}}{m^2}$ by Lemma~\ref{lemma:Eminvanishing}. Since $Q_2$ is bounded on $\bar J$ by Lemma~\ref{L:PAQHR}, the second summand above is bounded by  a constant multiple of $|\varphi'(E,L)|$. 
	Since $|\pa^{(1)}\ghat_0(m,E,L)| \lesssim \frac1{|m|}|\widehat{\pa^{(1)}\pa_{\theta}g_0}|\lesssim \frac1{|m|}(E-\EminL)^{-\frac12}$, the third term on the right-hand side above is bounded by 
	\begin{equation*}
		\big| \pa_q(\frac{\fhat_0(m,\cdot,\cdot)T}{\om'}) \big|
		\lesssim \frac1{|m|}(E-\EminL)^{-\frac12}|\varphi'(E,L)| + \frac1{|m|^2}(E_0-E)^{k-2}(L-L_0)^\ell (E-\EminL)^{\frac12},
	\end{equation*}
	and similarly,
	\begin{equation*}
		\big|\pa_L\Big(Q_1 \frac{\fhat_0(m,\cdot,\cdot)T}{\om'}\Big) \big|
		\lesssim \frac1{|m|}(E-\EminL)^{-\frac12}|\varphi'(E,L)| + \frac1{|m|^2}(E_0-E)^{k-1}(L-L_0)^{\ell-1} (E-\EminL)^{\frac12}.
	\end{equation*}
	The estimate~\eqref{E:Q1FBOUNDS} also follows by a bound similar to the latter one. 
	
	To verify the bound~\eqref{E:P2FBOUNDS}, we further distribute derivatives, as 
	\begin{align}
		\P_2(\frac{\fhat_0(m,\cdot,\cdot)T}{\om'})& =\A_1\A_1(\frac{\fhat_0(m,\cdot,\cdot)T}{\om'}) +\B_1\A_2(\frac{\fhat_0(m,\cdot,\cdot)T}{\om'}).\notag
	\end{align}
	We demonstrate the estimate for the first summand on the right, as the second follows from similar considerations. Expanding from~\eqref{DEF:AJOPS}, we observe
	\beqa\label{E:A1A1}
	{}&\A_1^2\Big(\frac{\fhat_0(m,\cdot,\cdot)T}{\om'}\Big)\\
	&=\pa_L\Big(Q_1\pa_L\big(Q_1\frac{\fhat_0 T}{\om'}\big)\Big)-\pa_q\pa_L\big(Q_1\frac{\fhat_0 T}{\om'}\big)-\pa_L\Big(Q_1\pa_q\big(\frac{\fhat_0 T}{\om'}\big)\Big)-\pa^2_q\big(\frac{\fhat_0 T}{\om'}\big)\\
	&\ \ \ +2\pi^2m^2\bigg(\pa_L\Big(Q_1\big(\frac{\pa_Er}{\om'}+\pa_Lr\big)Q_2\frac{\fhat_0 T}{\om'}\Big)-\pa_q\Big(\big(\frac{\pa_Er}{\om'}+\pa_Lr\big)Q_2\frac{\fhat_0 T}{\om'}\Big)\bigg)\\
	&\ \ \  +2\pi^2m^2\big(\frac{\pa_Er}{\om'}+\pa_Lr\big)\,\frac{Q_2}T\,\big(\partial_L(Q_1F)-\partial_qF \big)+4\pi^4m^4\big(\frac{\pa_Er}{\om'}+\pa_Lr\big)^2Q_2^2\frac{\fhat_0 T}{\om'}.
	\eeqa
	The first line is easily bounded using~\eqref{E:FZEROCONTROL} to control the derivatives of $\fhat_0$, using also that $T\in C^3$, $|\pa^{(j)}Q_i|\leq C(E-\EminL)^{\frac12-j}$ for $i,j=1,2$.  
	We therefore focus on the remaining terms. Starting with the middle line, we first apply~\eqref{ineq:paELrtildebounds} and \eqref{ineq:dQj} to estimate
	\begin{equation*}
		\Big|\pa^{(1)}\Big(\big(\frac{\widetilde{\pa_E r}}{\om'}+\widetilde{\pa_Lr}\big)Q_1Q_2\Big)\Big|\leq C(E-\EminL)^{-\frac12},
	\end{equation*}
	so that
	\beqa
	{}&\bigg|\pa^{(1)}\Big(Q_1 2\pi^2 m^2\big(\frac{\pa_E r}{\om'}+\pa_Lr\big)\frac{Q_2\fhat_0(m,\cdot,\cdot) T}{\om'}\Big)\bigg|\\
	{}&\leq\bigg|\pa^{(1)}\Big(Q_1 2\pi^2 m^2\big(\frac{(\pa_E r)_0}{\om'}+(\pa_Lr)_0\big)\frac{Q_2\fhat_0(m,\cdot,\cdot) T}{\om'}\Big)\bigg|\\
	{}&\hspace{38mm}+\bigg|\pa^{(1)}\Big(Q_1 2\pi^2 m^2\big(\frac{\widetilde{\pa_E r}}{\om'}+\widetilde{\pa_Lr}\big)\frac{Q_2\fhat_0(m,\cdot,\cdot) T}{\om'}\Big)\bigg|\\
	&\leq C\bigg(|m|(E-\EminL)^{-\frac12}(E_0-E)^{k-1}(L-L_0)^\ell+(E_0-E)^{k-2}(L-L_0)^\ell+(E_0-E)^{k-1}(L-L_0)^{\ell-1}\bigg),\nonumber
	\eeqa
	where we notice that $$\Big|\Big(\frac{m^2\widehat{g_0}(m,E,L)}{(E-\EminL)^{\frac12}}\Big)\Big|+(E-\EminL)^{\frac12}\Big|\pa^{(1)}\Big(\frac{m\widehat{g_0}(m,E,L)}{(E-\EminL)^{\frac12}}\Big)\Big|\leq C\|g_0\|_{W^{2,\infty}},$$
	by arguments similar to Lemma~\ref{L:IMPROVED}.
	The remaining $m^2$-terms can be estimated in the same way.
	
	Finally, we control the terms containing a power $m^4$ by simply estimating 
	\beqa
	\bigg|4\pi^4 m^4\big(\frac{\pa_E r}{\om'}+\pa_Lr\big)^2\frac{Q_2^2\fhat_0(m,\cdot,\cdot)}{(\om')^2}\bigg|\leq C\frac{m^2}{(E-\EminL)^{\frac12}}\frac{|m^2\widehat{f_0}(m,E,L)|}{(E-\EminL)^{\frac12}}.\nonumber
	\eeqa
	This concludes the proof of~\eqref{E:P2FBOUNDS}. To prove~\eqref{E:Q2FBOUNDS} we 
	note that by Lemma~\ref{L:PQ} we have $\Q_2 = \A_2\A_1+\B_2\A_2$, where we recall~\eqref{DEF:AJOPS} and~\eqref{DEF:BJOPS}. Focusing first on $\A_2\A_1$ we have
	\begin{align}
		& \A_2\A_1\Big(Q_1 \frac{\fhat_0(m,\cdot,\cdot)T}{\om'}\Big) =
		\pa_L\Big(Q_2 \pa_L\Big(Q_1 \frac{\fhat_0(m,\cdot,\cdot)T}{\om'}\Big) \Big) -\pa_L\Big(Q_2\pa_q(\frac{\fhat_0(m,\cdot,\cdot)T}{\om'})\Big)\notag\\
		& \qquad +2\pi^2 m^2\pa_L\Big(Q_2\big(\frac{\pa_E r}{\om'}+\pa_Lr\big)\frac{Q_2\fhat_0(m,\cdot,\cdot)}{\om'}\Big). \label{E:A2A1}
	\end{align}
	The first line is again easily bounded using~\eqref{E:FZEROCONTROL} to control the derivatives of $\fhat_0$, and is bounded by the first line of the right-hand side of~\eqref{E:Q2FBOUNDS}. To control the second line of~\eqref{E:A2A1}, we notice that it retains the same singular structure as the second line on the right-hand side of~\eqref{E:A1A1}. Similar argument applies to the term $\B_2\A_2$. Using the same ideas as above, we obtain~\eqref{E:Q2FBOUNDS}.
\end{proof}

\section{Decay in the 1+1-dimensional case}\label{S:DECAY}

%%%%%%%%%%%%%%%%%%%%%%%%%%%%%%%
%%%%%%%%%%%%%%%%%%%%%%%%%%%%%%%
%%%%%%%%%%%%%%%%%%%%%%%%%%%%%%%

The purpose of this section is to prove Theorem~\ref{T:MAIN} concerning
decay for the 1+1-dimensional case.

Our first goal is to obtain a quantitative decay estimate for the time-derivative of the spatial density $\rho_f(t,r)$ defined in~\eqref{E:DENSITYDEF}. 
To that end, we shall need to be careful in the vicinity of the elliptic point $r=r_\ast$ (see Figure~\ref{F:PW}), so we first introduce 
a pair of useful integral estimates.

\begin{lemma}
Let $\al,\beta<0$, $A>0$, $\de\in(0,\frac{1}{4}A)$. Then, for any $a\in[4\de,A]$, there exists a constant $C>0$, depending on $A,\al,\beta$, but independent of $\de,a$, such that 
\begin{align}
\int_\de^{a-\de}(a-x)^\al x^\beta\,\dif x\leq C\de^{\al+\beta+1}, \ \ \ \text{ if } \al+\beta<-1,\label{ineq:integral1}\\
\int_\de^{2\de}(a-x)^\al x^\beta\,\dif x\leq C\de^{\al+\beta+1}, \ \ \ \text{ if } \al+\beta\in(-1,0).\label{ineq:integral2}
\end{align}
\end{lemma}
The proof follows from direct evaluation of the integrals and the standard properties of hypergeometric functions.

%%%%%%%%%%%%%%%%%%%%%%%
%%%%%%%%%%%%%%%%%%%%%%%

\begin{prop}[Decay of macroscopic density]\label{P:RHOTDECAY}
Let $f$ be the solution to the pure transport equation~\eqref{E:PT0} with initial data  $f_0=|\varphi'(E)|g_0$ for some $g_0\in W^{1,\infty}(\Om)$ and suppose $k>\frac12$. 
Then there exists a constant $C=C\|g_0\|_{W^{1,\infty}(\overline\Om)}>0$, depending on $k$, such that for all $t\in\R$,
\begin{align}\label{E:DECAYRHOROOT}
\| \pa_t\rho_f(t,\cdot) \|_{L^\infty_r} \le C\begin{cases}(1+|t|)^{-\frac12}, & k\geq 1,\\
(1+|t|)^{-(k-\frac12)}, & k\in(\frac12,1).
\end{cases}
\end{align}
Moreover, if $f_0$ additionally satisfies the orthogonality condition $f_0\in \text{{\em ker}}\,\D^\perp$ (recall~\eqref{E:ORTH}),
there exists a constant $C>0$ such that 
\begin{align}\label{E:DECAYJUSTRHO}
\| \rho_f(t,\cdot) \|_{L^\infty_r} \le C\begin{cases}(1+|t|)^{-\frac12}, & k\geq 1,\\
(1+|t|)^{-(k-\frac12)}, & k\in(\frac12,1).
\end{cases}
\end{align}
\end{prop}

\begin{proof}
For the purpose of proving decay of $\pa_t\rho_f$
it is equivalent to work with the modified density
\begin{equation*}
\tilde\rho(t,r) \coloneqq \int f \diff w = \frac{r^2}\pi \rho_f(t,r).
\end{equation*}
Then
\begin{align}
	\pa_t\tilde\rho & = - \int \omega(E) \pa_\theta f \diff w \nonumber\\
	&= -  \int_{\Psi(r)}^{E_0} \omega(E)\,\frac{\pa_\theta f(t,\theta(r,E),E)+\pa_\theta f(t,1-\theta(r,E),E)}{\sqrt{2(E-\Psi(r))}} \diff E\eqqcolon\tilde\rho_{t1}+\tilde\rho_{t2},\label{E:TILDERHOPARTIALT}
\end{align} 
where we have used the change of variables $w\mapsto E$ for a given fixed $r$.
 An immediate consequence of~\eqref{E:TILDERHOPARTIALT} is that
\begin{equation}\label{E:BOUNDEDNESSRHO}
	\| \pa_t\rho_f(t,\cdot) \|_{L^\infty_r} \le C,\qquad t\in\R.
\end{equation}
As for the decay estimate~\eqref{E:DECAYRHOROOT}, 
we focus on the first term~$\tilde\rho_{t1}$ on the right-hand side of~\eqref{E:TILDERHOPARTIALT}.
The estimates for the second term~$\tilde\rho_{t2}$ follow in the same way by replacing the initial value~$f_0$ with $(\theta,E)\mapsto f_0(1-\theta,E)$.
We first recall the solution formula~\eqref{E:EXPLICIT}--\eqref{E:SEMIGROUP} to infer
\begin{align}
\pa_E\left(f_0\circ T_t\right) = (\theta_E - \omega'(E)t)\,\pa_\theta f_0\circ T_t + \pa_Ef_0\circ T_t,\label{E:PEFCIRCT}
\end{align}
which yields the identity
\begin{align}\label{E:RLRHO}
\pa_\theta f_0\circ T_t = \frac{-1}{\omega'(E)t} \left(\pa_E\left(f_0\circ T_t\right) - \pa_\theta f_0\circ T_t \,\theta_E- \pa_E f_0 \circ T_t \right).
\end{align}
Note that the left-hand side of~\eqref{E:RLRHO} is precisely $\pa_\theta f$ by~\eqref{E:EXPLICIT}
and therefore from~\eqref{E:TILDERHOPARTIALT} and~\eqref{E:RLRHO} we obtain
\begin{align}
\tilde\rho_{t1}(t,r) & = -  \int_{\Psi(r)}^{E_0} \omega(E) \frac{\pa_\theta f_0\circ T_t}{\sqrt{2(E-\Psi(r))}} \diff E \notag\\
& = \frac1t  \int_{\Psi(r)}^{E_0} \frac{\omega(E)}{\omega'(E)} \frac{\pa_E\left(f_0\circ T_t\right) - \pa_\theta f_0\circ T_t \theta_E- \pa_Ef_0 \circ T_t }{\sqrt{2(E-\Psi(r))}} \diff E .\label{E:TILDERHOPARTIALT2}
\end{align}
We now recall that $f_0=|\varphi'(E)|g_0$, where $g_0\in W^{1,\infty}(\overline{\Om})$, and so
$$|\pa_\th f_0|\leq C\sqrt{E-E_{\min}}(E_0-E)^{k-1}.$$
Let $\de=\frac1t>0$. Now, in the case that $E_0-\Psi(r)\leq 4\de$,  we easily obtain
\beqa\label{ineq:rhotbdlarger}
\big|\tilde\rho_{t1}(t,r)\big|  \leq&\, C  \int_{\Psi(r)}^{E_0}  \frac{(E_0-E)^{k-1}}{\sqrt{E-\Psi(r)}} \diff E\leq \begin{cases}
C\de^{k-\frac12}, & k\in(\frac12,1],\\
C\de^{\frac12}, & k>1,
\end{cases}
\eeqa
where we have applied~\eqref{ineq:integral2} in the first case.
 
We now suppose that $E_0-\Psi(r)\geq 4\de$. Let $\phi\in C_c^\infty(\Psi(r)+\frac{\de}{2},E_0-\frac{\de}{2})$ be such that $0\leq \phi\leq 1$,  $\phi(E)\equiv 1$ for $E\in(\Psi(r)+\de,E_0-\de)$, and $|\phi'(E)|\leq \frac4\de$.
 For convenience of notation, we define
\begin{equation*}
I_\de=(\Psi(r)+\frac{\de}{4},E_0-\frac\de4).
\end{equation*}
First we re-write \eqref{E:TILDERHOPARTIALT2} as 
\begin{align*}
\tilde\rho_{t1}(t,r) & =\frac2t  \int_{\Psi(r)}^{E_0} \frac{\phi(E)\omega(E)}{\omega'(E)} \frac{\pa_E\big(f_0\circ T_t\big) }{\sqrt{2(E-\Psi(r))}} \diff E \\
& -\frac2t  \int_{\Psi(r)}^{E_0} \frac{\phi(E)\omega(E)}{\omega'(E)} \frac{(\pa_\theta f_0)\circ T_t \theta_E+ (\pa_Ef_0 )\circ T_t }{\sqrt{2(E-\Psi(r))}} \diff E \\
& - \int_{\Psi(r)}^{E_0} (1-\phi(E))\omega(E) \frac{\pa_\theta f_0\circ T_t}{\sqrt{2(E-\Psi(r))}} \diff E \notag\\
& \eqqcolon I(t,r) + II(t,r) + III(t,r).
\end{align*}
The main term is $I(t,r)$, and so we first treat this. Observing that $\phi$ vanishes at both end-points of the interval, we integrate by parts with respect to $E$ and obtain
\beqa
I(t,r)=&\,-\frac2t  \int_{\Psi(r)}^{E_0}\bigg( \pa_E\Big(\frac{\phi(E)\omega(E)}{\omega'(E)}\Big) \frac{f_0\circ T_t }{\sqrt{2(E-\Psi(r))}} -\frac{\phi(E)\omega(E)}{\omega'(E)}\frac{f_0\circ T_t }{(2(E-\Psi(r)))^{\frac32}} \bigg)\diff E \nonumber
\eeqa
Recall that $T\in C^2(I)$ and hence $\om\in C^2(I)$ also by~\cite[Lemma 3.4]{HRSS2023}. Thus, we obtain the simple estimate 
$$\Big|\pa_E\Big(\frac{\phi(E)\omega(E)}{\omega'(E)}\Big)\Big|\leq \begin{cases}
C\de^{-1}, & E\in (\Psi(r)+\frac\de2,\Psi(r)+\de)\cup(E_0-\de,E_0-\frac\de2),\\
C, & E\in[\Psi(r)+\de,E_0-\de].
\end{cases}$$
Therefore, as $|f_0|\leq C(E_0-E)^{k-1}$, we have the estimate
\beqa\label{ineq:IBD}
\big|I(t,r)\big|\leq&\,\frac{C}{\de t}\bigg(\int_{\Psi(r)+\frac\de2}^{\Psi(r)+\de}(E_0-E)^{k-1}(E-\Psi(r))^{-\frac12}\,\diff E+\int_{E_0-\de}^{E_0-\frac\de2}(E_0-E)^{k-1}(E-\Psi(r))^{-\frac12}\,\diff E\bigg)\\
&+\frac{C}{t}\int_{\Psi(r)+\de}^{E_0-\de}(E_0-E)^{k-1}(E-\Psi(r))^{-\frac32}\,\diff E\\
\leq&\, C\begin{cases}
t^{-1}\de^{k-\frac32}=Ct^{-(k-\frac12)}, & k\in(\frac12,1),\\
t^{-1}\de^{-\frac12}=C t^{-\frac12}, & k\geq 1,
\end{cases}
\eeqa
where we have applied~\eqref{ineq:integral1}--\eqref{ineq:integral2} to estimate the integrals in the first case.

For term $II(t,r)$, we first obtain the estimate
$$\big|(\pa_\theta f_0)\circ T_t \theta_E+ (\pa_E f_0) \circ T_t \big|\leq\frac {C(E_0-E)^{k-1}}{\sqrt{E-\Psi(r)}} + C(E_0-E)^{k-2},$$
where we have used that $|\th_E|\leq\frac{C}{\sqrt{E-\Emin}\sqrt{E-\Psi(r)}}$  by Lemma~\ref{L:THETAREG}, $|\pa_\th f_0|\leq C(E_0-E)^{k-1}\sqrt{E-E_{\min}}$, and 
$$|\pa_Ef_0|\leq |\varphi'||\pa_Eg_0|+|\varphi''||g_0|\leq \frac{C}{\sqrt{E-E_{\min}}}+C(E_0-E)^{k-2},$$ along with $\Psi(r)\geq E_{\min}$.
Then
\beqa\label{ineq:IIBD}
|II(t,r)|\leq&\frac{C}{|t|}\int_{\Psi(r)+\frac{\de}{2}}^{E_0-\frac\de2}\bigg(\frac{(E_0-E)^{k-1}}{E-\Psi(r)}+(E_0-E)^{k-2}\bigg)\,\dif E\\
\leq&\,C\begin{cases} t^{-1}\de^{k-1}=t^{-k}, & k\in(\frac12,1)\\
t^{-1}|\log\de|=\frac{\log t}{t}, & k\geq 1,
\end{cases}
\eeqa
where we have again applied~\eqref{ineq:integral1}.

Finally, for $|III(t,r)|$, we again use $|\pa_\th f_0|\leq  C(E_0-E)^{k-1}$ to obtain the bound
\beqa\label{ineq:IIIBD}
|III(t,r)|\leq&\,C\int_{\Psi(r)}^{\Psi(r)+\de}\frac{(E_0-E)^{k-1}}{\sqrt{E-\Psi(r)}}\,\dif E+C\int_{E_0-\de}^{E_0-\frac\de2}\frac{(E_0-E)^{k-1}}{\sqrt{E-\Psi(r)}}\,\dif E\\
\leq&\, C\begin{cases}
\de^{k-\frac12}=t^{-(k-\frac12)}, & k\in(\frac12,1),\\
\de^{\frac12}=-t^{-\frac12}, & k\geq 1,
\end{cases}
\eeqa
where we have applied~\eqref{ineq:integral2}.

Thus, recalling \eqref{ineq:rhotbdlarger} and combining \eqref{ineq:IBD}--\eqref{ineq:IIIBD}, overall, we have obtained
\begin{equation*}
|\rho_t(t,r)|\leq C\begin{cases}
t^{-(k-\frac12)}, & k\in(\frac12,1),\\
t^{-\frac12}, & k\geq 1,
\end{cases}
\end{equation*}
which, combined with~\eqref{E:BOUNDEDNESSRHO}, yields~\eqref{E:DECAYRHOROOT}.

We next show~\eqref{E:DECAYJUSTRHO}. Since the assumption~\eqref{E:ORTH} is dynamically preserved by the solutions of~\eqref{E:PURETRANSPORT}, it follows that $f(t,\cdot,\cdot)\in (\text{ker}\D)^\perp$, and therefore there exists a unique $g(t,\cdot,\cdot)\in (\mathrm{ker}\,\D)^\perp$~\cite{HRS2022} such that 
\begin{equation*}
f = -\D g,
\end{equation*}
and moreover $g(t,\cdot,\cdot)$ is a solution of the transport equation~\eqref{E:PURETRANSPORT}. It follows that 
\begin{equation*}
\rho_f=- \int \omega(E) \pa_\theta g \diff w 
\end{equation*}
and we may now follow the proof of~\eqref{E:DECAYRHOROOT} verbatim.
\end{proof}

%%%%%%%%%%%%%%%%%%%%%%%%%%%%%%%
%%%%%%%%%%%%%%%%%%%%%%%%%%%%%%%

We next establish the decay for the gradient of the gravitational potential.
\begin{prop}\label{P:PARTIALRUFDECAY}
Let  $f_0=|\varphi'(E)|g_0$ where $g_0\in W^{2,\infty}(\Om)$ be such that $f_0\in \text{{\em ker}}\,\D^\perp$ and suppose $k>\frac12$. Let $f(t,\cdot,\cdot)$ be the solution to~\eqref{E:PT0} with initial data $f(0,\cdot,\cdot)=f_0(\cdot,\cdot)$. 
Then there
exists a constant $C=C\|g_0\|_{W^{2,\infty}}$ such that for all $t\in\R$,
\begin{equation}
\| \pa_rU_f(t,\cdot) \|_{L^\infty_r} \le \begin{cases} C(1+|t|)^{-\frac32}, & k\geq \frac32,\\
C(1+|t|)^{-k}, & k\in(\frac12,\frac32).
\end{cases} \label{E:PARTIALRUFDECAYEXTRA}
\end{equation}
\end{prop}

\begin{proof}
We prove the theorem in the case $k\geq 1$ as the remaining case, $k\in(\frac12,1)$ follows from similar arguments. Our starting point is the representation formula~\eqref{E:PARTIALRUFFORMULA}.
Switching to the action-angle variables and using the Parseval identity,  we arrive at 
\begin{align}
R^2 \pa_RU_{f}(t,R)  & = 4\pi^2 \int_I \int_{\mathbb S^1} h_R(r(\theta,E))\, f T(E) \diff\theta\diff E \notag\\
& = 4\pi^2 \sum_{m\in\Z^\ast} \int_I \fhat(t,m,E) \overline{\widehat{h_R}(m,E)} T(E)\,dE \notag \\
& = 4\pi \sum_{m\in\Z^\ast} \color{red} \color{black} \int_{\Psi(R)}^{E_0} e^{-2\pi mi \om(E) t} H_{f_0}(m,E) T(E)\diff E, \label{eq:paRU_finitial}
\end{align}
 where we recall the formula~\eqref{E:FHATREP}, use Lemma~\ref{L:GHFOURIER}, and recall the function $H_{f_0}(m,E)$ from~\eqref{E:HDEFNEW}. Using the nonstationary phase argument, we conclude that 
\begin{align}
&\int_{\Psi(R)}^{E_0} e^{-2\pi mi \om(E) t} H_{f_0}(m,E)  T(E)\diff E \notag\\
& = \frac1{2\pi i m t} \int_{\Psi(R)}^{E_0} e^{-2\pi mi \om(E) t} \frac{d}{dE}\left(H_{f_0}(m,E)\frac{T(E)}{\om'(E)}\right)\diff E \notag\\
& \ \ \ - \frac1{2\pi i m t} e^{-2\pi mi \om(E) t} \fhat_0(m,E)\sin(2\pi m \theta(R,E))  \frac{T(E)}{\om'(E)} \Big|^{E=E_0}_{E=\Psi(R)}.\label{eq:paRU_fphasem}
\end{align}
The boundary term at $E=\Psi(R)$ vanishes since there $\theta(R,E)=\theta(R,\Psi(R))=\frac12$ and therefore $\sin(2\pi m \theta(R,\Psi(R)))=0$. 

In the case $k=1$, the boundary term at $E_0$ is merely bounded by $C|m|^{-2}$ from~\eqref{E:ghatestimate} and the identity $f_0=\eps k(E_0-E)^{k-1}g_0$.  We then estimate the remaining integral term in this case by a further $C|m|^{-2}$ by noting that $\pa_EH_{f_0}=\eps \pa_E H_{g_0}$ and applying~\eqref{E:PARTIALEHNEWBOUND} to see that the integral converges.

When $k>1$, we notice that the boundary term at $E=E_0$ vanishes due to $f_0=\eps k(E_0-E)^{k-1}g_0$, and therefore only the integral remains in~\eqref{eq:paRU_fphasem}. 
To obtain~\eqref{E:PARTIALRUFDECAYEXTRA}, we first take $\de>0$ to be chosen later and separate the cases $E_0-\Psi(R)\leq 4\de$ and $E_0-\Psi(R)>4\de$.  In the former case, noting that $H_{f_0}=|\varphi'(E)|H_{g_0}$ and applying estimate~\eqref{E:PARTIALEHNEWBOUND}, we obtain
\beq\label{ineq:dEHf0}
\big|\frac{\dif}{\dif E} H_{f_0}(m,E)\big|\leq \frac{C}{|m|}(E-\Psi(R))^{-\frac12}+\frac{C}{|m|^2}(E_0-E)^{k-2}.
\eeq
Thus 
\beqa\label{ineq:dUfest1}
\bigg|&\int_{\Psi(R)}^{E_0} e^{-2\pi mi \om(E) t} H_{f_0}(m,E)  T(E)\diff E\bigg|\leq \frac{C}{ m^2 t} \int_{E_0-4\de}^{E_0}\Big((E-\Psi(R))^{-\frac12}+(E_0-E)^{k-2}\Big)\,\diff E\\
&\leq \frac{C}{m^2 t}\big(\de^{\frac12}+\de^{k-1}\big).
\eeqa
In the remaining case $E_0-\Psi(R)>4\de$, we introduce a cut-off function $\varphi_\de(E)\in C_c^\infty(\R)$ satisfying $0\leq\varphi_\de\leq 1$, $\varphi_\de(E)=0$ for $E\leq\Psi(R)+\de$ and $E\geq E_0-\de$, $\varphi_\de(E)=1$ for $\Psi(R)+2\de\leq E\leq E_0-2\de$, and $|\varphi_\de'(E)|\leq\frac{2}{\de}$. Then, arguing further from~\eqref{eq:paRU_fphasem},
\beqa
&\int_{\Psi(R)}^{E_0} e^{-2\pi mi \om(E) t} H_{f_0}(m,E)  T(E)\diff E \notag\\
&= \frac1{4\pi^2 m^2 t^2} \int_{\Psi(R)}^{E_0} e^{-2\pi mi \om(E) t}\pa_E\bigg(\varphi_\de(E) \frac{d}{dE}\Big(H_{f_0}(m,E)\frac{T(E)}{\om'(E)}\Big)\bigg)\diff E\\
&\ \ \ +\frac1{2\pi i m t} \int_{\Psi(R)}^{E_0} e^{-2\pi mi \om(E) t}\big(1-\varphi_\de(E)\big) \frac{d}{dE}\left(H_{f_0}(m,E)\frac{T(E)}{\om'(E)}\right)\diff E\\
& = : I_m+II_m.
\eeqa
Now, we observe that $H_{f_0}(m,E)=|\varphi'(E)|H_{g_0}(m,E)$, so that from Lemma~\ref{L:HIGHORDER}, specifically~\eqref{E:PARTIALEEH}, we obtain the bounds, for $\Psi(R)+\de\leq E\leq E_0-\de$,
\begin{align*}
\Big|&\frac{d}{dE}\Big( \frac{d}{dE}\Big(H_{f_0}(m,E)\frac{T(E)}{\om'(E)}\Big)\Big)\Big|\\
&\leq C(E-E_{\min})^{-\frac12}(E-\Psi(R))^{-1}+\frac{C}{|m|}(E_0-E)^{k-2}(E-\Psi(R))^{-\frac12}+\frac{C}{m^2}(E_0-E)^{k-3},\\
|&\varphi_\de'(E)|\Big|\frac{d}{dE}\Big(H_{f_0}(m,E)\frac{T(E)}{\om'(E)}\Big)\Big| \leq \frac{C}{\de|m|}(E-\Psi(R))^{-\frac12}+\frac{C}{m^2}(E_0-E)^{k-2},
\end{align*}
where we have used $|H_{g_0}(m,E)|\leq \frac{C}{m^2}$, $k>1$, and $T\in C^3(\bar I)$. The latter follows by arguments similar to~\cite[Lemma~3.4]{HRSS2023} and~\cite[Sc.~A.4]{St24},
where we use the boundary regularity $(\pa^4_r\Psi)(r)\sim (r-\Rmin)^{k-\frac32}$ and similarly when $r\sim\Rmax$.

Moreover,  we note the support of $\varphi_\de'(E)$ on $[\Psi(R),E_0]$ is contained in $[\Psi(R)+\de,\Psi(R)+2\de]\cup[E_0-2\de,E_0-\de]$, and therefore
\begin{align}
|I_m|\leq&\, \frac{C}{m^2 t^2}\bigg(\int_{\Psi(R)+\de}^{E_0-\de}\Big((E-\Psi(R))^{-\frac32}+(E_0-E)^{k-2}(E-\Psi(R))^{-\frac12}+(E_0-E)^{k-3}\Big)\,\diff E\notag\\
&\ \ \ \ \ +\frac{1}{\de}\int_{\Psi(R)+\de}^{\Psi(R)+2\de}\Big((E-\Psi(R))^{-\frac12}+(E_0-E)^{k-2}\Big)\,\diff E\notag\\
&\ \ \ \ \ +\frac{1}{\de}\int_{E_0-2\de}^{E_0-\de}\Big((E-\Psi(R))^{-\frac12}+(E_0-E)^{k-2}\Big)\,\diff E\bigg)\notag\\
\leq&\, \frac{C}{m^2t^2}\big(\de^{-\frac12}+\de^{k-2}\big)\label{ineq:dUfest2},
\end{align}
where, if $k\in(1,\frac32)$, we apply~\eqref{ineq:integral1} to estimate the middle term in the integrand on the first line, and if $k\geq \frac32$, we note this term gives rise to a uniformly bounded contribution.

On the other hand, we may apply~\eqref{ineq:dEHf0}
so that
\begin{align}
|II_m|\leq&\, \frac{C}{m^2t}\bigg(\int_{\Psi(R)}^{\Psi(R)+2\de}\Big((E-\Psi(R))^{-\frac12}+(E_0-E)^{k-2}\Big)\,\diff E\notag\\
&\ \ \ \ \ \ \ \ \ \ \ +\int_{E_0-2\de}^{E_0}\Big((E-\Psi(R))^{-\frac12}+(E_0-E)^{k-2}\Big)\,\diff E\bigg)\notag\\
\leq&\, \frac{C}{m^2t}\big(\de^{\frac12}+\de^{k-1}\big).\label{ineq:dUfest3}
\end{align}
Combining~\eqref{ineq:dUfest1}--\eqref{ineq:dUfest3} with~\eqref{eq:paRU_finitial} and choosing $\de=\frac1t$, we have obtained, for $k\geq \frac32$,
\begin{align}\label{E:PARTIALRUFDECAY}
\lv R^2 \pa_RU_{f}(t,R) \rv \le \frac {C}{t^{\frac32}} \sum_{m\in\Z_\ast}\frac1{m^2} \le \frac C{t^{\frac32}}, \ \ t>0,
\end{align}
while if $k\in(1,\frac32)$, we obtain instead
\begin{align}\label{E:PARTIALRUFDECAY2}
\lv R^2 \pa_RU_{f}(t,R) \rv \le \frac {C}{t^{k}} \sum_{m\in\Z_\ast}\frac1{m^2} \le \frac C{t^{k}}, \ \ t>0,
\end{align}
\end{proof}

%%%%%%%%%%%%%%%%%%%%%%%%
%%%%%%%%%%%%%%%%%%%%%%%%

\begin{remark}
It follows from the above proof that if $g_0\in W^{1,\infty}(\Om)$ then there exists a constant $C =C\|g_0\|_{W^{1,\infty}(\Om)}$ such that for any $t\in\mathbb R$
\begin{align}
\|\pa_r U_f(t,\cdot)\|_{L^\infty_r} \le \frac{C}{1+|t|}.
\end{align}
The boost in regularity of $g_0$ does give a faster decay by a factor of $|t|^{-\frac12}$ as shown in Proposition~\ref{P:PARTIALRUFDECAY}, but unlike the classical Landau damping~\cite{MoVi2011}, we do not expect
to obtain faster decay-in-time of the gravitational force $\pa_r U_f$ with more regular initial data $g_0$.
\end{remark}

%%%%%%%%%%%%%%%%%%%%%%%%
%%%%%%%%%%%%%%%%%%%%%%%%

\begin{corollary}\label{C:PARTIALRUFDECAY}
Under the assumptions of Theorem~\ref{T:MAIN}, there exists a constant $C=C\|g_0\|_{W^{2,\infty}}$ so that
\begin{align}
\|\pa_t\pa_r U_f(t,\cdot)\|_{L^2_r} & \le \frac C{1+t^{\frac32}}. \label{E:3OVER2}
\end{align}
\end{corollary}

%%%%%%%%%%%%%%%%%%%%%%%%

\begin{proof}
We note that due to~\eqref{E:PT0} and Remark~\ref{R:KERD} $\pa_tf\in\text{{\em ker}}\,(\D)^\perp$ and therefore the claim follows from Proposition~\ref{P:PARTIALRUFDECAY}.
\end{proof}

%%%%%%%%%%%%%%%%%%%%%%%%
%%%%%%%%%%%%%%%%%%%%%%%%

We next turn to the question of decay of the time-derivative of the gravitational potential $\pa_t U_f$.

%%%%%%%%%%%%%%%%%%%%%%%%%%%%%%%
%%%%%%%%%%%%%%%%%%%%%%%%%%%%%%%

\begin{prop}\label{P:DECAY3}
Let  $f_0=|\varphi'(E)|g_0$ for some $g_0\in W^{3,\infty}(\Om)$ and let $f(t,\cdot,\cdot)$ be the solution to~\eqref{E:PT0} with initial data $f(0,r,w)=f_0(r,w)$. Then, for $k\in(\frac12,2)$, there exists a constant 
$C=C\|g_0\|_{W^{3,\infty}}$ such that 
\begin{align}
\lv \pa_t U_f(t,r) \rv \le  C\,(1+t)^{-k}, \qquad t\in\R, \ r\neq r_\ast. \label{E:PARTIALTUFDECAYTRAPPEDk}
\end{align}
For $k\geq 2$, there
exists a constant $C=C(\|f_0\|_{W^{3,\infty}})$ such that 
\begin{align}
\lv \pa_t U_f(t,r) \rv \le  C\,\frac{1+|\log |r-r_\ast||}{1+t^2}, \qquad t\in\R, \ r\neq r_\ast. \label{E:PARTIALTUFDECAYTRAPPED}
\end{align}
Moreover, if $g_0\in W^{4,\infty}(\Om)$, $k\geq \frac52$, then in the outer region of the galaxy 
the pointwise decay rate are improved, i.e., there exists a constant $C>0$ such that
\begin{align}
\lv \pa_t U_f(t,r) \rv \le \frac{C  |r-r_\ast|^{-1}}{1+|t|^{\frac52}}, \qquad t\in\R, \ r>r_\ast. \label{E:PARTIALTUFDECAYEXTERIOR}
\end{align}
\end{prop}

%%%%%%%%%%%%%%%%%%%%%%%%%%%%%%%

\begin{proof}
Our starting point is the representation formula~\eqref{E:PARTIALTUFFORMULA}.
Switching to the action-angle variables and using the Parseval identity,  we arrive at 
\begin{align*}
\pa_tU_{f}(t,R)  & = 4\pi^2 \int_I \int_{\mathbb S^1} \frac{w}{r^2}\chi_{r(\theta,E)\ge R}\, f T(E) \diff\theta\diff E \notag\\
& = 4\pi^2 \sum_{m\in\Z^\ast} \int_I \fhat(t,m,E) \overline{\widehat{w g_R}(m,E)} T(E)\,dE \notag \\
& = 4\pi^2 \sum_{m\in\Z^\ast} \int_I e^{-2\pi mi \om(E) t} \fhat_0(m,E) \overline{\widehat{w g_R}(m,E)} T(E)\diff E \notag\\
& = 4\pi^2 \sum_{m\in\Z^\ast} \int_I e^{-2\pi mi \om(E) t}T(E) G_{f_0}(m,E)\diff E,
\end{align*}
where we recall the formula~\eqref{E:FHATREP}, and the definitions~\eqref{E:HDEF} and~\eqref{E:GRDEF} of~$G_{f_0}(m,E)$ and~$g_R$, respectively.
Let 
\begin{align}\label{E:LTMRDEF}
L(t,m;R)\coloneqq\int_I e^{-2\pi mi \om(E) t}  T(E) G_{f_0}(m,E)\diff E.
\end{align}
Our goal is to integrate-by-parts in~\eqref{E:LTMRDEF} to recover algebraic decay of $L(t,m;R)$. In the case that $k\in(\frac12,2)$, arguments entirely analogous to the proof of Proposition~\ref{P:PARTIALRUFDECAY}, but using now estimates for $G_{f_0}$ in the place of the corresponding estimates for $H_{f_0}$, yields~\eqref{E:PARTIALTUFDECAYTRAPPEDk}. We therefore focus on~\eqref{E:PARTIALTUFDECAYTRAPPED}--\eqref{E:PARTIALTUFDECAYEXTERIOR}.

Suppose that $k\geq 2$. If $R\in(r_\ast,\Rmax]$, we have, by~\eqref{E:WGREXPLICIT},
\begin{align}
L(t,m;R)  &= \frac{-1}{2\pi i m t}\int_{\Psi(R)}^{E_0} \frac{d}{dE}\left(e^{-2\pi m i\omega(E) t}\right) \frac{ T(E)}{\omega'(E)} G_{f_0}(m,E) \diff E \notag\\
& = \frac{-1}{2\pi i m t} \left[e^{-2\pi m i\omega(E) t} \frac{G_{f_0}(m,E)T(E)}{\omega'(E)}\right]^{E=E_0}_{E=\Psi(R)} \notag\\
& \ \ \ \  +  \frac{1}{2\pi i m t}\int_{\Psi(R)}^{E_0} e^{-2\pi m i\omega(E) t} \frac{d}{dE}\left(\frac{T(E)G_{f_0}(m,E)}{\omega'(E)}\right) \diff E\notag\\
& = \frac{1}{2\pi i m t}\int_{\Psi(R)}^{E_0} e^{-2\pi m i\omega(E) t} \frac{d}{dE}\left(\frac{T(E)G_{f_0}(m,E)}{\omega'(E)}\right) \diff E. \label{E:RLV1}
\end{align}
Notice that the boundary terms in the second equality vanishes due to
\[
\fhat_0(m,E_0) =|\varphi'(E)|\widehat{g_0}(m,E)= 0
\]
as $k\geq 2$ and $\theta(R,\Psi(R)) = \frac12$ by~\eqref{E:WVANISH}.
Integrating-by-parts again we  arrive at
\begin{align}
L(t,m;R)  &= \frac1{4\pi^2 m^2 t^2} \left[e^{-2\pi m i\omega(E) t} \frac1{\om'(E)}\frac{d}{dE}\left(\frac{T(E)G_{f_0}(m,E)}{\omega'(E)}\right) \right]^{E=E_0}_{E=\Psi(R)} \notag\\
& \ \ \ \ -  \frac1{4\pi^2 m^2 t^2} \int_{\Psi(R)}^{E_0} e^{-2\pi m i\omega(E) t} \frac d{dE}\left(\frac1{\om'(E)}\frac{d}{dE}\left(\frac{T(E)G_{f_0}(m,E)}{\omega'(E)}\right)\right) \diff E.\label{E:LTMRDECAY}
\end{align}
Note that the lower boundary term vanishes by using the formula~\eqref{E:EXTERIOR0} and 
the equality $\Theta(R,\Psi(R))=\frac12$ since $R>r_\ast$. On the other hand, if $k=2$, the upper boundary term is bounded, while it vanishes if $k>2$ as $\pa_E f_0=-\varphi''g_0+\varphi'\pa_E g_0$, and $\widehat{g_R}$ is a bounded function.

 From the assumptions on $\om(E)=T(E)^{-1}$, the identity $G_{f_0}=|\varphi'(E)|G_{g_0}$, and crucially Lemma~\ref{L:HINTEGRABLE}, we conclude that the argument inside the integral on the right-most side of~\eqref{E:LTMRDECAY} is bounded by
 \begin{multline*}
 	\bigg|\frac d{dE}\left(\frac1{\om'(E)}\frac{d}{dE}\left(\frac{T(E)G_{f_0}(m,E)}{\omega'(E)}\right)\right)\bigg|\\
 	\leq \frac C{|m|} (E-\Emin)^{-\frac12}(E-\Psi(R))^{-\frac12}+\frac{C}{|m|^2}(E_0-E)^{k-3},
 \end{multline*}
 where the second term is present only if $k>2$ due the the vanishing of $\frac{\dif^3}{\dif E^3}\varphi'(E)$ if $k=2$.
In particular, this is integrable since $R\neq r_*$. In fact, the following bound holds
\begin{equation*}
\int_{\Psi(R)}^{E_0} (E-\Emin)^{-\frac12}(E-\Psi(R))^{-\frac12} \diff E \le C \left(1+\lv \log (\Psi(R)-\Emin)\rv\right),
\end{equation*}
which can be easily verified by a direct calculation.  Since $\Psi(R) = \Emin + \frac12\Psi''(r_\ast)(R-r_\ast)^2 + O(|R-r_\ast|^3)$ and $\frac12\Psi''(r_\ast)\neq 0$ (recall that $\Emin$ is a non-degenerate elliptic point), 
we conclude that there exists a constant $C$ such that 
\begin{align}\label{E:LTMR1}
\lv L(t,m;R)\rv \le  C\,\frac{1+\lv\log|R-r_\ast|\rv}{|m|^3\,t^2}, \qquad m\in \Z^\ast,\, R\in(r_\ast,\Rmax],\, t\in\R^\ast.
\end{align}

 If $R\in[\Rmin,r_\ast)$, we repeat the above strategy, but keep in mind that the integration domain is all of $I$ and it splits naturally into regions 
 $I = (\Emin,\Psi(R)]\cup(\Psi(R),E_0]$. As in~\eqref{E:RLV1}, we have 
\begin{align}
	& L(t,m;R)  \notag\\
	& = -\frac{1}{2\pi i m t}\left(\int_{\Emin}^{\Psi(R)}+\int_{\Psi(R)}^{E_0}\right) \frac{d}{dE}\left(e^{-2\pi m i\omega(E) t}\right) \frac{ T(E)}{\omega'(E)} G_{f_0}(m,E) \diff E\notag\\
	& = \frac{-1}{2\pi i m t} \left[e^{-2\pi m i\omega(E) t} \frac{G_{f_0}(m,E)T(E)}{\omega'(E)}\right]^{E=E_0}_{E=\Psi(R)} 
	- \frac{1}{2\pi i m t} \left[e^{-2\pi m i\omega(E) t} \frac{G_{f_0}(m,E)T(E)}{\omega'(E)}\right]^{E=\Psi(R)}_{E=\Emin}\notag\\  
	& \ \ \ \ +  \frac{1}{2\pi i m t}\left(\int_{\Emin}^{\Psi(R)}+\int_{\Psi(R)}^{E_0}\right) e^{-2\pi m i\omega(E) t} \frac{d}{dE}\left(\frac{T(E)G_{f_0}(m,E)}{\omega'(E)}\right) \diff E\notag\\
	& = \frac{1}{2\pi i m t}\left(\int_{\Emin}^{\Psi(R)}+\int_{\Psi(R)}^{E_0}\right) e^{-2\pi m i\omega(E) t} \frac{d}{dE}\left(\frac{T(E)G_{f_0}(m,E)}{\omega'(E)}\right) \diff E. \label{E:INTERIOR}
\end{align}
Here, the boundary term vanishes due to $\fhat(m,E_0)=\fhat(m,\Emin)=0$ (recall $m\in\Z^\ast$ and Lemma~\ref{lemma:Eminvanishing}), and the two boundary terms at $E=\Psi(R)$ appear with
exactly the opposite sign. We now run the integration-by-parts argument again and obtain just like above
\begin{align}
 L(t,m;R)
& =  \frac1{4\pi^2 m^2 t^2} \left[e^{-2\pi m i\omega(E) t} \frac1{\om'(E)}\frac{d}{dE}\left(\frac{T(E)G_{f_0}(m,E)}{\omega'(E)}\right) \right]^{E=E_0}_{E=\Psi(R)} \notag\\
& \ \ \ \  +  \frac1{4\pi^2 m^2 t^2} \left[e^{-2\pi m i\omega(E) t} \frac1{\om'(E)}\frac{d}{dE}\left(\frac{T(E)G_{f_0}(m,E)}{\omega'(E)}\right) \right]^{E=\Psi(R)}_{E=\Emin}\notag\\
 & \ \ \ \ -  \frac1{4\pi^2 m^2 t^2} \int_{\Psi(R)}^{E_0} e^{-2\pi m i\omega(E) t} \frac d{dE}\left(\frac1{\om'(E)}\frac{d}{dE}\left(\frac{T(E)G_{f_0}(m,E)}{\omega'(E)}\right)\right) \diff E \notag\\
& \ \ \ \ -  \frac1{4\pi^2 m^2 t^2} \int_{\Emin}^{\Psi(R)} e^{-2\pi m i\omega(E) t} \frac d{dE}\left(\frac1{\om'(E)}\frac{d}{dE}\left(\frac{T(E)G_{f_0}(m,E)}{\omega'(E)}\right)\right) \diff E \notag\\
 & = - \frac1{4\pi^2 m^2 t^2} \left[e^{-2\pi m i\omega(E) t} \frac1{\om'(E)}\frac{d}{dE}\left(\frac{T(E)G_{f_0}(m,E)}{\omega'(E)}\right) \right]_{E=\Emin} \notag\\
& \ \ \ \ -   \frac1{4\pi^2 m^2 t^2} \int_{\Psi(R)}^{E_0} e^{-2\pi m i\omega(E) t} \frac d{dE}\left(\frac1{\om'(E)}\frac{d}{dE}\left(\frac{T(E)G_{f_0}(m,E)}{\omega'(E)}\right)\right) \diff E \notag\\
& \ \ \ \  -  \frac1{4\pi^2 m^2 t^2} \int_{\Emin}^{\Psi(R)} e^{-2\pi m i\omega(E) t} \frac d{dE}\left(\frac1{\om'(E)}\frac{d}{dE}\left(\frac{T(E)G_{f_0}(m,E)}{\omega'(E)}\right)\right) \diff E.
 \label{E:INTERIOR1}
\end{align}
Here we see that, unlike~\eqref{E:LTMRDECAY}, not all the boundary terms vanish. The boundary term at $E=\Emin$ survives  and 
 we use the bound~\eqref{E:PARTIALEHBOUNDLOW} to infer that this boundary term is bounded by $\frac C{|m|^3 t^2}$. The remaining boundary terms
 do indeed vanish if $k>2$ or are bounded if $k=2$; at $E=E_0$ we use $\pa_E\fhat_0(m,E_0)=k\eps(E_0-E)^{k-2}\widehat{g_0}(m,E)$. At $E=\Psi(R)$ the boundary terms appear with the opposite sign and the function $\frac d{dE}H(m,\cdot)$
is  continuous at $E =\Psi(R)>\Emin$ due to~\eqref{E:EXTERIOR0} and $\Theta(R,\Psi(R))=0$ for $R\in[\Rmin,r_\ast)$.

Finally, just like above, using~\eqref{E:PARTIALEHBOUND}--\eqref{E:PARTIALEEHBOUNDLOW} we conclude 
 \begin{align}\label{E:LTMR2}
 \lv L(t,m;R)\rv \le  C\,\frac{1+\lv\log|R-r_\ast|\rv}{|m|^3\,t^2}, \qquad m\in \Z^\ast,\,R\in[\Rmin,r_\ast),\, t\in\R^\ast.
 \end{align} 
 Since $\pa_tU_f =4\pi^2\sum_{m\in\Z^\ast} L(t,m;R)$, it follows from~\eqref{E:LTMR1} and~\eqref{E:LTMR2} that there exists a constant $C$ such that~\eqref{E:PARTIALTUFDECAYTRAPPED} holds.
 
{\em Proof of~\eqref{E:PARTIALTUFDECAYEXTERIOR}.} 
To obtain the improved bound~\eqref{E:PARTIALTUFDECAYEXTERIOR}, we suppose $k\geq \frac52$ and we rely on the following key observation. As a consequence of Lemma~\ref{L:HINTEGRABLE}, the integrand on the right-most side of~\eqref{E:LTMRDECAY} features an inverse square root singularity of the form
$(E-\Psi(R))^{-\frac12}$ at the leading order, near $E = \Psi(R)$. More precisely, since $T\in C^3(\bar I)$ the most singular term in the integrand is of the form 
\begin{equation*}
e^{-2\pi m i\omega(E) t} \frac{T(E)}{\om'(E)^2}\frac{\sqrt{E-\Psi(R)}\frac{d^2}{dE^2}G_{f_0}(m,E)}{\sqrt{E-\Psi(R)}},
\end{equation*}
where the absolute value of the numerator $|\sqrt{E-\Psi(R)}\frac{d^2}{dE^2}G_{f_0}(m,E)|$ is bounded by $\frac{C}{|m|\sqrt{E-\Emin}}+\frac{C}{|m|^2}(E_0-E)^{k-3}$ by Lemma~\ref{L:HINTEGRABLE} and therefore non-singular when $R\neq r_\ast$.
To establish the higher decay, we again employ a cut-off argument as in the proof of Proposition~\ref{P:PARTIALRUFDECAY}, using now the key estimate $$|\frac{d}{dE}(\frac{d^2}{dE^2}G_{f_0}(m,E))|\leq C(E-\Emin)^{-1}(E-\Psi(R))^{-1}+C(E_0-E)^{k-4}$$ from Lemma~\ref{L:HIGHORDER}. Comparing to the analogous estimate in the proof of Proposition~\ref{P:PARTIALRUFDECAY}, we see that we have an additional singular factor of $(E-\Emin)^{-\frac12}$ supported on the region $E>\Psi(R)$. This precisely contributes the additional singular term $|R-r_*|^{-1}$ in the final bound.
 \end{proof}

%%%%%%%%%%%%%%%%%%%%%%%%%%%%%%%
%%%%%%%%%%%%%%%%%%%%%%%%%%%%%%%

\section{Decay in the radially symmetric case}\label{S:MAIN2}

The purpose of this section is to prove Theorem~\ref{T:MAIN2} using the results from Section~\ref{S:AAREG}, in particular, Section~\ref{S:AAREG3D}. 
We first derive the following additional lemma concerning the integrability properties of the Green's function $h_R(\theta,E,L)$.

\begin{lemma}\label{L:hRweightedL2}
 For any $R\in[\Rmin,\Rmax]$,  the Green's function $h_R$ satisfies the weighted $L^2$ bound
\beq
\int_{J}\int_{\mathbb{S}^1}\Big|\frac{h_R(\theta,E,L)^2}{E-\EminL}\Big| \,\dif \th\,\dif E\,\dif L\leq C.\nonumber
\eeq 
\end{lemma}

\begin{proof}
As the estimate is most singular when $R$ is a trapping radius, i.e.,~$R=r_{L_*}$ for some $L_*\in[L_0,L_{\max}]$, we focus on this case. Radially Taylor expanding the effective potential yields
\begin{equation*}
	\Psi_L(R)=\EminL+\frac{\alpha_L^2}2\,(R-r_L)^2+\mathcal O(R-r_L)^3.
\end{equation*}
Because $L\mapsto r_L$ is linear~\cite[(2.2.61)]{St23}, it hence follows that there exists a constant $C>0$ depending on~$L_0$ and~$\Lmax$ such that
\begin{equation}\label{I:PSILDIFF}
	\frac1C|L-L_*|^2\leq \big|\Psi_L(R)-\Psi_L(r_L)\big|\leq C |L-L_*|^2,\qquad L\in[L_0,\Lmax].
\end{equation}

We recall that $h_R$ is identically $1$ on the interval $\th\in[-\theta(R,E,L),\theta(R,E,L)]$ when $E\geq \Psi_L(R)$ and vanishes otherwise. Thus
\begin{align*}
	\int_{J}\int_{\mathbb{S}^1}\Big|\frac{h_R(\theta,E,L)^2}{E-\EminL}\Big| \,\dif \th\,\dif E\,\dif L&=\int_{L_0}^{L_{\max}}\int_{\Psi_L(R)}^{E_0}2\theta(R,E,L)(E-\EminL)^{-1}\,\dif E\,\dif L\\
	&\leq  C+ \int_{L_0}^{L_{\max}}\big|\log(\Psi_L(R)-\EminL)\big|\,\dif L\leq C,
\end{align*}
where we have used that $0\leq\theta(R,E,L)\leq\frac12$  and~\eqref{I:PSILDIFF} to see that the final integral is finite.
\end{proof}

{\em Proof of Theorem~\ref{T:MAIN2}.}
We argue in the case $\ell>0$, $k>1$, as the remaining cases are somewhat simpler.

{\em Proof of~\eqref{E:DECAYL}.}
To estimate $\pa_R U_f$ we use the representation formula
\begin{align}
R^2 \pa_RU_{f}(t,R)  & = 4\pi^2 \int_{L_0}^{L_{\max}}\int_{\EminL}^{E_0}\int_{\mathbb S^1} h_R(r(\theta,E,L))\, f T(E,L) \diff\theta\diff E\,\diff L \notag\\
& = 4\pi^2 \sum_{m\in\Z^\ast} \int_{L_0}^{L_{\max}}\int_{\EminL}^{E_0} \fhat(t,m,E,L) \overline{\widehat{h_R}(m,E,L)} T(E,L)\diff E\,\diff L \notag \\
& = 4\pi^2 \sum_{m\in\Z^\ast} \int_{L_0}^{L_{\max}} \int_{\Psi_L(R)}^{E_0} e^{-2\pi mi \om(E,L) t}  \fhat_0(m,E,L) \widehat{h_R}(m,E,L)  T(E,L)\diff E\,\diff L, \notag\\
& = 4\pi^2  \sum_{m\in\Z^\ast} \int_{L_0}^{L_{\max}} \int_{\om_R(L)}^{\om_0(L)} e^{-2\pi mi q t}  \widehat{h_R}(m,q,L) \frac{\fhat_0(m,q,L)T}{\om'} \diff q\,\diff L\label{E:NSP2}
\end{align}
where we recall~\eqref{E:FHATREPL} and Lemma~\ref{L:GHFOURIER}.
 Notice that the $m=0$ terms vanish because $f_0\perp\ker(\D)$ and recall Remark~\ref{R:KERD}. 

%{\em Proof of part (a).}
As the greatest restriction on the decay rate occurs when $R$ is a trapping radius, i.e.,~$R=r_{L_\ast}$ for some $L_\ast\in[L_0,\Lmax]$, we prove the decay bound for this range of $R$. In the remaining case, the proof is considerably easier.

We first apply Lemma~\ref{L:RECURSION1} (specifically~\eqref{E:INTID1}) and obtain
\begin{multline}
\int_{L_0}^{L_{\max}} \int_{\Psi_L(R)}^{E_0} e^{-2\pi mi \omega(E,L) t} \widehat{h_R}(m,E,L) \frac{\fhat_0 T}{\omega'} \diff E\diff L\\
 = -\frac{1}{2\pi i m t} \int_{L_0}^{L_{\max}} \int_{\om_R(L)}^{\om_0(L)} e^{-2\pi mi q t}\Bigg[ \widehat{h_R}  \A_1(\frac{\fhat_0 T}{\om'})
+ w\cos(2\pi m \theta) \A_2(\frac{\fhat_0 T}{\om'})\Bigg] \diff q \diff L.\label{E:DRUFDECAYELFIRST}
\end{multline}
We now introduce two cut-off functions, one in $E$ and one in $L$ as follows. Let $0<\de\ll1$ to be chosen later. In the case that $L_*-L_0<2\de$, we define a cut-off $\phi_\de\in C_c^\infty((L_0-1,L_{\max}+1))$ such that $\phi_\de=0$ on $B_{4\de}(L_\ast)$, $\phi_\de=1$  on $[L_0,L_{\max}]\setminus B_{8\de}(L_\ast)$, $0\leq \phi_\de\leq 1$, $|\phi_\de'(L)|\leq\frac1\de$. Note that such a cut-off vanishes on the interval of radius $2\de$ around $L_0$ also.

In the case that $L_*-L_0\geq 2\de$, we instead take $\phi_\de\in C_c^\infty((L_0-1,L_{\max}+1))$ such that $\phi_\de=0$ on $B_{\de/2}(L_\ast)\cup B_{\de/2}(L_0)$, $\phi_\de=1$ on $[L_0,L_{\max}]\setminus (B_{\de}(L_\ast)\cup B_{\de}(L_0))$  and again $0\leq\phi_\delta\leq1$ and $|\phi_\delta'|\leq\frac3\delta$. 

Finally, we introduce the cut-off $\psi_\de\in C_c^\infty((E_{\min}^{L_0}-1,E_0+1))$ such that $\psi_\de =0$ on $B_{\de}(E_0)$, $\psi_\de=1$ on $[E^{L_0}_{\min},E_0]\setminus B_{2\de}(E_0)$, $0\leq \psi_\de\leq 1$, $|\psi_\de'|\leq\frac2\de$.

We split the integral on the right hand side of~\eqref{E:DRUFDECAYELFIRST} as 
\begin{align}
&\int_{L_0}^{L_{\max}} \int_{\Psi_L(R)}^{E_0} e^{-2\pi mi \omega(E,L) t} \widehat{h_R}(m,E,L) \frac{\fhat_0 T}{\omega'} \diff E\diff L \notag\\
& = -\frac{1}{2\pi i m t} \int_{L_0}^{L_{\max}} \int_{\om_R(L)}^{\om_0(L)} \phi_\de(L)\psi_\de(E)e^{-2\pi mi q t}\Bigg[ \widehat{h_R}  \A_1(\frac{\fhat_0 T}{\om'})
+ w\cos(2\pi m \theta) \A_2(\frac{\fhat_0 T}{\om'})\Bigg] \diff q \diff L\notag\\
&\ \ \ -\frac{1}{2\pi i m t} \int_{L_0}^{L_{\max}} \int_{\om_R(L)}^{\om_0(L)} (1-\phi_\de(L)\psi_\de(E))e^{-2\pi mi q t}\Bigg[ \widehat{h_R}  \A_1(\frac{\fhat_0 T}{\om'})
+ w\cos(2\pi m \theta) \A_2(\frac{\fhat_0 T}{\om'})\Bigg] \diff q \diff L\notag\\
&=\colon \tilde I(m)+\tilde J(m). \label{E:IJDEF}
\end{align}
Treating the terms $\tilde I$ and $\tilde J$ separately, we first integrate by parts once more in $\tilde I$, using the definition of the cut-offs to see the vanishing condition at $L_0$ and $E_0$, respectively, and Lemma~\ref{L:RECURSION1}, and obtain
\begin{multline*}
\tilde I(m)=\frac{1}{(2\pi i m t)^2} \int_{L_0}^{L_{\max}} \int_{\om_R(L)}^{\om_0(L)}\,e^{-2\pi mi q t}\Bigg[ \widehat{h_R}  \Big(\A_1\Big(\phi_\de\psi_\de\A_1(\frac{\fhat_0 T}{\om'})\Big)+\B_1\Big(\phi_\de\psi_\de\A_2(\frac{\fhat_0 T}{\om'})\Big)\Big)\\
+ w\cos(2\pi m \theta) \Big(\A_2\Big(\phi_\de\psi_\de\A_1(\frac{\fhat_0 T}{\om'})\Big)+ \B_2\Big(\phi_\de\psi_\de\A_2(\frac{\fhat_0 T}{\om'})\Big)\Big)\Bigg] \diff q \diff L.
\end{multline*}

Distributing derivatives, applying~\eqref{E:FZEROCONTROL} (with $N=2$, $j=0,1,2$), and applying Fubini's theorem to return to $(E,L)$-coordinates, this yields the bound
\beqa
|\tilde I(m)|\leq&\,\frac{C}{m^2t^2}\bigg[\int_{L_0}^{L_{\max}} \int_{\Psi_L(R)}^{E_0}\big(|\widehat{h_R}|+|w|\big)|\phi_{\de}'(L)|\psi_\de(E)\Big(\Big|\A_1\big(\frac{\fhat_0 T}{\om'}\big)\Big|+\Big|\A_2\big(\frac{\fhat_0 T}{\om'}\big)\Big|\Big)\,\dif E\,\dif L\\
&+\int_{L_0}^{L_{\max}} \int_{\Psi_L(R)}^{E_0}\big(|\widehat{h_R}|+|w|\big)\phi_{\de}(L)|\psi_\de'(E)|\Big(\Big|\A_1\big(\frac{\fhat_0 T}{\om'}\big)\Big|+\Big|\A_2\big(\frac{\fhat_0 T}{\om'}\big)\Big|\Big)\,\dif E\,\dif L\\
&+\int_{L_0}^{L_{\max}} \int_{\Psi_L(R)}^{E_0}\phi_\de(L)\psi_\de(E)|\widehat{h_R}|\Big|\P_2\big(\frac{\fhat_0 T}{\om'}\big)\Big|\dif E\,\dif L\\
&+\int_{L_0}^{L_{\max}} \int_{\Psi_L(R)}^{E_0}\phi_\de(L)\psi_\de(E)|w|\Big|\Q_2\big(\frac{\fhat_0 T}{\om'}\big)\Big|\dif E\,\dif L\bigg]\\
=\colon&\, \tilde I_1+\tilde I_2+\tilde I_3+\tilde I_4, \nonumber
\eeqa
where we have applied the Leibniz rule to distribute the derivative contributions from the operators $\A_1$, $\A_2$, $\B_1$ and $\B_2$ and used that $Q_1$ and $Q_2$ are bounded functions by Lemma~\ref{L:PAQHR}.

Before bounding $\tilde I_1$--$\tilde I_4$, we first give a useful auxiliary estimate on a key integral:
\beqa\label{INEQ:AUX}
\int_{L_0}^{L_{\max}}\phi_\de(L)\big(\Psi_L(R)-\Psi_L(r_L)\big)^{-\frac12}\,\dif L\leq C\int_{L_0}^{L_{\max}}\phi_\de(L)|L-L_*|^{-1}\,\dif L\leq C|\log\de|,
\eeqa
by~\eqref{I:PSILDIFF} and the definition of $\phi_\de$.

 Now, to bound $\tilde I_1$, we recall that the support of $\phi_\de'(L)$ is  at most  of size $8\de$, while as $k>1$, all of the $E$-integrands are integrable, so that, applying~\eqref{E:P1FBOUNDS}--\eqref{E:Q1FBOUNDS},
 \begin{align}
 	|\tilde I_1|\leq&\frac{C}{m^2t^2}+\frac{C}{m^2t^2} \int_{L_0}^{L_{\max}}\int_{\Psi_L(R)}^{E_0}|\phi_\de'(L)|\psi_\de(E)\Big(\frac{1}{|m|}(E-E^L_{\min})^{-\frac12}\nonumber\\&\hspace*{15.5em}+\frac{1}{m^2}(E_0-E)^{k-2}+\frac1{m^2}(L-L_0)^{\ell-1}\Big)\,\dif E\,\dif L\nonumber\\
 	\leq&\frac{C}{m^2t^2}+\frac{C}{m^2t^2}\int_{L_0}^{L_{\max}}|\phi_\de'(L)|\big(1+(L-L_0)^{\ell-1}\big)\,\dif L\leq C\frac{(1+\de^{\ell-1})}{m^2t^2},\label{ineq:I1}
 \end{align}
where we note that the support of $\phi_\de'$ contains an annulus (pair of intervals) centred on $L_\ast$, which may or may not be close to $L_0$.
 The integral containing $(E-E^L_{\min})^{-\frac12}$ has been estimated similarly to Lemma~\ref{L:hRweightedL2}. 

To bound $\tilde I_2$, we recall that  the support of $\psi_\de'$ is of size $\leq2\de$, that on the support we have $E_0-E>\de$, and that $|\psi_\de'|\leq\frac2\de$, so that
\begin{align}
	|\tilde I_2|&\leq \frac{C}{m^2t^2}+\frac{C}{m^2t^2} \int_{L_0}^{L_{\max}}\int_{\Psi_L(R)}^{E_0}\phi_\de(L)|\psi_\de'(E)|\Big(\frac{1}{|m|}(E-E^L_{\min})^{-\frac12}\nonumber\\&\hspace*{17em}+\frac{1}{m^2}(E_0-E)^{k-2}+\frac1{m^2}(L-L_0)^{\ell-1}\Big)\,\dif E\,\dif L\nonumber\\
	&\leq\frac{C}{m^2t^2}+\frac{C}{m^2t^2} \int_{L_0}^{L_{\max}}\phi_\de(L)\Big((\Psi_L(R)-\Psi_L(r_L))^{-\frac12}+(L-L_0)^{\ell-1}+\de^{k-2}\Big)\,\dif L\leq\nonumber\\ &\leq C\frac{|\log\de|+\de^{k-2}}{m^2t^2},\label{ineq:I2}
\end{align}
where we have bounded $(E-E^L_{\min})^{-\frac12}\le (\Psi_L(R)-\Psi_L(r_L))^{-\frac12}$ and then applied~\eqref{INEQ:AUX}.

Next, to bound $\tilde I_3$, we observe the following  simple integral estimate: 
\begin{equation}
\int_{E_{\min}^L}^{E_0-\de}(E-E^L_{\min})^{-\frac12}(E_0-E)^{k-2}\,\dif E\leq C(1+\delta^{k-\frac32}).
 \label{E:INTSCALE}
 \end{equation}
 Thus, estimating $\P_2$ with~\eqref{E:P2FBOUNDS} (and  neglecting superfluous powers of $m$), we have
 \beqa
 |\tilde I_3|\leq&\,\frac{C}{m^2t^2}\int_{L_0}^{L_{\max}} \int_{\Psi_L(R)}^{E_0}\phi_\de(L)\psi_\de(E)|\widehat{h_R}|\Big((E-E^L_{\min})^{-\frac32}+(E-E_{\min}^L)^{-\frac12}\big((L-L_0)^{\ell-1}+(E_0-E)^{k-2}\big)\\
&\hspace{30mm}+(L-L_0)^{\ell-2}+(L-L_0)^{\ell-1}(E_0-E)^{k-2}+(E_0-E)^{k-3}\\
&\hspace{30mm}+|m||\varphi'|(E-\EminL)^{-\frac12}+m^2|\varphi'|(E-\EminL)^{-\frac12}\Big(\frac{|m^2\widehat{g_0}(m,E,L)|}{(E-\EminL)^{\frac12}}\Big)\Big)\dif E\,\dif L.
\eeqa
We estimate the first two lines of the right hand side first, using the estimate $|\widehat{h_R}|\leq 1$ to bound
\beqa\label{ineq:I3I}
{}&\frac{C}{m^2t^2}\int_{L_0}^{L_{\max}} \int_{\Psi_L(R)}^{E_0}\phi_\de(L)\psi_\de(E)|\widehat{h_R}|\Big((E-E^L_{\min})^{-\frac32}+(E-E_{\min}^L)^{-\frac12}\big((L-L_0)^{\ell-1}+(E_0-E)^{k-2}\big)\\
&\hspace{30mm}+(L-L_0)^{\ell-2}+(L-L_0)^{\ell-1}(E_0-E)^{k-2}+(E_0-E)^{k-3}\Big)\,\dif E\,\dif L\\
& \leq\frac{C}{m^2t^2}\int_{L_0}^{L_{\max}}\phi_\de(L)\Big(\big(\Psi_L(R)-\Psi_L(r_L)\big)^{-\frac12}+(L-L_0)^{\ell-1}+1+\de^{k-\frac32}+(L-L_0)^{\ell-2}+\de^{k-2}\Big)\,\dif L\\
 &\leq C\frac{|\log\de|+\de^{k-2}+\de^{\ell-1}}{m^2t^2},
 \eeqa
 where we have applied~\eqref{INEQ:AUX} and~\eqref{E:INTSCALE}. To estimate the final terms in $\tilde I_3$, we employ the improved estimate $|\widehat{h_R}|\leq \frac1{|m|}$ (recall $\widehat{h_R}(m,E,L)=\frac{\sin(2\pi m \theta(R,E,L))}{\pi m}$) in the first term and H\"older's inequality in the second to see
 \begin{multline}
 	\frac{C}{m^2t^2}\int_{L_0}^{L_{\max}} \int_{\Psi_L(R)}^{E_0}\phi_\de(L)\psi_\de(E)|\widehat{h_R}|\Big(\frac{|m||\varphi'|}{(E-\EminL)^{\frac12}}+\frac{m^2|\varphi'|}{(E-\EminL)^{\frac12}}\,\frac{|m^2\widehat{g_0}(m,E,L)|}{(E-\EminL)^{\frac12}}\Big)\dif E\,\dif L\\
 	\leq\frac{C}{m^2t^2}+\frac{C}{t^2}\Big\|\frac{\widehat{h_R}(m,\cdot,\cdot)}{(E-\EminL)^{\frac12}}\Big\|_{L^2(J)}\Big\|\frac{\widehat{\pa^2_\th g_0}(m,\cdot,\cdot)|}{(E-\EminL)^{\frac12}}\Big\|_{L^2(J)}.\label{ineq:I3II}
 \end{multline}
 In order to bound $\tilde I_4$, we observe by comparing~\eqref{E:P2FBOUNDS} to~\eqref{E:Q2FBOUNDS} that  most terms in~$\tilde I_4$ can also be bounded using the estimate~\eqref{ineq:I3I}  with the exception of one term which must be bounded separately:
 \begin{align}
 	{}& \frac{C}{m^2t^2}\int_{L_0}^{L_{\max}} \int_{\Psi_L(R)}^{E_0}\phi_\de(L)\psi_\de(E)|\varphi'|\Big|m\pa^{(1)}\Big(\frac{m\widehat{g_0}(m,E,L)}{(E-\EminL)^{\frac12}}\Big)\Big|\,\dif E\,\dif L\nonumber\\
 	&\leq \frac{C}{|m|t^2}\Big(\int_{L_0}^{L_{\max}}\big|\log(\Psi_L(R)-\Psi_L(r_L))\big|\,\dif L\Big)^{\frac12}\Big\|(E-\EminL)^{\frac12}\pa^{(1)}\Big(\frac{m\widehat{g_0}(m,E,L)}{(E-\EminL)^{\frac12}}\Big)\Big\|_{L^2(J)}\nonumber\\
 	&\leq \frac{C}{|m|t^2}\Big\|(E-\EminL)^{\frac12}\pa^{(1)}\Big(\frac{m\widehat{g_0}(m,E,L)}{(E-\EminL)^{\frac12}}\Big)\Big\|_{L^2(J)};\label{ineq:I4}
 \end{align}
 see the proof of Lemma~\ref{L:hRweightedL2} for the latter estimate.

 Thus,
 combining the estimates~\eqref{ineq:I1}--\eqref{ineq:I4}, we have obtained
 \begin{multline*}
 	|\tilde I(m)|\leq C\frac{|\log\de|+\de^{k-2}+\de^{\ell-1}}{m^2t^2}+\frac{C}{t^2}\Big\|\frac{\widehat{h_R}(m,\cdot,\cdot)}{(E-\EminL)^{\frac12}}\Big\|_{L^2(J)}\Big\|\frac{\widehat{\pa^2_\th g_0}(m,\cdot,\cdot)|}{(E-\EminL)^{\frac12}}\Big\|_{L^2(J)}\\
 	+\frac{C}{|m|t^2}\Big\|(E-\EminL)^{\frac12}\pa^{(1)}\Big(\frac{m\widehat{g_0}(m,E,L)}{(E-\EminL)^{\frac12}}\Big)\Big\|_{L^2(J)}.
 \end{multline*}
 Summing in $m$ and applying the Plancherel theorem, we therefore obtain
 \beqa
 \sum_{m\in\Z_*}|\tilde I(m)|\leq&\, \frac{C}{t^2}\bigg(\sum_{m\in\Z_*}\frac{|\log\de|+\de^{k-2}+\de^{\ell-1}}{m^2}\\
 &+\Big(\sum_{m\in\Z_*}\Big\|\frac{\widehat{h_R}(m,\cdot,\cdot)}{(E-\EminL)^{\frac12}}\Big\|_{L^2(J)}^2\Big)^{\frac12}\Big(\sum_{m\in\Z_*}\Big\|\frac{\widehat{\pa^2_\th g_0}(m,\cdot,\cdot)|}{(E-\EminL)^{\frac12}}\Big\|_{L^2(J)}^2\Big)^{\frac12}\\
 &+\Big(\sum_{m\in\Z_*}\frac{1}{m^2}\Big)^{\frac12}\Big(\sum_{m\in\Z_*}\Big\|(E-\EminL)^{\frac12}\pa^{(1)}\Big(\frac{m\widehat{g_0}(m,E,L)}{(E-\EminL)^{\frac12}}\Big)\Big\|_{L^2(J)}^2\Big)^{\frac12}\bigg)\\
 \leq&\,C\frac{|\log\de|+\de^{k-2}+\de^{\ell-1}}{t^2}, \label{E:ITILDEBOUND}
 \eeqa
 where we have applied Lemma~\ref{L:hRweightedL2} and used the estimates
 \beqa
 \Big(\sum_{m\in\Z_*}\Big\|(E-\EminL)^{\frac12}\pa^{(1)}\Big(\frac{m\widehat{g_0}(m,E,L)}{(E-\EminL)^{\frac12}}\Big)\Big\|_{L^2(J)}^2\Big)^{\frac12}&\leq C\|g_0\|_{W^{2,\infty}(\Om)},\\
 \Big(\sum_{m\in\Z_*}\Big\|\frac{\widehat{\pa^2_\th g_0}(m,\cdot,\cdot)|}{(E-\EminL)^{\frac12}}\Big\|_{L^2(J)}^2\Big)^{\frac12}&\leq C\|g_0\|_{W^{2,\infty}(\Om)},\nonumber
 \eeqa
 which both follow from the Plancherel theorem and the regularity of $r$ and $w$ as functions of $E,L$, as given by Lemma~\ref{L:RELCONTROL}; see Lemma~\ref{L:IMPROVED} for similar arguments.
 
 Turning now to the integral $\tilde J(m)$ from~\eqref{E:IJDEF}, we simply estimate using~\eqref{E:P1FBOUNDS}--\eqref{E:Q1FBOUNDS} with $N=2$, $j=0,1$, and split the domain of integration into $A_1\cup A_2$, where
 \begin{equation*}
 	A_1\coloneqq\{(E,L)\,|\,E\in(\Psi_L(R),E_0),\,L\in(L_0,L_{\max})\}\cap \{\phi_\de(L)=1\},
 \end{equation*}
 and
 \begin{equation*}
 	A_2\coloneqq\{(E,L)\,|\,E\in(\Psi_L(R),E_0),\,L\in(L_0,L_{\max})\}\cap\{\phi_\de(L)<1\}.
 \end{equation*}
 Then, applying the estimates~\eqref{E:P1FBOUNDS}--\eqref{E:Q1FBOUNDS} and $|\widehat{h_R}|\leq1$, we obtain
 \begin{multline*}
 	|\tilde J(m)|\leq\frac{C}{m^2t}\Big(\int_{A_1}+\int_{A_2}\Big)(1-\phi_\de(L)\psi_\de(E))\Big((E-E_{\min}^L)^{-\frac12}(E_0-E)^{k-1}(L-L_0)^{\ell}\\
 	+(E-E_{\min}^L)^{\frac12}\Big((E_0-E)^{k-2}(L-L_0)^{\ell}+(E_0-E)^{k-1}(L-L_0)^{\ell-1}+1\Big)\Big)\,\dif E\,\dif L.
 \end{multline*}
 For the integral on $A_1$, we observe that, on this set, the support of $1-\phi_\de\psi_\de$ is contained in the set $E_0-E<2\de$, so that
 \beqa
{}& \frac{C}{m^2t}\int_{A_1}(1-\psi_\de(E))\Big((E-E_{\min}^L)^{-\frac12}(E_0-E)^{k-1}(L-L_0)^{\ell}\\
 &\hspace{4em}+(E-E_{\min}^L)^{\frac12}\Big((E_0-E)^{k-2}(L-L_0)^{\ell}+(E_0-E)^{k-1}(L-L_0)^{\ell-1}+1\Big)\Big)\,\dif E\,\dif L\\
 &\leq \frac{C}{m^2t}\int_{L_0}^{L_{\max}}\de^{k-1}(L-L_0)^{\ell-1}\,\dif L\leq \frac{C\de^{k-1}}{m^2t}.\nonumber
 \eeqa
  In the same way, we uniformly estimate the $E$-integrals in the integral on~$A_2$ to arrive at 
 \begin{align}
{}& \frac{C}{m^2t}\int_{A_2}(1-\phi_\de(L)\psi_\de(E))\Big((E-E_{\min}^L)^{-\frac12}(E_0-E)^{k-1}(L-L_0)^{\ell}\notag\\
 &\hspace{4em}+(E-E_{\min}^L)^{\frac12}\Big((E_0-E)^{k-2}(L-L_0)^{\ell}+(E_0-E)^{k-1}(L-L_0)^{\ell-1}+1\Big)\Big)\,\dif E\,\dif L\notag\\
&\leq\frac{C}{m^2t} \int_{L_0}^{L_{\max}}(1-\phi_\delta(L))(1+(L-L_0)^{\ell-1})\,\dif L\leq \frac{C(\de+\de^{\ell})}{m^2t};\nonumber
\end{align}
 in the last estimate we used the structure of $\{\phi_\de\neq1\}$. 
Thus, upon summing over $m\in\Z_\ast$ we obtain
\begin{equation*}
	\sum_{m\in\Z_\ast}|\tilde J(m)|\leq \frac{C\de^{k-1}}{t}+ \frac{C(\de+\de^{\ell})}{t}.
\end{equation*}	
Choosing $\de=\frac1t$ above and in~\eqref{E:ITILDEBOUND}, we conclude
\beq
|\pa_R U_f(t,R)|\leq C\begin{cases}
\frac{1+|\log t|}{t^2}, & k\geq 2,\, \ell\geq 1,\\
\frac{1}{t^{\min\{\ell+1,k\}}}, &\text{else},\nonumber
\end{cases}
\eeq
thus proving~\eqref{E:DECAYL}.% part (a) of Theorem~\ref{T:MAIN2}.

{\em Proof of~\eqref{E:DECAYBETTER1}.}
When the data $f_0$ is supported strictly away from the trapping set $\pa J_{\text{trap}}$ the operators $\P_j$ and $\Q_j$ applied to $\frac{\fhat_0 T}{\om'}$
experience no singularities. The only regularity limitation can therefore come from the vacuum boundary $\pa J_{\text{vac}}$, but since the data is strictly supported away from $\pa J_{\text{vac}}$ and all the coefficients appearing in $\P_j,\Q_j$ are smooth in $J$, we may integrate-by-parts in~\eqref{E:NSP2} arbitrary many times to conclude~\eqref{E:DECAYBETTER1}.

%\bcr
{\em Proof of~\eqref{E:RHODECAYL}.}
For the purpose of proving decay of $\rho_f$, as in the $1+1$-dimensional case, it is sufficient to prove the decay of $\pa_t\rho_f$, due to the orthogonality assumption $f_0\in \text{{\em ker}}\,(\D)^\perp$. For convenience,
we work with the modified density
\begin{equation*}
\tilde\rho(t,r) \coloneqq \iint f \diff w\,\diff L = \frac{r^2}\pi \rho_f(t,r).
\end{equation*}
Then
\begin{align}
&	\pa_t\tilde\rho(t,r)  = - \int_{L_0}^{L_{\max}}\int \omega(E,L) \pa_\theta f \diff w\,\diff L \nonumber\\
	&= -  \int_{L_0}^{L_{\max}} \int_{\Psi_L(r)}^{E_0} \om(E,L)\,\frac{\pa_\theta f(t,\theta(r,E,L),E,L)+\pa_\theta f(t,1-\theta(r,E,L),E,L)}{w(r,E,L)} \diff E\,\diff L\eqqcolon\tilde\rho_{t1}+\tilde\rho_{t2},\label{E:TILDERHOPARTIALTL}
\end{align} 
where we have  used the change of variables $(w,L)\mapsto (E,L)$ for a given fixed $r$.
 To prove the  decay estimate~\eqref{E:RHODECAYL}, 
we focus on the first term~$\tilde\rho_{t1}$ on the right-hand side of~\eqref{E:TILDERHOPARTIALTL}.
The estimates for the second term~$\tilde\rho_{t2}$ follow in the same way by replacing the initial value~$f_0$ with $(\theta,E,L)\mapsto f_0(1-\theta,E,L)$.

We make a further change of variables, recalling the definition of $q=\om(E,L)$ from above and now setting $p=2(E-\Psi_L(r))$ (i.e.~$p=w^2$). The Jacobian matrices for the coordinate transformation are then
$$\frac{\pa(q,p)}{\pa(E,L)}=\begin{pmatrix}
\om' & \pa_L\om \\ 2 & -\frac{1}{r^2}
\end{pmatrix},\qquad \frac{\pa(E,L)}{\pa(q,p)}=\frac{1}{1+2\frac{r^2\pa_L\om}{\om'}}\begin{pmatrix}
\frac{1}{\om'} & \frac{r^2\pa_L\om}{\om'} \\ \frac{2r^2}{\om'} & -r^2
\end{pmatrix}.$$
From the smallness assumption on $\pa_L\om$, these are clearly invertible. We denote by $\mathcal K$ the domain $ J\cap\{E\geq \Psi_L(r)\}$.

We now recall the solution formula~\eqref{E:EXPLICITL}, \eqref{E:SEMIGROUPL} to infer that in the $(p,q)$ coordinates,
\begin{align}
\pa_q\left(f_0\circ T_t\right) = (\theta_q - t)\,\pa_\theta f_0\circ T_t + \pa_q f_0\circ T_t,\label{E:PEFCIRCT}
\end{align}
which yields the identity
\begin{align}\label{E:RLRHOL}
\pa_\theta f_0\circ T_t = \frac{-1}{t} \left(\pa_q\left(f_0\circ T_t\right) - \pa_\theta f_0\circ T_t \,\theta_q- \pa_q f_0 \circ T_t \right).
\end{align}
Note that the left-hand side of~\eqref{E:RLRHOL} is precisely $\pa_\theta f$ by~\eqref{E:EXPLICIT}
and therefore from~\eqref{E:TILDERHOPARTIALTL} and~\eqref{E:RLRHOL} we obtain
\begin{align}
\tilde\rho_{t1}(t,r) & = -  \int_{\mathcal K} q \frac{\pa_\theta f_0\circ T_t}{\sqrt{p}} \Big(\frac{\om'}{r^2}+2\pa_L\om\Big)^{-1}\diff q\,\diff p \notag\\
& = \frac1t \int_{\mathcal K} q \Big(\frac{\om'}{r^2}+2\pa_L\om\Big)^{-1} \frac{\pa_q\left(f_0\circ T_t\right) - \pa_\theta f_0\circ T_t \theta_q- \pa_qf_0 \circ T_t }{\sqrt{p}} \diff q\,\diff p.\label{E:TILDERHOPARTIALT2L}
\end{align}
We now introduce three cut-off functions. We introduce the cut-off parameter $\bar\de>0$ to be chosen later and let $\phi_{\bar\de}\in C_c^\infty([-1,E_0-\Psi_{L_0}(r)+1])$ be such that $0\leq \phi_{\bar \de}\leq 1$, 
\beqs
\phi_{\bar \de}(p)=\begin{cases}
0, & p\leq{\bar \de},\\
1, & 2{\bar \de}\leq  p\leq E_0-\Psi_{L_0}(r),
\end{cases}
\eeqs and with $|\phi_{\bar \de}'(p)\big|\leq\frac{C}{{\bar \de}}$.
Now we set 
$$F(p,q)=\frac{1}{\frac{\om'}{r^2}+2\pa_L\om}\in C^1(\overline{ \mathcal K}),$$
and re-write \eqref{E:TILDERHOPARTIALT2L} as 
\begin{align*}
\tilde\rho_{t1}(t,r) & =-\frac2t  \int_{\mathcal K} \phi_{\bar \de}(p) qF\frac{\pa_q\big(f_0\circ T_t\big) }{\sqrt{p}} \diff q\,\diff p \\
& \quad -\frac2t \int_{\mathcal K} \phi_{\bar \de}(p) qF \frac{(\pa_\theta f_0)\circ T_t \theta_q+ (\pa_qf_0 )\circ T_t }{\sqrt{p}} \diff q\,\diff p \\
& \quad - \int_{\mathcal K} qF(1-\phi_{\bar \de}(p)) \frac{\pa_\theta f_0\circ T_t}{\sqrt{p}} \diff q\,\diff p \notag\\
& \eqqcolon I(t,r) + II(t,r) + III(t,r).
\end{align*}
We focus first on $I(t,r)$. Observing that, due to the cut-off function $\phi_{\bar \de}$, the integrand vanishes on $\pa J_{\text{trap}}$. 
We integrate by parts with respect to $q$ to get
\begin{align}
I(t,r)=&\,\frac2t  \int_{\mathcal K} \phi_{\bar \de} (q\pa_q F+F)\frac{f_0\circ T_t}{\sqrt{p}} \diff q\,\diff p \notag\\
=
& \frac2t\int_{L_0}^{L_{\max}}\int_{\Psi_L(r)}^{E_0} \phi_{\bar \de} (\om\frac{\pa_q F}{F}+1)\frac{f_0\circ T_t}{\sqrt{2(E-\Psi_L(r))}}\,\diff E\,\diff L,\label{E:ILINTER}
\end{align}
where we have used the assumptions $k>1,\ell>0$, and~\eqref{E:INITIALL} to make sure that the boundary terms corresponding to the vacuum boundary $E(p,q)=E_0$ and $L(p,q)=L_0$ vanish. 
To bound the right hand side of~\eqref{E:ILINTER}, we estimate as
\beqa
\bigg|\frac2t&\,\int_{L_0}^{L_{\max}}\int_{\Psi_L(r)}^{E_0} \phi_{\bar \de} (\om\frac{\pa_q F}{F}+1)\frac{f_0\circ T_t}{\sqrt{2(E-\Psi_L(r))}}\,\diff E\,\diff L\bigg|\\
\leq&\, \frac Ct\int_{L_0}^{L_{\max}}\int_{\Psi_L(r)}^{E_0}\frac{(E_0-E)^{k-1}(L-L_0)^{\ell}}{\sqrt{2(E-\Psi_L(r))}}\,\dif E\,\dif L\leq \frac{C}{t},\nonumber
\eeqa 
where we have again used that $k>1$, $\ell>0$. Thus we have obtained
\beq\label{E:ILEST}
\big|I(t,r)\big|\leq \frac{C}{t}.
\eeq
To bound $II(t,r)$, we first note that 
$$\big|\pa_q\th\big|\leq\big|\th_E|\Big|\frac{\pa E}{\pa q}\Big|+|\th_L|\Big|\frac{\pa L}{\pa q}\Big|\leq C(E-\EminL)^{-\frac12}(E-\Psi_L(r))^{-\frac12},$$
while $|\pa_\th f_0|= |\varphi'||\pa_\th g_0|\leq C|\varphi'|(E-\EminL)^{\frac12}$, and hence
\beqa\label{E:IILEST}
\big|II(t,r)\big|\leq&\,\frac{C}{t}\int_{L_0}^{L_{\max}}\int_{\Psi_L(r)+{\bar \de}}^{E_0}\frac{1}{E-\Psi_L(r)}\,\dif E\leq \frac{C|\log {\bar \de}|}{t}.
\eeqa
Finally, to control $III(t,r)$, we revert to $(E,L)$ coordinates and observe that
\begin{align}
\big| III(t,r)\big| &\le \int_{\{\phi_{\bar \de}<1\}} \Big|\om F(1-\phi_{\bar \de}) \frac{\pa_\theta f_0\circ T_t}{w}\Big|\,\diff (E,L)\leq C\int_{\{\phi_{\bar \de}<1\}}  \frac{1}{\sqrt{E-\Psi_L(r)}}\,\diff(E,L) \notag\\
& \leq C\bar \de^{\frac12},\label{E:IIILEST}
\end{align}
where we have used  that the domain has measure at most $\bar\de$. 
Combining~\eqref{E:ILEST},~\eqref{E:IILEST} and~\eqref{E:IIILEST},  we have obtained
\beq
|\rho_t(t,r)|\leq C\Big(\frac{|\log\bar \de|}{t}+\bar\de^{\frac12}\Big).\nonumber
\eeq
Taking $\bar \de=\frac1{t^2}$, we conclude the claimed estimate.
%\ec
%%%%%%%%%%%%%%%%%%%%%%%%%%%%%%%%%%%%%%%%%
%%%%%%%%%%%%%%%%%%%%%%%%%%%%%%%%%%%%%%%%%

\end{document}